\documentclass[10pt]{article}
\usepackage{latexsym,amsmath,amssymb,graphics,amscd}
\makeatletter
\@addtoreset{figure}{section}
\def\thefigure{\thesection.\@arabic\c@figure}
\def\fps@figure{h,t}
\@addtoreset{table}{bsection}

\def\thetable{\thesection.\@arabic\c@table}
\def\fps@table{h, t}
\@addtoreset{equation}{section}

\makeatother

\begin{document}

\newtheorem{theorem}{Theorem}[section]
\newtheorem{definition}[theorem]{Definition}
\newtheorem{lemma}[theorem]{Lemma}
\newtheorem{remark}[theorem]{Remark}
\newtheorem{proposition}[theorem]{Proposition}
\newtheorem{corollary}[theorem]{Corollary}
\newtheorem{example}[theorem]{Example}
\newtheorem{examples}[theorem]{Examples}

\newcommand{\bfi}{\bfseries\itshape}

\newsavebox{\savepar}
\newenvironment{boxit}{\begin{lrbox}{\savepar}
\begin{minipage}[b]{15.8cm}}{\end{minipage}
\end{lrbox}\fbox{\usebox{\savepar}}}

\makeatletter
\title{{\bf The reduced spaces of a symplectic Lie group action}}
\author{Juan-Pablo Ortega$^1$ and Tudor S. Ratiu$^2$}
\addtocounter{footnote}{1}
\footnotetext{Centre National de la Recherche Scientifique, D\'epartement de
Math\'ematiques de Besan\c con, Universit\'e de Franche-Comt\'e, UFR des Sciences
et Techniques. 16, route de Gray. F-25030 Besan\c con cedex. France. {\texttt
Juan-Pablo.Ortega@math.univ-fcomte.fr }}
\addtocounter{footnote}{1}
\footnotetext{Section de Math\'ematiques and Centre Bernoulli. \'Ecole
Polytechnique F\'ed\'erale de Lausanne. CH-1015 Lausanne.  Switzerland.
{\texttt Tudor.Ratiu@epfl.ch}}

\date{}
\makeatother
\maketitle

\addcontentsline{toc}{section}{Abstract}

\begin{abstract}  There exist three main approaches to reduction
associated to canonical Lie group actions on a symplectic manifold,
namely, foliation reduction, introduced by Cartan, Marsden-Weinstein
reduction, and optimal reduction, introduced by the authors. When the
action is free, proper, and admits a momentum map these three approaches
coincide. The goal of this paper is to study the general case of a
symplectic action that does not admit a momentum map and one needs to
use its natural generalization, a cylinder valued momentum map
introduced by Condevaux, Dazord, and Molino~\cite{condevaux dazord and
molino}. In this case it will be shown that the three reduced spaces
mentioned above do not coincide, in general. More specifically, the
Marsden-Weinstein reduced spaces are not symplectic but Poisson and
their symplectic leaves are given by the optimal reduced spaces.
Foliation reduction produces a symplectic reduced space whose Poisson
quotient by a certain Lie group associated to the group of symmetries of
the problem equals the Marsden-Weinstein reduced space. We illustrate
these constructions with concrete examples, special emphasis being given
to the reduction of a magnetic cotangent bundle of a Lie group in the 
situation when  the magnetic term ensures  the non-existence  of the
momentum map for the lifted action. The precise relation of the cylinder
valued momentum map with group valued momentum maps for Abelian Lie
groups is also given.
\end{abstract}

\tableofcontents

\section{Introduction}

Let $(M, \omega )$ be a connected paracompact symplectic manifold acted 
upon properly and canonically by a Lie group $G$. In this paper it is 
assumed that the $G$-action is free; the non-free case is the subject of
\cite{stratified spaces}. Let $\mathfrak{g}$ be  the Lie algebra of $G$ 
and $\mathfrak{g}^\ast$ its dual. Assume for the moment that the 
action admits a standard  equivariant momentum map
$\mathbf{J} :M \rightarrow \mathfrak{g}^\ast$. There are three main
approaches to the symmetry reduction of $(M,\omega)$ by
$G$ that yield, up to connected components, the same spaces:
\begin{itemize}
\item {\bf Foliation reduction}~\cite{cartan}{\bf :} consider the fiber
$\mathbf{J}^{-1}(\mu) $ and the characteristic distribution $D=\ker T
\mathbf{J}\cap  (\ker T \mathbf{J}) ^\omega$ on it; the upper index
$\omega$ on a vector subbundle of  $TM$ denotes the $\omega$-orthogonal
complement. The symplectic structure of $(M,\omega )$ drops naturally to
the leaf space $\mathbf{J}^{-1}(\mu)/ D$.
\item {\bf Marsden-Weinstein reduction}~\cite{mwr}{\bf :} let $G _\mu$
be the isotropy subgroup of the element $\mu \in  \mathfrak{g}^\ast$
with respect to the coadjoint action of  $G $ on $\mathfrak{g}^\ast$.
The orbit manifold $\mathbf{J}^{-1}(\mu)/G_{\mu}$ inherits from 
$(M, \omega)$ a  natural symplectic form $\omega _\mu$ uniquely
characterized by  the expression  $i^\ast_\mu\omega =
\pi_\mu^\ast\omega_\mu$, with $i_\mu: \mathbf{J}^{-1}(\mu)
\hookrightarrow M$ the inclusion and $\pi_\mu : \mathbf{J}^{-1}(\mu)
\rightarrow \mathbf{J}^{-1}(\mu)/G _\mu$ the projection.
\item {\bf Optimal reduction}~\cite{optimal, symplectic reduced}{\bf :}  
let $A_{G}'$ be the distribution on $M$ defined by $ A_{G}':=\{ X_f \mid
f \in C^\infty(M)^G\}$. The distribution $A_{G}'$ is smooth and
integrable in the sense of Stefan and Sussmann~\cite{stefan, stefan b,
sussman}. The  optimal momentum map ${\cal J}:M \longrightarrow M /
A_{G}'$ is defined as the canonical projection onto the leaf space of
$A_{G}'$ which is, in most cases, not even a Hausdorff topological
space, let alone a manifold. For any $g\in G$, the map
$\Psi_g(\rho)={\cal J}(g\cdot m)\in M/A_{G}'$ defines a continuous
$G$-action on $M/A_{G}'$ with respect to which ${\cal J}$ is
$G$-equivariant. 
The orbit space $M _\rho:={\cal J}^{-1}(\rho)/
G _\rho$ is a smooth symplectic regular quotient manifold with symplectic form
$\omega_\rho$ characterized by $
\pi_\rho^\ast\omega_\rho=i _\rho^\ast \omega $, where $\pi _\rho:
\mathcal{J}^{-1}(\rho) \rightarrow 
\mathcal{J}^{-1}(\rho)/G_{\rho} $  is the projection and $i _\rho:
\mathcal{J}^{-1}(\rho) \hookrightarrow M $ the inclusion.
\end{itemize}
These reduction theorems are important for symmetric Hamiltonian
dynamics since the flow associated to  a $G$-invariant Hamiltonian
function projects to a Hamiltonian flow on the symplectic reduced spaces.

\medskip

Our goal in this paper is to carry out the regular reduction procedure
for \textit{any\/} symplectic action, even when a momentum map does not
exist. As will be shown, the three approaches to reduction yield spaces
that are, in general, distinct but that are non-trivially related to each
other in very interesting ways. Our results are based on a key
construction of Condevaux, Dazord, and Molino~\cite{condevaux dazord and
molino} naturally generalizing the standard momentum map to a
{\bfi  cylinder valued  momentum map} $\mathbf{K}:M \rightarrow
\mathbb{R}^a \times \mathbb{T}^b $, $a,b \in  \Bbb N $, that
\textit{always\/} exists for \textit{any\/} symplectic Lie group action.
The cylinder $\mathbb{R}^a \times \mathbb{T}^b$ is obtained as the
quotient $\mathfrak{g}^\ast /\overline{{\mathcal H}}  $, with ${\mathcal
H} $ a discrete subgroup of $(\mathfrak{g}^\ast, +)$ which is the
holonomy of  a flat connection on the trivial principal fiber bundle 
$\pi: M \times \mathfrak{g}^\ast \rightarrow M$ with
$(\mathfrak{g}^\ast,+) $ as Abelian structure group. This flat connection
is constructed using exclusively the canonical $G$-action and the
symplectic form $\omega$ on $M$ thereby justifying the name 
{\bfi  Hamiltonian holonomy\/} for ${\mathcal H}$.
\medskip

\noindent{\bf The main result.}  {\it Let $(M, \omega)$ be a
connected paracompact symplectic manifold and $G$ a Lie group acting
freely and properly on it by symplectic diffeomorphisms. Let
$\mathbf{K}:M \rightarrow \mathfrak{g}^\ast/\overline{{\mathcal H}}$ be
a cylinder valued momentum map for this action. Then
$\mathfrak{g}^\ast/\overline{{\mathcal H}}$ carries a natural Poisson
structure and there exists a smooth $G$-action on it with respect to
which $\mathbf{K}$ is equivariant and Poisson. Moreover:
\smallskip

\noindent {\bf (i)} The Marsden-Weinstein reduced space $M^{[\mu]}
:=\mathbf{K} ^{-1}([\mu])/G_{[\mu]}$, $[\mu]\in \mathfrak{g}^\ast/
\overline{{\mathcal H}} $, has a natural Poisson structure inherited
from the symplectic structure $(M, \omega)$ that is, in general,
degenerate. $M ^{[\mu]}$ will be referred to as the {\bfi  Poisson
reduced space}.
\smallskip

\noindent {\bf (ii)} The optimal reduced spaces can be naturally identified with
the symplectic leaves of $M ^{[\mu]}$.
\smallskip

\noindent {\bf  (iii)} The reduced spaces obtained by foliation reduction
equal the orbit spaces $M _{[\mu]}:= \mathbf{K}^{-1}([\mu])/N_{[\mu]}$,
where $N$ is a normal connected Lie subgroup of $G$ whose Lie algebra is
the annihilator $\mathfrak{n}:= \left({\rm Lie}\left(\overline{{\cal
H}}\right)\right)^\circ \subset  \mathfrak{g}$ of 
${\rm Lie}\left(\overline{{\cal H}}\right) \subset \mathfrak{g}^\ast$ in
$\mathfrak{g}$. The manifolds $M _{[\mu]}$ will be referred to as the
{\bfi  symplectic reduced spaces}. 
\smallskip

\noindent {\bf (iv)}The quotient group $H _{[\mu]}:= G _{[\mu]}/
N_{[\mu]}$ acts canonically on  $M _{[\mu]} $ and the quotient Poisson
manifold $M _{[\mu]}/ H _{[\mu]} $ is Poisson diffeomorphic to
$M^{[\mu]}$.}

As will be shown in the course of this paper, one of the reasons behind
the existence of the three distinct reduced manifolds is the
non-closedness of the discrete Hamiltonian holonomy ${\cal H}$ (as the
holonomy group of a flat connection). In fact, $\overline{{\cal H}}$
measures in some sense the  degree of degeneracy of the Poisson
structure of the Marsden-Weinstein reduced space $M  ^{[\mu]}$. Moreover,
when ${\mathcal H}$ is closed, the three reduction approaches yield
(up to connected components) the same symplectic space.

The present paper deals only with free actions. In our forthcoming
paper~\cite{stratified spaces} we will study the situation in which this
hypothesis has been dropped.

The contents of the paper are as follows. Section~\ref{The cylinder
valued momentum map section} introduces and presents in detail the
properties of the cylinder valued momentum map. Section~\ref{The
equivariance properties of the cylinder valued  map section} studies the
invariance properties of the Hamiltonian holonomy $ {\mathcal H} $ and
constructs a natural action on the target space of the
cylinder valued momentum map with respect to which the cylinder valued
momentum map is equivariant. This action is an essential
ingredient for reduction. Section~\ref{Poisson structures on and
extensions} defines a Poisson structure on the target space of the
cylinder valued momentum map  with respect to which this map is 
Poisson. It also provides a careful study of this Poisson structure and
explicitly characterizes its symplectic leaves. This section
also contains a general discussion on central extensions of Lie algebras
and groups, their actions, and their role in the characterization of the
symplectic leaves of affine and projected affine Lie-Poisson structures
on duals of Lie algebras. Apart from its intrinsic interest, this
information on central extensions will be heavily used in the example
of Section~\ref{Example: Magnetic cotangent bundles of Lie groups
section}. Section~\ref{The reduction theorems section} contains a
detailed statement and proof of the reduction results announced
above. Section~\ref{Example: Magnetic cotangent bundles of Lie groups
section}  contains an in-depth study of an example that illustrates some
of the main results in the paper. The cotangent bundle of a Lie group is
considered, but with a symplectic structure that is the sum of the
canonical one and of an invariant magnetic term, whose value at the
identity does not integrate to a group two-cocycle. This modification
destroys, in general, the existence of a standard momentum map for the
lift of left translations and forces the use of all the developments in
the paper. This section contains an interesting generalization of the
classical result that states that the coadjoint orbits endowed
with their canonical Kostant-Kirillov-Souriau symplectic structure are 
symplectic reduced spaces of the cotangent bundle of the corresponding
Lie group. The paper concludes with an appendix that specifies the
relation, in the context of Abelian Lie group actions, of the cylinder
valued momentum map and the so called Lie group valued momentum maps.
 
\medskip 

\noindent {\bf Notations and general assumptions.} {\bf Manifolds:} In
this paper all manifolds are finite dimensional.
{\bf Group actions:} The image of a point $m$ in a manifold $M$ under a
group action $\Phi:G \times  M \rightarrow M $ is denoted interchangeably
by $\Phi(g,m)= \Phi_g (m)= g \cdot  m$, for any $g \in  G$. The symbol
$L _g: G \rightarrow G $ (respectively $R _g:G \rightarrow  G $) denotes
left (respectively right) translation on $G$ by the group element $g 
\in G $. The group orbit containing $m \in  M $ is denoted by
$G\cdot m$  and its tangent space by $T _m( G \cdot m )$ or
$\mathfrak{g}\cdot  m$.  The Lie algebra of the group $G$ is usually
denoted by $\mathfrak{g}$. Given any $\xi \in 
\mathfrak{g}$, the symbol $\xi _M $ denotes the infinitesimal generator
vector field associated to $\xi$ defined by $\xi_M (m)=
\left.\frac{d}{dt}\right|_{t=0}\exp \, t \xi \cdot  m
$, for any $m \in M  $. A {\bfi right (left) Lie algebra action} of $\mathfrak{g}$
on $M$ is a Lie algebra (anti)homomorphism $\xi \in
\mathfrak{g}\longmapsto \xi_M \in \mathfrak{X}(M)$ such that the mapping $(m,
\xi)\in M
\times \mathfrak{g}\longmapsto \xi_M (m) \in TM$ is smooth. If $\mathfrak{g}$ acts
on a symplectic manifold $(M , \omega)$ we say that the $\mathfrak{g}$-action is
{\bfi   canonical} when $\pounds _{\xi_M}\omega=0 $, for any $\xi \in
\mathfrak{g}$.  {\bf The Chu map:} Given a symplectic manifold
$(M,
\omega)$ acted canonically upon by a Lie algebra
$\mathfrak{g}$,  the Chu map
$\Psi: M
\rightarrow Z ^2 (\mathfrak{g})$ is defined by the expression $\Psi (m)(\xi,
\eta):= \omega (m) (\xi_M(m), \eta_M(m))$, for any $m \in  M $, $\xi, \eta, \in
\mathfrak{g}$.

\section{The cylinder valued momentum map}
\label{The cylinder valued momentum map section}

In this section we define carefully the cylinder valued momentum map and 
study its elementary properties. This construction, first
introduced by Condevaux, Dazord, and Molino  in~\cite{condevaux dazord
and molino} under the name of ``reduced momentum map", is the key stone
of the main results in this paper.

The following notations will be used throughout  this work. If
$\langle \cdot, \cdot \rangle : W^\ast \times W \rightarrow \mathbb{R}$
is a nondegenerate duality pairing and $V \subset W$, define the {\bfi
annihilator\/} subspace $V^\circ : = \{ \alpha\in W ^\ast \mid \langle
\alpha, v \rangle = 0 \; \text{for all}\; v \in V \} \subset W ^\ast$ and
similarly for a subset of $W^\ast$. If $(S, \omega)$ is a symplectic
vector space and $U \subset S $, define the $\omega$-{\bfi orthogonal
subspace\/} $U^\omega : = \{s \in S \mid \omega(s, u) = 0 \; \text{for
all} \; u \in U \}$.

Let $(M, \omega )$ be a connected and paracompact symplectic manifold
and let $\mathfrak{g}$ be a Lie algebra that acts canonically on $M$.
Take the Cartesian product $M \times 
\mathfrak{g}^\ast$  and let
$\pi:M\times\mathfrak{g}^\ast\rightarrow M$ be the projection onto
$M$. Consider
$\pi$ as the bundle map of the trivial principal fiber bundle $(M \times
\mathfrak{g}^\ast, M,
\pi, \mathfrak{g}^\ast)$ that has $(\mathfrak{g}^\ast,+) $ as Abelian structure
group. The group
$(\mathfrak{g}^\ast,+) $ acts on $M \times \mathfrak{g}^\ast $ by $\nu \cdot (m,
\mu):=(m, \mu- \nu)$, with $m \in M $ and $\mu, \nu \in \mathfrak{g}^\ast$. Let
$\alpha
\in \Omega^1(M \times \mathfrak{g}^\ast; \mathfrak{g}^\ast)$ be the connection
one-form defined by
\begin{equation}
\label{definition of alpha connection}
\langle \alpha(m , \mu) (v _m, \nu), \xi\rangle:=
(\mathbf{i}_{\xi_M} \omega) (m)
(v _m) -\langle \nu, \xi \rangle,
\end{equation} where $(m, \mu)\in M \times \mathfrak{g}^\ast $, $(v _m, \nu) \in T
_m M \times
\mathfrak{g}^\ast $,  $\langle\cdot , \cdot \rangle $ denotes the natural pairing
between
$\mathfrak{g}^\ast $ and $\mathfrak{g}$, and $\xi_M $ is the infinitesimal
generator vector field associated to $\xi\in \mathfrak{g}$ defined by 
$
\xi_M(m)= \left.\frac{d}{dt}\right|_{t=0}\exp \,t \xi \cdot  m$, $m \in M $. It is
easy to check that
$\alpha$  is a flat connection.  For
$(z,
\mu)
\in M
\times
\mathfrak{g}^\ast$, let $(M \times \mathfrak{g}^\ast)(z, \mu) $ be the holonomy
bundle through
$(z, \mu) $ and let $ {\mathcal H}(z , \mu)$ be the holonomy group of $\alpha$
with reference point $(z, \mu) $ (which is an Abelian discrete subgroup of
$\mathfrak{g}^\ast$ by the flatness of $\alpha$). The Reduction
Theorem~\cite[Theorem 7.1, page 83]{kobayashi nomizu 1} guarantees that the
principal bundle 
$((M
\times
\mathfrak{g}^\ast)(z, \mu),M, \pi|_{(M \times \mathfrak{g}^\ast)(z,
\mu)},{\mathcal H}(z , \mu))
$ is a reduction of the principal bundle $(M \times \mathfrak{g}^\ast, M, \pi,
\mathfrak{g}^\ast)$; it is here that we used the paracompactness  of 
$M$ since it is  a technical hypothesis in the  Reduction Theorem. To
simplify notation, we will write $(\widetilde{M}, M, \widetilde{p},
{\mathcal H}) $ instead of $((M \times \mathfrak{g}^\ast)(z, \mu),M,
\pi|_{(M \times \mathfrak{g}^\ast)(z, \mu)},{\mathcal H}(z , \mu)) $. Let
$\widetilde{\mathbf{K}}: \widetilde{M} \subset M \times
\mathfrak{g}^\ast\rightarrow \mathfrak{g}^\ast$ be the projection into the
$\mathfrak{g}^\ast$-factor. 

Let $\overline{{\mathcal H}}$ be the closure of ${\mathcal H} $ in
$\mathfrak{g}^\ast$. Since $\overline{{\mathcal H}}$ is a closed subgroup of
$(\mathfrak{g}^\ast, +)$, the quotient $C:= \mathfrak{g}^\ast/ \overline{{\mathcal
H}}$ is a cylinder (that is, it is isomorphic to the Abelian Lie group
$\mathbb{R}^a \times
\mathbb{T}^b$ for some $ a, b \in \mathbb{N}$). Let
$\pi_C:
\mathfrak{g}^\ast\rightarrow
\mathfrak{g}^\ast/\overline{{\mathcal H}}=C$ be the projection. Define
$\mathbf{K}: M \rightarrow  C $ to be the map that makes the following diagram
commutative:
\begin{equation}
\label{diagram commutative cylinder valued momentum map}
\begin{CD}
\widetilde{M}@>\widetilde{\mathbf{K}}>>\mathfrak{g}^\ast\\ @V\widetilde{p} VV	
@VV\pi_C V\\ M@>\mathbf{K}>>\mathfrak{g}^\ast/ \overline{ {\mathcal H}}.
\end{CD}
\end{equation} In other words, ${\bf K}$ is defined by  $ \mathbf{K}(m)= \pi_C
(\nu) $, where
$\nu \in
\mathfrak{g}^\ast $ is any element  such that $(m, \nu) \in \widetilde{M } $. This
is a good definition because if we have two points $(m, \nu), (m, \nu') \in
\widetilde{M } $, this implies that $(m, \nu), (m, \nu') \in \widetilde{p}
^{-1}(m) $ and, as
${\mathcal H} $  is the structure group of the principal fiber bundle
$\widetilde{p}: \widetilde{M} \rightarrow M $, there exists an element $\rho \in
{\mathcal H} $ such that $\nu'= \nu+ \rho$. Consequently, $\pi_C (\nu)= \pi_C(\nu')
$. 

We will refer to $ \mathbf{K}:M \rightarrow  \mathfrak{g}^\ast/ \overline{{\mathcal
H}}=:C
$ as a {\bfi  cylinder valued momentum map} associated to the canonical
$\mathfrak{g}$-action on $(M, \omega)$.
The cylinder valued  momentum  map is a strict generalization of the standard
(Kostant-Souriau) momentum map since it is easy to prove (see for
instance~\cite[Proposition 5.2.10]{hsr}) that the $G$-action has a standard
momentum map if and only if the holonomy group ${\mathcal H} $ is trivial. In such
a case the cylinder valued  momentum map {\it is } a standard momentum
map.

Notice that we refer to {\it ``a''} and not to {\it ``the''}  cylinder valued 
momentum map since each choice of the holonomy bundle of the
connection~(\ref{definition of alpha connection}) defines such a map. In
order to see  how the definition of ${\bf K}$ depends on the choice of
the holonomy bundle $\widetilde{M } $ take $\widetilde{M _1}$ and
$\widetilde{M _2} $  two holonomy bundles of $(M \times
\mathfrak{g}^\ast, M, \pi, \mathfrak{g}^\ast)$. We now notice three
things. First, there exists $\tau \in \mathfrak{g}^\ast$ such that
$\widetilde{M _2}= R _\tau (\widetilde{M _1})$, where
$R _\tau (m, \mu):=(m, \mu+ \tau)$, for any $(m, \mu)\in M \times
\mathfrak{g}^\ast$. Second, since $(\mathfrak{g}^\ast,+)$ is Abelian all the
holonomy groups based at any point are the same and hence the  projection $\pi_C:
\mathfrak{g}^\ast \rightarrow
\mathfrak{g}^\ast/ \overline{{\mathcal H}} $ in~(\ref{diagram commutative cylinder
valued momentum map})  does not depend on the choice of $\widetilde{M } $; in view
of this remark we will refer to ${\cal H}$ as the {\bfi  Hamiltonian holonomy} of
the $G$-action on $(M , \omega  )$. Third,
$\pi_C$ is a group homomorphism. Let now
$\widetilde{p}_{\widetilde{M} _i}: \widetilde{M} _i\rightarrow M $,
$\widetilde{\mathbf{K}}_{\widetilde{M} _i}: \widetilde{M} _i \rightarrow
\mathfrak{g}^\ast$, and $\mathbf{K} _{\widetilde{M} _i}: M \rightarrow
\mathfrak{g}^\ast$ be the maps in the diagram~(\ref{diagram commutative cylinder
valued momentum map}) constructed using the holonomy bundles $\widetilde{M } _i $,
$i \in \{1,2\} $. Let $m \in M $. By definition $ \mathbf{K}_{\widetilde{M} _2} 
(m)= \mathbf{K}_{\widetilde{M} _2}(\widetilde{p} _{\widetilde{M} _2}(m, \nu))$,
where $(m, \nu) \in \widetilde{M} _2$. Since $\widetilde{M _2}= R _\tau
(\widetilde{M _1})$ there exists an element $\nu' \in
\mathfrak{g}^\ast$ such that $(m, \nu') \in \widetilde{M} _1$ and $(m , \nu)= (m ,
\nu' +
\tau)$. Hence, 
\begin{align*}
\mathbf{K}_{\widetilde{M} _2}  (m)&= \mathbf{K}_{\widetilde{M} _2}(\widetilde{p}
_{\widetilde{M} _2}(m, \nu))
	=\mathbf{K}_{\widetilde{M} _2}(\widetilde{p} _{\widetilde{M} _2}(m, \nu' +\tau))\\
	&=
\pi_C(\widetilde{\mathbf{K}}_{\widetilde{M} _2}(m, \nu'+ \tau)) = \pi_C(\nu' +
\tau)=\pi_C(\nu') + \pi_C( \tau)= \mathbf{K}_{\widetilde{M} _1} (m)+ \pi_C (\tau).
\end{align*} 
Since in the previous chain of equalities the point $m \in M
$ is arbitrary and $
\tau \in \mathfrak{g}^\ast$ depends only on $\widetilde{M} _1 $ and $\widetilde{M
} _2 $ we have that
\[
\mathbf{K}_{\widetilde{M} _2}=\mathbf{K}_{\widetilde{M} _1}+ \pi_C(\tau).
\] 
The following proposition summarizes the elementary properties of the
cylinder valued  momentum map. 

\begin{proposition}
\label{properties of the cylinder valued momentum map in list} Let $(M, \omega)$
be a connected and paracompact symplectic manifold and $\mathfrak{g}$
 a Lie algebra acting canonically on it. Then any  cylinder
valued momentum map $ \mathbf{K}: M \rightarrow C $ associated to this action has
the following properties:
\begin{description}
\item [(i)] ${\bf K}$ is a smooth  map that satisfies Noether's Theorem, that is,
for any $\mathfrak{g}
$-invariant function $h \in  C^\infty(M)^{\mathfrak{g}}:=\{f \in  C^\infty(M)\mid
\mathbf{d}h (\xi _M) = 0\text{ for all }\xi\in \mathfrak{g}\}
$, the flow $F _t 
$ of its associated Hamiltonian vector field $X _h$ satisfies the
identity
\[
\mathbf{K} \circ F _t= \mathbf{K}| _{{\rm Dom}(F _t)}.
\]
\item [(ii)] For any $ v _m \in T _mM $, $m \in M $, we have the
relation 
\[ T _m \mathbf{K} ( v _m) = T _\mu \pi_C ( \nu),
\] where $\mu \in \mathfrak{g}^\ast  $ is any element such that $\mathbf{K} (m)=
\pi_C (\mu)
$ and $\nu \in \mathfrak{g}^\ast $ is uniquely determined by:
\begin{equation}
\label{equation that determines nu}
\langle \nu, \xi\rangle=(\mathbf{i}_{\xi_M} \omega)(m) (v _m),
\quad\text{ for any }\quad
\xi\in \mathfrak{g}.
\end{equation}
\item [(iii)] 
$\ker (T _m \mathbf{K})= \left( \left({\rm Lie}(\overline{{\mathcal H}})\right)
^{\circ}\cdot m\right) ^\omega$.
\item [(iv)]{\bf Bifurcation Lemma:}
\[{\rm range}\, (T _m \mathbf{K})= T _\mu \pi_C \left((\mathfrak{g}_{m})^\circ
\right),\]  where
$\mu \in \mathfrak{g}^\ast  $ is any element such that 
$\mathbf{K} (m)= \pi_C (\mu)$.
\end{description}
\end{proposition}

\begin{remark}
\normalfont Later on in Theorem~\ref{cylinder valued momentum map symplectic
dynamics reduction} we will show that the cylinder valued momentum map 
remains constant along the flow of functions that are  less invariant than those
in part {\bf (i)} of the previous proposition.
\end{remark}

\noindent\textbf{Proof.\ \ } Since $\mathfrak{g}^\ast/
\overline{{\mathcal H}}$ is a homogeneous manifold, the canonical
projection $\pi_C: \mathfrak{g}^\ast \rightarrow \mathfrak{g}^\ast/
\overline{{\mathcal H}}$ is a surjective submersion.  Moreover, 
by \eqref{diagram commutative cylinder valued momentum map},
$\mathbf{K} \circ \widetilde{p}= \pi_C \circ \widetilde{\mathbf{K}}$ is
a smooth map. Thus, since $\widetilde{p} $ is a surjective submersion,
it follows that the map $\mathbf{K} $ is necessarily smooth.

We start by proving {\bf (ii)}. Let $m \in  M  $ and $(m, \mu) \in \widetilde{p}
^{-1}(m)$. If $v _m= T _{(m, \mu)} \widetilde{p} (v _m, \nu) $ then \eqref{diagram
commutative cylinder valued momentum map} gives
\[ T _m \mathbf{K} (v _m) =T _m \mathbf{K}\left(T _{(m, \mu)} \widetilde{p} (v _m,
\nu)\right) =T _\mu \pi_C \left(T _{(m, \mu)} \widetilde{ \mathbf{K}} (v _m,\nu)
\right) = T _\mu \pi_C ( \nu).
\]

\medskip

\noindent {\bf (i)} We now check that ${\bf K}$ satisfies Noether's condition. Let
$h \in C^\infty(M)^{\mathfrak{g}} $ and let $F _t$ be the flow of the associated
Hamiltonian vector field $X _h $. Using the expression for the derivative $T
_m\mathbf{K} $ in {\bf (ii)} it follows that $T _m \mathbf{K}( X _h(m)) = T _\mu
\pi_C(\nu)$, where
$\mu \in \mathfrak{g}^\ast  $ is any element such that $\mathbf{K} (m)= \pi_C (\mu)
$ and $\nu \in \mathfrak{g}^\ast $ is uniquely determined by
\[
\langle \nu, \xi\rangle=(\mathbf{i}_{\xi_M} \omega)(m) (X _h (m)) =-\mathbf{d}h
(m)(\xi_M (m)) = \xi_M[h](m)=0,
\] for all $\xi\in \mathfrak{g} $, which proves that $ \nu=0 $ and, consequently,
$T _m
\mathbf{K} (X _h(m))=0 $, for all $m \in M $. Finally, as 
\[
\frac{d}{dt} (\mathbf{K} \circ F _t) (m)=T _{F _t (m)} \mathbf{K} \left( X _h(F _t
(m)) \right)=0, 
\] 
we have $\mathbf{K} \circ F _t= \mathbf{K}| _{{\rm Dom}(F _t)}$, as
required.

\medskip

\noindent {\bf (iii)} Due to the expression in {\bf (ii)}, a vector $v _m \in \ker
T _m \mathbf{K} $ if and only if the unique element $\nu\in \mathfrak{g}^\ast$
determined by~(\ref{equation that determines nu})  satisfies $T _\nu \pi_C (\nu)
=0  $, that is, $\nu\in {\rm Lie}(\overline{{\mathcal H}})$. Equivalently, we have
that
$\langle \nu, \xi\rangle=0 $, for any $\xi \in ({\rm Lie}(\overline{{\mathcal
H}}))^\circ \subset (\mathfrak{g}^\ast)^\ast = \mathfrak{g} $ which, in
terms  of $v _m$, yields that $ (\mathbf{i}_{\xi_M} \omega)(m)(v_m)
=0 $ for any $\xi \in ({\rm Lie}(\overline{{\mathcal H}}))^\circ $. This
can obviously be rewritten by saying that
$v _m \in \left( \left({\rm Lie}(\overline{{\mathcal H}})\right) ^{\circ}\cdot
m\right) ^\omega$.

\medskip

\noindent {\bf (iv)} We start by checking that ${\rm range}\, (T _m
\mathbf{K})\subset   T_\mu\pi_C \left((\mathfrak{g}_{m})^\circ\right)$. Let  $T _m
\mathbf{K} ( v _m )\in {\rm range}\, (T _m \mathbf{K}) $. Let $\nu \in
\mathfrak{g}^\ast $ be the element determined by~(\ref{equation that determines
nu}) which hence satisfies $T _m
\mathbf{K}(v_m) = T_\mu \pi_C (\nu) $. Now, notice that for any $\xi \in
\mathfrak{g}_{m}$ we have that 
\[
\langle \nu, \xi\rangle= \omega(m)(\xi_M (m), v _m)=0
\] which implies that $\nu \in (\mathfrak{g}_{m})^{\circ } $. This proves the
inclusion 
${\rm range}\, (T _m \mathbf{K})\subset  T_\mu \pi_C \left((\mathfrak{g}_{m})^\circ
\right)$. Hence, the equality will be proven if we show that 
\begin{equation}
\label{equality of dimensions range} {\rm rank}\,( T _m\mathbf{K})=\dim \left(T_\mu
\pi_C \left((\mathfrak{g}_{m})^\circ\right)\right).
\end{equation} On one hand we can use the equality in {\bf (iii)} to obtain 
\begin{align} {\rm rank}\,( T _m\mathbf{K}) &= \dim M- \dim (\ker T _m\mathbf{K})=
\dim M-\dim M+\dim\left(
\left({\rm Lie}(\overline{{\mathcal H}})\right) ^{\circ}\cdot m\right)\notag\\
	&= \dim\left( \left({\rm Lie}(\overline{{\mathcal H}})\right) ^{\circ}\right)-\dim
\left(\mathfrak{g}_{m}\cap \left({\rm Lie}(\overline{{\mathcal H}})\right)
^{\circ}\right).\label{half of the range equality}
\end{align} On the other hand,
\begin{align*}
\dim (&T_\mu\pi_C \left( (\mathfrak{g}_{m})^\circ\right)) = \dim
(\mathfrak{g}_{m}) ^{\circ}-\dim (\ker T _\mu
\pi_C|_{(\mathfrak{g}_{m})^{\circ}})= \dim \mathfrak{g}- \dim
\mathfrak{g}_{m}-\dim (\ker T _\mu \pi_C\cap (\mathfrak{g}_{m})^{\circ})\notag\\
	&= \dim \mathfrak{g}- \dim \mathfrak{g}_{m}-\dim ({\rm Lie}(\overline{{\mathcal
H}})\cap (\mathfrak{g}_{m})^{\circ})\notag\\
	&= \dim \mathfrak{g}- \dim \mathfrak{g}_{m}-\dim({\rm Lie}(\overline{{\mathcal
H}}))-\dim (\mathfrak{g}_{m})^{\circ}+ \dim ({\rm Lie}(\overline{{\mathcal H}})+
(\mathfrak{g}_{m})^{\circ})\notag\\
	&=-\dim({\rm Lie}(\overline{{\mathcal H}}))+\dim \mathfrak{g}- \dim ([{\rm
Lie}(\overline{{\mathcal H}})+ (\mathfrak{g}_{m})^{\circ}]^\circ)=\dim\left({\rm
Lie}(\overline{{\mathcal H}})) ^{\circ}\right)-\dim
\left(\mathfrak{g}_{m}\cap \left({\rm Lie}(\overline{{\mathcal H}})\right)
^{\circ}\right),\notag
\end{align*} which coincides with~(\ref{half of the range equality}), thereby
establishing~(\ref{equality of dimensions range}).  \quad $\blacksquare$

\begin{proposition}[The cylinder valued momentum map and restricted actions] Let
$(M, \omega )$ be a connected and paracompact symplectic manifold, $\mathfrak{g}$
a Lie algebra  acting symplectically on it, and  $\mathbf{K}_{\mathfrak{g}}:M
\rightarrow \mathfrak{g}^\ast/\overline{{\cal H}}_{\mathfrak{g}}$ an
associated cylinder valued momentum map. Let $\mathfrak{h}$ be
a Lie subalgebra of $\mathfrak{g}$ and let
$ {\cal H}_{\mathfrak{h}} $ be the Hamiltonian holonomy of the
$\mathfrak{h}$-action. Then
\begin{description} 
\item  [(i)] $i ^\ast ({\mathcal H}_{\mathfrak{g}})\subset {\mathcal
H}_{\mathfrak{h}} $, where
$i ^\ast : \mathfrak{g}^\ast \rightarrow \mathfrak{h} ^{\ast} $ is the dual of the
inclusion $i:
\mathfrak{h}\hookrightarrow \mathfrak{g}$. Hence there is a unique Lie group
epimorphism
$\overline{i ^\ast }: \mathfrak{g}^\ast/ \overline{{\cal H}}_{\mathfrak{g}}
\rightarrow 
\mathfrak{h}^\ast / \overline{{\cal H}} _{\mathfrak{h}} $ such that $\overline{i
^\ast } \circ
\pi _{\mathfrak{g}^\ast }= \pi_{\mathfrak{h}^\ast } \circ  i ^\ast  $, with
$\pi_{\mathfrak{h}^\ast }:
\mathfrak{h} ^{\ast} \rightarrow  \mathfrak{h} ^{\ast} / \overline{{\cal
H}}_{\mathfrak{h}} $ and  $\pi _{\mathfrak{g}^\ast}: \mathfrak{g}^\ast\rightarrow
\mathfrak{g}^\ast/ \overline{{\cal H}}_{\mathfrak{g}} $ the natural projections.
\item [(ii)] Let $\widehat{i ^\ast }:M \times   \mathfrak{g}^\ast \rightarrow  M
\times 
\mathfrak{h} ^{\ast}$ be the map given by $\widehat{i ^\ast }(m, \mu):= (m, i
^\ast (\mu))$,
$(m, \mu) \in  M \times   \mathfrak{g}^\ast$, and  $\widetilde{ M _{\mathfrak{g}}}
$ the holonomy bundle used in the construction of $\mathbf{K}_{\mathfrak{g}} $ and
that contains the point $(m _0, \mu _0)$. Let   $\widetilde{ M _{\mathfrak{h}}} $
be the holonomy bundle for the
$H$-action containing the point $(m, i ^\ast (\mu)) $. Then
\begin{equation}
\label{inclusion is holonomies}
\widehat{i ^\ast }(\widetilde{M _{\mathfrak{g}}})\subset  \widetilde{M
_{\mathfrak{h}}}.
\end{equation}
\item [(iii)] Let  $\mathbf{K}_{\mathfrak{h}}:M \rightarrow  \mathfrak{h} ^{\ast}/
\overline{{\cal H}} _{\mathfrak{h}} $  be the $\mathfrak{h}$-cylinder valued
momentum map constructed using $\widetilde{ M _{\mathfrak{h}}} $. Then
\begin{equation}
\label{relation cylinder valued momentum map restricted}
\mathbf{K} _{\mathfrak{h}}= \overline{i ^\ast } \circ \mathbf{K} _{\mathfrak{g}}.
\end{equation} 
\end{description}
\end{proposition}

\noindent\textbf{Proof.\ \ }{\bf (i)} Let $\mu \in  {\mathcal H}_{\mathfrak{g}}$
and $c (t)
\subset M $ be a loop in $M$  such that  $c (0)=c (1) =m _0 $ whose horizontal
lift $(c (t) ,
\mu (t) ) $ is such that  $\mu (0)= \mu _0 $  and $ \mu (1)- \mu _0= \mu  $. The
horizontality of $(c (t) ,
\mu (t) ) $ means that for any $\xi\in \mathfrak{g} $ the equality $\langle \mu'
(t),
\xi\rangle= \omega (c (t))(\xi _M (c (t)), c' (t))$  holds. Consequently $ i ^\ast
\mu (1)- i ^\ast \mu _0= i ^\ast \mu $ and $\langle i ^\ast \mu' (t),
\eta\rangle= \omega (c (t))(\eta _M (c (t)), c' (t))$, for any $\eta \in
\mathfrak{h}$ which proves that $(c (t),  i ^\ast  \mu (t)) $  is the
$\mathfrak{h}$-horizontal lift of $c (t) $ passing through $(m _0, i ^\ast  \mu
_0)$ and hence that $i ^\ast  \mu \in  {\mathcal H}_{\mathfrak{h}}$. The rest of
the statement is a straightforward verification.

\smallskip

\noindent {\bf  (ii)} Let $(z, \nu)\in  \widetilde{M_{\mathfrak{g}}}$. By
definition there exists a piecewise smooth $\mathfrak{g}$-horizontal curve $(m
(t), \mu (t))
$ such that  $(m (0), \mu (0))=(m _0, \mu _0)$ and $(m (1), \mu (1))=(z, \nu)$. An
argument similar to the one that we just used in the proof of  ${\bf (i)} $ shows
that the $\mathfrak{g}$-horizontality of $(m (t), \mu (t)) $ implies the
$\mathfrak{h}$-horizontality of
$(m (t), i ^\ast  \mu (t))= \widehat{i ^\ast } \left( (m (t) , \mu (t))\right)$ 
and hence $(m (1), i ^\ast \mu (1))= \widehat{ i ^\ast } (z, \nu) \in 
\widetilde{M _{\mathfrak{h}}} $.

\smallskip

\noindent {\bf (iii)} Let  $ m \in M $ arbitrary. For some $\mu \in 
\mathfrak{g}^\ast$ such that $(m, \mu) \in \widetilde{M _{\mathfrak{g}}} $ we have
that 
\[
\left(\overline{i ^\ast } \circ \mathbf{K}_{\mathfrak{g}} \right)(m)
=\left(\overline{i ^\ast }
\circ \pi_{\mathfrak{g}^\ast } \right)(\mu)= \left( \pi_{\mathfrak{h}^\ast } \circ
i ^\ast 
\right) (\mu).
\] On the other hand $\mathbf{K}_{\mathfrak{h}} (m)= \pi _{\mathfrak{h}^\ast
}(\nu) $, for some
$\nu \in  \mathfrak{h} ^{\ast}$ such that  $(m , \nu)\in  \widetilde{M
_{\mathfrak{h}}} $. Since by~(\ref{inclusion is holonomies}) $(m, i ^\ast 
(\mu))\in \widetilde{ M _{\mathfrak{h}}}$ we have that
$\mathbf{K}_{\mathfrak{h}}(m) = \pi _{\mathfrak{h}}(i ^\ast  (\mu) )= \left(
\overline{i ^\ast } \circ \mathbf{K}_{\mathfrak{g}}\right) (m)$, as required. \quad
$\blacksquare$

\section{The equivariance properties of the cylinder valued  map}
\label{The equivariance properties of the cylinder valued  map section}

Suppose now that the $\mathfrak{g}$-Lie algebra action on $(M, \omega)$ 
considered in the previous section is obtained from a canonical action
of the Lie group $G$ on $(M, \omega)$ by taking the  infinitesimal
generators of all elements  in $\mathfrak{g}$. The main goal of this
section is the construction of a
$G$-action on the target space of the cylinder valued momentum map $\mathbf{K}:M
\rightarrow
 \mathfrak{g}^\ast/\overline{{\mathcal H}} $ with respect to which it is
$G$-equivariant. The following paragraphs generalize to the context of the
cylinder valued momentum map the strategy followed by Souriau~\cite{souriau}  for
the standard momentum map. We start with an important fact about the Hamiltonian
holonomy of a symplectic action.

\begin{proposition}
\label{the Hamiltonian holonomy is Ad invariant} Let $(M, \omega )$ be a connected
and paracompact symplectic manifold and $\Phi: G \times  M \rightarrow  
M$ a symplectic Lie group action. Then the Hamiltonian holonomy
${\cal H}$ of the action is invariant under the coadjoint action, that is,
\begin{equation}
\label{the Hamiltonian holonomy is Ad invariant inclusion}
\mbox{\rm Ad}^\ast _{g ^{-1}} ({\mathcal H})\subset {\mathcal H},
\end{equation} for any $g
\in  G $.
\end{proposition}

\noindent\textbf{Proof.\ \ } Let $ \nu \in  {\mathcal H} $ arbitrary. By
definition of the holonomy group there exists a loop $c:[0,1] \rightarrow M $ at a
point $m \in M $, that is, $c (0)= c (1)=m $, whose horizontal lift
$\widetilde{c}(t)=(c (t) , \mu (t))$ satisfies the relations
$\widetilde{c} (0)=(m,
\mu)$ and $\widetilde{c} (1)=(m, \mu+ \nu)$, for some
$\mu
\in 
\mathfrak{g}^{\ast}$. We now show that $\mbox{\rm Ad}^\ast _{g ^{-1}} \nu \in 
{\mathcal H} $, for any $g \in G $. Take the loop $d :[0,1]\rightarrow M  $  at
the point $g \cdot m $ defined by $d(t):= \Phi_g(c (t)) $. We will prove the claim
by showing that the horizontal lift
$\widetilde{  d} $ of  $d$ is given by 
\begin{equation}
\label{horizontal lift of d}
\widetilde{d}(t)= \left( \Phi _g(c (t)), \mbox{\rm Ad}_{g ^{-1}} ^\ast  \mu
(t)\right).
\end{equation} If this is the case  then $\mbox{\rm Ad}^\ast _{g ^{-1}} \nu \in 
{\mathcal H} $ necessarily  since
$\widetilde{d}(0)=(g \cdot m, \mbox{\rm Ad}^\ast _{g ^{-1}} \mu)$ and
$\widetilde{d} (1)=(g
\cdot m, \mbox{\rm Ad}^\ast _{g ^{-1}} \mu+\mbox{\rm Ad}^\ast _{g ^{-1}} \nu)$. In
order to establish~(\ref{horizontal lift of d}) it suffices to check that
$\widetilde{d} (t) $ is horizontal. Notice first that
\[
\widetilde{d}' (t)= \left(T_{c (t)}\Phi _g \left(c'(t)\right),
\mbox{\rm Ad} ^\ast _{ g ^{-1}} \mu' (t)
\right).
\] Then, for any $\xi \in  \mathfrak{g}$, we have
\begin{align*} 
{\rm \mathbf{i}}_{\xi _M } \omega\left(\Phi _g(c  (t))
\right)&\left(T_{c (t)} \Phi _g (c' (t))\right)+ \left\langle
\mbox{\rm Ad} ^\ast  _{g ^{-1}} \mu ' (t),
\xi\right\rangle\\
	&=\omega\left(\Phi _g(c  (t))\right)\left(T_{c (t)}\Phi _g \left(
\left(\mbox{\rm Ad} _{g ^{-1}} \xi \right) _M (c  (t))\right), T _{c (t)}
\Phi _g (c' (t))\right)+
\left\langle \mu ' (t), \mbox{\rm Ad}_{g ^{-1}} \xi\right\rangle
\\
	&= \omega (c (t))\left(\left(\mbox{\rm Ad}_{g ^{-1}}\xi \right)_M (c (t)), c' (t) 
\right)+
\left \langle \mu ' (t), \mbox{\rm Ad}_{g ^{-1}} \xi\right\rangle=0
\end{align*} 
because of the symplectic character of the action and the fact that
the curve $\widetilde{c}(t)
$ is horizontal. \quad $\blacksquare$

\begin{remark}
\label{inclusions with Hamiltonian holonomy lie}
\normalfont The connected component of the identity $\overline{{\mathcal H}} _0  $
of $\overline{{\mathcal H}}  $ is a vector space. Consequently,
\begin{equation*} {\rm Lie} \left( \overline{{\mathcal H}}\right)=
\overline{{\mathcal H}} _0 \subset
\overline{{\mathcal H}}.
\end{equation*} In order to prove this recall that any Abelian connected Lie
group, like $\overline{{\mathcal H}} _0 $, is isomorphic to
$\mathbb{T}^a\times
\mathbb{R}^b$, for some $a,b \in \mathbb{N}$. Since
$\overline{{\mathcal H}} _0$ is a closed Lie subgroup of $(\mathfrak{g}^\ast, + )$
it cannot contain any compact subgroup and hence $a=0 $ necessarily.
\end{remark}

\begin{corollary} In the hypotheses of Proposition~\ref{the Hamiltonian holonomy
is Ad invariant} the following statements hold
\begin{description}
\item [(i)] $\mbox{\rm Ad} ^\ast _{g ^{-1}} \left( \overline{ {\mathcal
H}}\right)\subset
\overline{{\mathcal H}} $, for any $g  \in G $.
\item [(ii)] $\mbox{\rm Ad} ^\ast _{g ^{-1}} \left( {\rm Lie}\left(\overline{
{\mathcal H}}\right)\right)\subset {\rm Lie}\left(\overline{{\mathcal H}}\right)
$, for any $g  \in G $.
\item [(iii)]  There is a unique group action $\mathcal{A}d^\ast : G
\times 
\mathfrak{g}^\ast/\overline{{\mathcal H}} \rightarrow  \mathfrak{g}^\ast/
\overline{{\mathcal H}}
$ such that  for any $g \in G $
\begin{equation}
\label{funny ad definition}
\mathcal{A}d^\ast_{ g ^{-1}} \circ \pi _C= \pi _C \circ \mbox{\rm Ad} ^\ast  _{g
^{-1}}.
\end{equation} The map $\pi_C : \mathfrak{g}^\ast\rightarrow \mathfrak{g}^\ast/
\overline{{\mathcal H}} $ is the projection. We will refer to $\mathcal{A}d^\ast$
as the {\bfi  projected coadjoint action} of  $G$ on $\mathfrak{g}^\ast/
\overline{{\mathcal H}}  $.
\end{description}
\end{corollary}

\noindent\textbf{Proof.\ \ } {\bf (i)} By~(\ref{the Hamiltonian holonomy is Ad
invariant inclusion}) and the continuity of the coadjoint action we have
$\mbox{\rm Ad} ^\ast _{g ^{-1}} 
\left( \overline{ {\mathcal H}}\right)\subset \overline{\mbox{\rm Ad} ^\ast _{g
^{-1}} 
\left(   {\mathcal H}\right)}\subset
\overline{{\mathcal H}} $.

\medskip

\noindent {\bf (ii)} The inclusion~(\ref{the Hamiltonian holonomy is Ad invariant
inclusion}) guarantees that the restricted map $\mbox{\rm Ad}^\ast _{g ^{-1}}:
\overline{{\mathcal H}} \rightarrow \overline{{\mathcal H}} $ is a Lie group
homomorphism so it induces a Lie algebra homomorphism (which is itself) $\mbox{\rm
Ad}^\ast _{g ^{-1}}:{\rm Lie}
\left(\overline{{\mathcal H}}\right) \rightarrow {\rm Lie
}\left(\overline{{\mathcal H}}\right)
$. In particular this implies the statement.

\medskip

\noindent {\bf (iii)} The map $\mathcal{A}d^\ast_{ g ^{-1}} $ given by
$\mathcal{A}d^\ast_{ g ^{-1}}(\mu+ \overline{{\mathcal H}}):= \mbox{\rm Ad}^\ast
_{g ^{-1}}\mu+ \overline{{\mathcal H}}
$ is well defined by part {\bf (i)} and satisfies~(\ref{funny ad definition}). A
straightforward verification shows that the map  $\mathcal{A}d^\ast$ defines an
action. \quad $\blacksquare$

\medskip

As we will see in the following paragraphs, the results that we just proved allow
us to reproduce in the context of the cylinder valued   momentum map the
techniques introduced by Souriau~\cite{souriau} to study the equivariance
properties of the standard momentum map.

\begin{proposition}
\label{cocycles for cylinder valued momentum map} Let $(M, \omega )$ be a 
connected and paracompact symplectic manifold and $\Phi: G \times  M
\rightarrow   M$ a symplectic Lie group action. Let $\mathbf{K}:M
\rightarrow  \mathfrak{g}^\ast/\overline{{\mathcal H}} $ be  a cylinder 
valued  momentum map for this action. Define $\overline{\sigma}: G\times 
M \rightarrow \mathfrak{g}^\ast /\overline{{\mathcal H}}$ by
\begin{equation}
\label{definition of the non-equivariance cocycle}
\overline{\sigma}(g,m):= \mathbf{K}(\Phi_g( m))-\mathcal{A}d^\ast_{g ^{-1}}
\mathbf{K} (m). 
\end{equation} Then:
\begin{description}
\item [(i)] The map $ \overline{\sigma} $ does not depend on the points $m \in  M
$ and hence it defines a map  
$\sigma: G
\rightarrow  \mathfrak{g}^\ast /\overline{{\mathcal H}} $. 
\item [(ii)] The map $\sigma: G \rightarrow  \mathfrak{g}^\ast
/\overline{{\mathcal H}} $ is a group valued one-cocycle, that is, for any $g,h
\in G $, it satisfies the equality
\[
\sigma(gh)= \sigma(g)+ \mathcal{A}d^\ast _{g ^{-1}} \sigma (h).
\]
\item [(iii)] The map 
\begin{equation*}
\begin{array}{cccc}
\Theta:& G \times  \mathfrak{g}^\ast /\overline{{\mathcal H}}&\longrightarrow &
\mathfrak{g}^\ast /\overline{{\mathcal H}}\\
	&(g, \pi_C(\mu))&\longmapsto&
\mathcal{A}d^\ast _{g ^{-1}}(\pi_C(\mu))+ \sigma (g),
\end{array}
\end{equation*} defines a $G$-action on $\mathfrak{g}^\ast /\overline{{\mathcal
H}}$ with respect to which the cylinder valued  momentum map $\mathbf{K}$ is
$G$-equivariant, that is, for any $g  \in G $, $m \in M $,
\[
\mathbf{K}(\Phi_g(m))= \Theta_g(\mathbf{K}(m)).
\]
\item [(iv)] The infinitesimal
generators of the affine $G$-action on $\mathfrak{g}^\ast/\overline{{\mathcal H}}$
are given by the expression
\begin{equation}
\label{infinitesimal generators of the affine}
\xi_{\mathfrak{g}^\ast/ \overline{{\mathcal H}}} (\pi _C (\mu))=-T _\mu \pi _C
\left(\Psi (m)(\xi, \cdot )
\right),
\end{equation} for any $\xi \in  \mathfrak{g}$, $(m, \mu)\in  \widetilde{M} $, and
where $\Psi:M \rightarrow Z ^2(\mathfrak{g})$ is the Chu map defined by 
$\Psi(\xi, \eta):= \omega\left(\xi_M, \eta_M\right)$, for any $\xi,
\eta\in \mathfrak{g}$.
\end{description} 
We will refer to $\sigma: G \rightarrow  \mathfrak{g}^\ast
/\overline{{\mathcal H}} $ as the {\bfi  non-equivariance one-cocycle} of the 
cylinder valued    momentum map
$\mathbf{K}:M \rightarrow \mathfrak{g}^\ast/ \overline{{\mathcal H}} $ and to
$\Theta$ as the {\bfi  affine $G$-action} on $\mathfrak{g}^\ast
/\overline{{\mathcal H}}$ induced by $\sigma$.
\end{proposition}

\noindent\textbf{Proof.\ \ } {\bf  (i)} For any $ g  \in G $ define the map
$\tau_g:M \rightarrow
\mathfrak{g}^\ast/
\overline{{\mathcal H}} $ by $\tau _g (m):= \overline{\sigma} (g,m)$. We will
prove the claim by showing that $\tau _g $ is a constant map. Indeed, for any
point $z \in  M $ and any vector $v _z\in T _z M $ 
\begin{equation}
\label{first approach by derivatives} 
T _z \tau_g (v _z) 
= T_{g \cdot z } \mathbf{K} \left(T _z \Phi_g (v _z)\right) -
T_{\mathbf{K}(z)}
\mathcal{A}d^\ast _{g ^{-1}} \left(T _z \mathbf{K} (v _z)\right).
\end{equation} 
Recall now that by part~{\bf (ii)} of Proposition~\ref{properties
of the cylinder valued momentum map in list}, $T_z \mathbf{K} (v_z) =
T_\mu \pi_C (\nu)$, where $\mu \in \mathfrak{g}^\ast$ is any element
such that $\mathbf{K} (z)= \pi_C (\mu)$ and $\nu \in \mathfrak{g}^\ast$
is uniquely determined by the equality $\langle \nu,
\xi\rangle=(\mathbf{i}_{\xi_M} \omega)(z)(v_z)$,  for any
$\xi\in \mathfrak{g}$. Equivalently, the relation between $v _z $ and $\nu  $ can
be expressed by saying that the pair $(v _z, \nu) \in T_{(z, \mu)}\widetilde{M} $,
where $\widetilde{M} $ is the holonomy bundle of the connection $\alpha$
in~(\ref{definition of alpha connection}) used in the definition of the cylinder
valued   momentum map ${\bf K}$. Let now $\mu' \in  \mathfrak{g}^\ast$ be such
that $(g \cdot  z, \mu') \in \widetilde{M}$. We now show that $(v_z,
\nu) \in T_{(z, \mu)}\widetilde{M}$ implies that $(T_z \Phi_g (v_z),
\mbox{\rm Ad}^\ast _{g ^{-1}}\nu) \in T_{(g \cdot z,
\mu')}\widetilde{M}$ due to the symplectic character of the $G$-action.
Indeed, for any $\xi \in \mathfrak{g}$ we have
\begin{align*}
 {\rm {\bf i}}_{\xi_M}\omega (g \cdot  z)\left(T _z\Phi_g (v_z)\right) 
&= \omega(g \cdot z)\left(T_z\Phi_g \left(\mbox{\rm Ad}_{g ^{-1}}
\xi \right)_M (z) ,T _z\Phi_g(v _z)\right)
=\omega( z)\left( 
\left(\mbox{\rm Ad}_{g ^{-1}} \xi \right)_M (z) , v _z\right)\\
&= \langle \nu, \mbox{\rm Ad}_{g ^{-1}} \xi\rangle =\langle  \mbox{\rm
Ad}^\ast _{g ^{-1}}\nu,
\xi\rangle,
\end{align*} 
which proves that $(T _z \Phi_g (v_z), \mbox{\rm Ad}^\ast _{g ^{-1}}\nu)
\in T_{(g \cdot z, \mu')}\widetilde{M}$ and hence that $T_{g \cdot z }
\mathbf{K} \left(T _z \Phi_g (v_z)\right) = T_{\mu'}\pi_C \left( 
\mbox{\rm Ad}^\ast _{g ^{-1}}\nu \right)$. If we use this fact to
compute the derivatives in the right hand side of~(\ref{first approach by
derivatives})  we obtain  
\begin{align} 
T _z \tau_g(v _z) 
&=  T_{\mu'}\pi_C\left(\mbox{\rm Ad}^\ast _{g ^{-1}}\nu\right) 
- T_{\mathbf{K}(z)} \mathcal{A}d ^\ast_{g ^{-1}}
\left(T _\mu \pi _C(\nu)\right)
=T_{\mu'}\pi_C \left(\mbox{\rm Ad}^\ast _{g^{-1}}\nu\right)
-T_\mu\left( \pi_C \circ \mbox{\rm Ad}^\ast_{g^{-1}} \right)(\nu)
\notag\\
& = T_{\mu'}\pi_C \left(\mbox{\rm Ad}^\ast _{g^{-1}}\nu\right)
-T_{\mu}\pi_C \left(\mbox{\rm Ad}^\ast _{g ^{-1}}\nu\right)=0.
\label{needed for trivial dependence} 
\end{align} 
The last equality follows from the fact that $\pi _C  $ is an
Abelian Lie group homomorphism. Indeed, for any $\rho, \mu \in 
\mathfrak{g}^{\ast}$
\[ 
T _\mu \pi_C (\rho) = \left.\frac{d}{dt}\right|_{t=0}\pi_C(\mu + t
\rho)=
\left.\frac{d}{dt}\right|_{t=0} \left( \pi _C(\mu)+ \pi_C(t \rho) \right)=
\left.\frac{d}{dt}\right|_{t=0} \pi _C(t \rho)= T _0 \pi _C (\rho),
\] 
and hence $T _\mu \pi_C=T _0 \pi_C$, for any $\mu \in \mathfrak{g}^\ast$.
Finally, the expression~(\ref{needed for trivial dependence}) shows 
that  $T \tau_g=0 $, for any $g  \in G$.  As
$M$ is, by hypothesis, connected this guarantees that
$\overline{\sigma}(g,m)=\overline{\sigma}(g,m')$ for any $m,m' \in M $, 
which proves the claim.

\medskip

\noindent {\bf (ii)} Using the definition~(\ref{definition of the non-equivariance
cocycle}) at the point $h \cdot m  $ we obtain
$
\sigma(g)= \mathbf{K}(gh \cdot  m)- \mathcal{A}d^\ast _{ g ^{-1}} \mathbf{K} (h
\cdot  m)$.
If we now use the point $m \in  M $ we can write
$
\sigma(h)= \mathbf{K}( h \cdot  m)- \mathcal{A}d^\ast _{ h ^{-1}} \mathbf{K} (m)$.
Consequently,
\begin{align*}
\sigma(g)+ \mathcal{A}d^\ast _{g ^{-1}} \sigma (h)
& =\mathbf{K}(gh \cdot m)-\mathcal{A}d ^\ast _{g ^{-1}} \mathbf{K}(h
\cdot   m)+\mathcal{A}d ^\ast _{g ^{-1}}
\mathbf{K}(h \cdot   m)-\mathcal{A}d ^\ast _{g ^{-1}}\mathcal{A}d ^\ast _{h ^{-1}}
\mathbf{K}( m)\\ 
&=\mathbf{K}(gh \cdot   m)-\mathcal{A}d ^\ast _{(gh) ^{-1}}
\mathbf{K}( m)=\sigma(gh),
\end{align*} 
as required. 

\medskip

\noindent {\bf (iii)} It is a straightforward consequence of the definition. 

\medskip

\noindent {\bf (iv)} By definition
\begin{align}
\label{preliminary infinitesimal thing}
\xi_{\mathfrak{g}^\ast/ \overline{{\mathcal H}}} (\pi _C (\mu))&=
\left.\frac{d}{dt}\right|_{t=0}
\Theta _{\exp\, t \xi}(\pi_C (\mu)) \nonumber\\
&= \left.\frac{d}{dt}\right|_{t=0}
\left(\pi_C(\mbox{\rm Ad} ^\ast  _{\exp (-t \xi)} \mu)+ \sigma(\exp t
\xi) \right)
= -T _\mu \pi _C (\mbox{\rm ad}^\ast  _\xi \mu) + T _e \sigma (\xi).
\end{align}
Additionally, if $\mathbf{K}(m)= \pi _C (\mu) $ then
\[
T_e \sigma (\xi) = \left.\frac{d}{dt}\right|_{t=0}\left(\mathbf{K}(\exp
t  \xi \cdot   m)- \mathcal{A}d ^\ast _{\exp(-t \xi)}\mathbf{K} (m) 
\right)
= T_m \mathbf{K} \left(\xi _M (m)\right) + T _\mu \pi _C \left(\mbox{\rm
ad} ^\ast  _\xi \mu\right).
\] 
By~(\ref{equation that determines nu}), $T _m \mathbf{K} \left(
\xi_M(m)\right) = T _\mu \pi _C (\nu)$ with $\nu \in 
\mathfrak{g}_{m}^{\ast}$ uniquely determined by the expression $\langle
\nu,
\eta\rangle=(\mathbf{i}_{\eta_M} \omega)(m) \left(\xi_M(m)\right) =
\Psi(m)(\eta,\xi)$,  for any  
$\eta\in \mathfrak{g}$. Hence~(\ref{preliminary infinitesimal
thing}) yields~(\ref{infinitesimal generators of the affine}).
\quad
$\blacksquare$

\section{Poisson structures on $\mathfrak{g}^\ast/ \overline{{\mathcal H}}$}
\label{Poisson structures on and extensions}

In the following theorem we present a Poisson structure on the target space of the
cylinder valued momentum map $\mathbf{K}:M \rightarrow \mathfrak{g}^\ast/
\overline{{\mathcal H}}$ with respect to which this mapping becomes a Poisson map.
We also see that the symplectic leaves of this Poisson structure can be described
as the orbits of the affine action introduced in the previous section with respect
to a subgroup of $G$ whose definition is related to the non-closedness of the
Hamiltonian holonomy ${\cal H}$ as a subspace of $\mathfrak{g}^{\ast}$.

\begin{theorem}
\label{The Poisson structure of target} Let $(M, \omega)$ be a connected 
paracompact symplectic manifold acted symplectically upon by the Lie 
group $G$. Let $\mathbf{K}:M \rightarrow 
\mathfrak{g}^\ast/\overline{{\mathcal H}} $ be a cylinder valued
momentum map for this action with non-equivariance cocycle
$\sigma:M \rightarrow  \mathfrak{g}^\ast/ \overline{{\mathcal H}} $ and
defined using the holonomy bundle $\widetilde{M} \subset M \times 
\mathfrak{g}^\ast$. The bracket
$\{\cdot,\cdot\}_{\mathfrak{g}^\ast/\overline{{\mathcal H}}} :
C^{\infty}(\mathfrak{g}^\ast/\overline{{\mathcal H}})\times
C^{\infty}(\mathfrak{g}^\ast/\overline{{\mathcal H}}) \rightarrow
\mathbb{R} $ defined by 
\begin{equation}
\label{bracket on target space}
\{ f , g\}_{\mathfrak{g}^\ast/\overline{{\mathcal H}}}(\pi_C (\mu)) =
\Psi (m) \left(\frac{\delta(f \circ \pi _C)}{\delta \mu}, \frac{\delta(g
\circ \pi _C)}{\delta \mu} \right),
\end{equation} 
where $f , g \in C^{\infty}(\mathfrak{g}^\ast/\overline{{\mathcal H}}) $,
$(m, \mu)\in  \widetilde{M} $, $\pi _C: \mathfrak{g}^\ast \rightarrow
\mathfrak{g}^\ast/\overline{{\mathcal H}} $ is the projection, and $\Psi: M
\rightarrow Z ^2 (\mathfrak{g})$ is the Chu map, defines a Poisson 
structure on $\mathfrak{g}^\ast/ \overline{{\mathcal H}} $ such that 
\begin{itemize}
\item [{\rm {\bf (i)}}] $\mathbf{K}:M \rightarrow 
\mathfrak{g}^\ast/\overline{{\mathcal H}} $ is a Poisson map.
\item [{\rm {\bf (ii)}}] The annihilator $ \left( {\rm Lie} \left(
\overline{{\mathcal H}}\right)\right)^{\circ} \subset \mathfrak{g}$ of
${\rm Lie} \left(
\overline{{\mathcal H}}\right) $ in $\mathfrak{g}^\ast$ is an ideal in
$\mathfrak{g}$. Let $N \subset G $ be a connected normal Lie subgroup of  $G$
whose Lie algebra is  $\mathfrak{n}: = \left( {\rm Lie} \left(
\overline{{\mathcal H}}\right)\right)^{\circ}$. The symplectic leaves of
$\left(\mathfrak{g}^\ast/\overline{{\mathcal H}},  \{
\cdot ,\cdot\}_{\mathfrak{g}^\ast/\overline{{\mathcal H}}}\right)$ are 
the orbits of the affine $N$-action on
$\mathfrak{g}^\ast/\overline{{\mathcal H}}$ induced by
$\sigma:G \rightarrow \mathfrak{g}^\ast/\overline{{\mathcal H}}$.
\item  [{\rm {\bf (iii)}}] For any $[\mu]:= \pi _C (\mu) \in 
\mathfrak{g}^\ast/ \overline{{\mathcal H}}$, the symplectic form 
$\omega^+_{N \cdot [\mu]} $ on the affine orbit $N \cdot [ \mu]$
induced by the Poisson structure~(\ref{bracket on target space}) is
given by 
\begin{align*}
\omega^+_{N \cdot [\mu]}
([\mu])\left(\xi_{\mathfrak{g}^\ast /\overline{{\cal H}}}([\mu]) ,
\eta_{\mathfrak{g}^\ast / \overline{{\cal H}}}([\mu])\right) 
&= \omega^+_{N \cdot [ \mu]} ([\mu])\left(-T_{\mu}\pi _C \left(\Psi
(m)(\xi, \cdot)\right), -T_{\mu}\pi _C \left(\Psi (m)(\eta, \cdot)
\right)\right) \\
&= \Psi (m) (\xi, \eta),
\end{align*} 
for any  $\xi, \eta \in  \mathfrak{n} $, $(m, \mu) \in  \widetilde{M}$.
\end{itemize}
\end{theorem}

The proof of this theorem requires several preliminary considerations. 

\begin{lemma}
\label{holonomies are killers} Let ${\cal H}$ be the Hamiltonian holonomy in the
statement of the previous theorem. Then
\begin{equation}
\label{inclusion with bracket bracket}
\overline{{\mathcal H}} \subset [\mathfrak{g}, \mathfrak{g}]^{\circ}.
\end{equation} Moreover
\begin{equation}
\label{inclusion with bracket bracket with Lie} [\mathfrak{g},
\mathfrak{g}]\subset  \left( {\rm Lie} \left( \overline{{\mathcal
H}}\right)\right)^{\circ}
\end{equation} and hence $\left( {\rm Lie} \left( \overline{{\mathcal
H}}\right)\right)^{\circ} $ is an ideal of $\mathfrak{g}$.
\end{lemma}

\noindent\textbf{Proof.\ \ } Let $\mu \in  {\mathcal H} $ be arbitrary.
By definition, there exists a loop $c (t)  $ in $M$ such that $c (0)=
c(1)=m$ and a horizontal lift $\gamma (t):= (c (t), \mu (t)) \in 
\widetilde{M} \subset M \times \mathfrak{g}^\ast$ such that 
$\mu(1)- \mu(0)= \mu $. Since  $\gamma (t)$ is horizontal we have 
${\bf  i}_{\xi_M} \omega \left(c' (t)\right) = \langle  \mu' (t),
\xi\rangle$, for any $\xi \in  \mathfrak{g}$. Since the $G$-action is
symplectic, the infinitesimal generator vector fields $\xi_M$  and
$\eta_M$  are locally Hamiltonian, for any $\xi, \eta \in \mathfrak{g}$,
and hence $[\xi, \eta]_M=-[\xi_M, \eta _M]$ is globally Hamiltonian. Let
$f \in  C^\infty(M)$  be such that $X _f=[\xi, \eta]_M = -[\xi_M,
\eta_M]$. The relation $\mu= \mu(1)- \mu(0)=\int _0 ^1 \mu' (t) d t$
implies thus
\begin{align*}
\langle \mu, &[ \xi, \eta]\rangle 
=\int _0 ^1\langle \mu' (t), [ \xi, \eta]\rangle d t
=\int _0 ^1 \mathbf{i}_{[\xi , \eta]_M}\omega \left(c' (t)\right) d t
=\int _0 ^1\omega (c(t)) ([\xi, \eta]_M(c (t)), c' (t) )d t\\
&=\int _0 ^1 \omega (c(t)) (X _f(c (t)), c' (t) )d t
=\int _0 ^1 \mathbf{d} f (c (t))\left(c' (t)\right) d t= f (c  (1))-f (c
(0))=f (m)- f  (m)=0.
\end{align*} 
This shows that ${\mathcal H} \subset
[\mathfrak{g},\mathfrak{g}]^{\circ} $ and hence that
$\overline{{\mathcal H}} \subset \overline{[\mathfrak{g},
\mathfrak{g}]^{\circ}}=[\mathfrak{g},\mathfrak{g}]^{\circ}$. The
inclusion~(\ref{inclusion with bracket bracket with Lie}) is a
consequence of~(\ref{inclusion with bracket bracket})  and
Remark~\ref{inclusions with Hamiltonian holonomy lie}. \quad
$\blacksquare$ 

\paragraph{Projected Poisson structures.} The Poisson structure in
Theorem~\ref{The Poisson structure of target} will be obtained as the projection
onto the quotient
$\mathfrak{g}^\ast/\overline{{\mathcal H}} $ of an affine Lie-Poisson structure on
$\mathfrak{g}^\ast$. The next proposition proves the existence of this projected
Poisson structure. 

\begin{proposition}
\label{projected lie Poisson structures proposition} Let ${\mathcal H}$ be the
Hamiltonian holonomy in the statement of Theorem~\ref{The Poisson structure of
target}. Let $\Sigma: \mathfrak{g} \times  \mathfrak{g} \rightarrow \mathbb{R}$ be
an arbitrary Lie algebra two-cocycle on $\mathfrak{g}$ and $\{ \cdot , \cdot
\}_{\pm}^\Sigma$ the associated affine Lie-Poisson bracket on $\mathfrak{g}^\ast$
defined by 
\begin{equation}
\label{affine lie Poisson bracke sigma}
\{ f , g \}_{\pm}^\Sigma(\mu):=\pm \left\langle \mu, \left[\frac{\delta f}{ \delta
\mu},\frac{\delta g}{
\delta \mu}\right]\right\rangle\mp \Sigma \left( \frac{\delta f}{ \delta
\mu},\frac{\delta g}{
\delta \mu}\right), \qquad f,g \in C^{\infty}(\mathfrak{g}^\ast).
\end{equation} The action $\phi: \overline{{\mathcal H}} \times  \mathfrak{g}^\ast
\rightarrow \mathfrak{g}^\ast
$ given by $\phi_\nu (\mu):= \mu + \nu $  on  $(\mathfrak{g}^\ast, \{ \cdot , \cdot
\}_{\pm}^\Sigma) $ is free, proper, and Poisson and hence it induces a unique
Poisson structure
$\{ \cdot , \cdot \}_{\mathfrak{g}^\ast/\overline{{\mathcal H}}}^{\pm \Sigma}$ on
$\mathfrak{g}^\ast/\overline{{\mathcal H}} $ such that 
\begin{equation}
\label{projected Poisson structure target}
\{F,H\}^{\pm \Sigma}_{\mathfrak{g}^\ast/\overline{{\mathcal H}}}\circ \pi_C=\{ F
\circ  \pi _C, H\circ
\pi _C \}_{\pm}^\Sigma, \qquad F,H \in
C^{\infty}(\mathfrak{g}^\ast/\overline{{\mathcal H}}).
\end{equation} where $\pi _C: \mathfrak{g}^\ast  \rightarrow
\mathfrak{g}^\ast/\overline{{\mathcal H}} $ is the projection.
\end{proposition}

\noindent\textbf{Proof.\ \ }We first note that if $\mu , \nu \in 
\mathfrak{g}^\ast$ and $f
\in C^{\infty}(\mathfrak{g}^\ast)$ then 
\[
\frac{\delta(f \circ  \phi_\nu)}{\delta \mu}= \frac{\delta f }{ \delta(\mu+ \nu)}.
\] With this equality in mind we can write, for any $f,g \in
C^{\infty}(\mathfrak{g}^\ast) $, $\mu \in \mathfrak{g}^\ast $, and
$\nu\in \overline{{\mathcal H}}$
\begin{align*}
\{ \phi_\nu^\ast f ,\phi_\nu^\ast g \}_{\pm}^\Sigma(\mu)&=
\pm \left\langle \mu, \left[\frac{\delta (f\circ  \phi_\nu)}{ \delta
\mu},\frac{\delta( g\circ  \phi_\nu)}{
\delta \mu}\right]\right\rangle\mp \Sigma 
\left( \frac{\delta (f\circ  \phi_\nu)}{\delta \mu}, 
\frac{\delta(g\circ  \phi_\nu)}{\delta \mu}\right)\\
&=\pm \left\langle \mu, \left[\frac{\delta f}{ \delta (\mu+\nu)},
\frac{\delta g}{\delta (\mu+ \nu)}\right]\right\rangle\mp \Sigma 
\left(\frac{\delta f}{ \delta (\mu+ \nu)},\frac{\delta g}{\delta (\mu+
\nu)}\right)\\
&=\{f ,g \}_{\pm}^\Sigma(\phi_\nu(\mu))- \left\langle \nu,
\left[\frac{\delta f}{
\delta (\mu+ \nu)},\frac{\delta g}{
\delta (\mu+ \nu)}\right]\right\rangle
=\{  f , g \}_{\pm}^\Sigma(\phi_\nu(\mu))
\end{align*} since $\left\langle \nu, \left[\frac{\delta f}{\delta (\mu+
\nu)},\frac{\delta g}{\delta (\mu+ \nu)}\right]\right\rangle =0 $ by
Lemma~\ref{holonomies are killers}. This shows that the action $\phi$ is
Poisson. It is clearly free and proper.
\quad $\blacksquare$ 

\paragraph{Central extensions and their actions.} The description of the
symplectic leaves of the projected Poisson structures
$\left(\mathfrak{g}^\ast/\overline{{\mathcal H}}, \{ \cdot ,
\cdot
\}_{\mathfrak{g}^\ast/\overline{{\mathcal H}}}^{\pm \Sigma}\right)$ requires the
use of certain facts on central extensions that we present in the following
paragraphs.

Let $\mathfrak{g}_\Sigma:= \mathfrak{g}\oplus \mathbb{R}$ be the one-dimensional
central extension of the Lie algebra $\mathfrak{g}$ determined by the cocycle
$\Sigma$ whose bracket is given by the expression
\[
\left[(\xi,s),(\eta,t)\right]:= \left([\xi, \eta], - \Sigma(\xi, \eta) \right).
\] Let $G _\Sigma $ be the  connected simply connected Lie group whose Lie algebra
is
$\mathfrak{g}_\Sigma $. Let $\pi_{\mathfrak{g}}: \mathfrak{g}_\Sigma \rightarrow 
\mathfrak{g}$ be the projection and let $\pi_G: G _\Sigma \rightarrow G $ be the
unique Lie group homomorphism whose derivative is $\pi_{\mathfrak{g}}$ and makes
the diagram
\unitlength=0.4mm
\begin{center}
\begin{picture}(-22,80)
\put(123,55){$\{0\}$}
     \put(123,5){$\{e\}$}
    \put(-30,55){$\mathfrak{g}_\Sigma:= \mathfrak{g}\oplus \mathbb{R}$}
     \put(-11,5){$G _\Sigma$}
     \put(60,55){$\mathfrak{g}$}
     \put(60,5){$G$}
     \put(-96,55){$\mathbb{R}$}
     \put(-96,5){$\ker \pi_G$}
     \put(-165,55){$\{0\}$}
     \put(-165,5){$\{e\}$}
     \put(66,30){$\exp_{\mathfrak{g}}$}
       \put(-4,30){$\exp_{\mathfrak{g}_\Sigma}$}
         \put(-87,30){$\left.\exp_{\mathfrak{g}_\Sigma}\right|_{\mathbb{R}}$}
     \put(29,14){$\pi_G$}
     \put(20,63){$\pi_\mathfrak{g}=T _e\pi _G$}
\put(-150,57){\vector(1,0){48}}
     \put(-90,57){\vector(1,0){57}}
     \put(17,57){\vector(1,0){37}}
     \put(69,57){\vector(1,0){49}}
      \put(-8,50){\vector(0,-1){37}}
      \put(-93,50){\vector(0,-1){37}}
     \put(62,50){\vector(0,-1){37}}
     \put(-150,8){\vector(1,0){48}}
     \put(-70,8){\vector(1,0){57}}
     \put(7,8){\vector(1,0){47}}
     \put(69,8){\vector(1,0){49}}
     \end{picture}
\end{center} 
commutative. We notice that the connected component of the identity
$(\ker \pi _G)_0 $ of $\ker
\pi_G  $ equals
\begin{equation}
\label{thing with kernel and exponential} 
(\ker\pi_G)_0
=\{\exp_{\mathfrak{g}_\Sigma}(0,a)\mid a \in \mathbb{R} \}.
\end{equation} 
Indeed, for any $a \in \mathbb{R}$, $\pi_G
(\exp_{\mathfrak{g}_\Sigma}(0,a) )=\exp _{\mathfrak{g}} (0)=e $. This shows that
the one-dimensional connected Lie subgroup
$\{\exp_{\mathfrak{g}_\Sigma}(0,a)\mid a \in \mathbb{R} \}$ of $G _\Sigma$  is
included in
$(\ker \pi _G)_0 $. Since $\dim(\ker \pi _G)=1 $ it follows that
$\{\exp_{\mathfrak{g}_\Sigma}(0,a)\mid a
\in \mathbb{R} \}$ is open and hence closed in $(\ker \pi _G)_0 $ thus they are
equal. Additionally,
\begin{equation}
\label{kernel and centralizer} (\ker \pi _G)_0 \subset (Z(G _\Sigma)) _0,
\end{equation} where $(Z(G _\Sigma)) _0 $ denotes the connected component of the
identity of the center $Z(G _\Sigma)$ of $G _\Sigma$. To see  this, note
that any element $g \in G _\Sigma $  can be written as
$g = \exp_{\mathfrak{g}_\Sigma}(\xi_1, a_1)
\cdots\exp_{\mathfrak{g}_\Sigma}(\xi_1, a _1)$,
$\xi_1, \ldots, \xi_n \in  \mathfrak{g} $, $a _1, \ldots, a _n \in  
\mathbb{R}$ and hence, by the Baker-Campbell-Hausdorff formula, the
elements of the form  $\exp_{\mathfrak{g}_\Sigma}(0, a)$ commute with
every factor in the decomposition of $g$. 

Furthermore, if $Z (G)=\{e\} $ then~(\ref{kernel and centralizer}) is an
equality, that is, 
\begin{equation}
\label{kernel and centralizer equality} 
(\ker \pi _G)_0 = (Z(G _\Sigma)) _0.
\end{equation} 
In order to prove this fact recall that, by general theory, the map
$\exp _{\mathfrak{g}_\Sigma}$ is an isomorphism between the center of
$\mathfrak{g}_\Sigma  $ and $Z(G _\Sigma)_0 $. Since
$\mathfrak{g}$ has no center the dimension of $Z(G _\Sigma)$ equals one.
Additionally, since by \eqref{thing with kernel and
exponential} and~(\ref{kernel and centralizer}), $(\ker \pi _G)_0$ is  a
one-dimensional connected subgroup of 
$(Z(G _\Sigma)) _0$ the equality~(\ref{kernel and centralizer equality}) follows. 

\begin{proposition}
\label{properties of mu and more} Let $G$ be  a Lie group with  Lie algebra
$\mathfrak{g}$, $\Sigma$  a Lie algebra two-cocycle on $\mathfrak{g}$, and
$\mathfrak{g}_\Sigma 
$ the one-dimensional central extension of  $\mathfrak{g}$ determined by $\Sigma$.
Let $G _\Sigma $ be the connected and simply connected Lie group whose Lie algebra
is
$\mathfrak{g}_\Sigma$. There exists a smooth map $\mu_\Sigma:G _\Sigma
\longrightarrow
\mathfrak{g}^\ast$ such that  for any $g,h \in  G _\Sigma$, $(\xi, s)\in
\mathfrak{g}_\Sigma $,
$(\nu, a)\in \mathfrak{g}^\ast_\Sigma$ we have
\begin{itemize}
\item [{\rm {\bf (i)}}] $\mbox{\rm Ad} _g(\xi,s )= \left(\mbox{\rm Ad}_{\pi _G (g)}\xi,
s+\langle
\mu_\Sigma (g) , \xi\rangle 
\right)$.
\item [{\rm {\bf (ii)}}] $\mbox{\rm Ad}^\ast _{g ^{-1}}(\nu,a)=
\left(\mbox{\rm Ad}^\ast _{\pi_G(g) ^{-1}}
\nu+a \mu_\Sigma (g  ^{-1}),a\right)$.
\item [{\rm {\bf (iii)}}] $\mu_\Sigma(gh)= \mu_\Sigma(h)+ \mbox{\rm
Ad}^\ast _{\pi _G(h)}\mu_\Sigma (g)$.
\item [{\rm {\bf (iv)}}] $\mu_\Sigma(e)=0 $. 
\item [{\rm {\bf (v)}}] $\mu_{\Sigma}(g ^{-1})=- \mbox{\rm Ad}^\ast
_{\pi_G(g)^{-1}}\mu_{\Sigma}(g)$.
\item [{\rm {\bf (vi)}}] $T _g \mu_{\Sigma} \left(T _eL _g(\xi,s)
\right)=- \Sigma(\xi, \cdot )+ \mbox{\rm ad}^\ast _\xi \mu_{\Sigma} (g)$.
\end{itemize} 
In this statement we have identified $\mathfrak{g}^\ast_\Sigma $
with $\mathfrak{g}^\ast\oplus
\mathbb{R}  $ by using the pairing $\langle(\nu, a), (\xi,s) \rangle:=\langle \nu,
\xi\rangle+as
$. We will refer to $\mu_\Sigma $ as the {\bfi  extended
$\mathfrak{g}^\ast$-valued one-cocycle associated to} $\Sigma$.
\end{proposition}

\noindent\textbf{Proof.\ \ }We first notice that
\begin{align*}
\pi_{\mathfrak{g}}(\mbox{\rm Ad}_g(\xi,s))
&= \left.\frac{d}{dt}\right|_{t=0}
\pi_G(g \exp _{\mathfrak{g}_\Sigma} t (\xi,s)g^{-1})
=\left.\frac{d}{dt}\right|_{t=0} \pi_G(g)
\pi_G(\exp_{\mathfrak{g}_\Sigma} t (\xi,s)) \pi_G(g ^{-1})\\
	&=\left.\frac{d}{dt}\right|_{t=0} \pi_G(g) \exp _{\mathfrak{g}} (t \xi)
\pi_G(g ^{-1})= \mbox{\rm Ad}_{\pi_G (g)} \xi.
\end{align*} Consequently, there exists a smooth function $f :G _\Sigma \times
\mathfrak{g}\times \mathbb{R}
\longrightarrow\mathbb{R} $ such that  
\[
\mbox{\rm Ad}_g(\xi,s)= \left(\mbox{\rm Ad}_{\pi_G (g)} \xi, f(g, \xi,s) \right).
\] Note that for $g$ fixed  the map $f(g, \cdot , \cdot )$ is linear by the
linearity of $\mbox{\rm Ad}_g $, hence there exist smooth maps $\mu_\Sigma:G
_\Sigma \rightarrow  \mathfrak{g}^\ast$ and
$s _\Sigma:G _\Sigma \rightarrow \mathbb{R}$ such that 
\begin{equation}
\label{almost if f equations} f(g, \xi,s)=\langle \mu_{\Sigma}(g), \xi\rangle+ s
_\Sigma (g) s.
\end{equation} Additionally, by~(\ref{thing with kernel and exponential})
and~(\ref{kernel and centralizer}) the elements of the form $\exp
_{\mathfrak{g}_\Sigma} (0,s)$ belong to the center $Z(G _\Sigma)$ of $G _\Sigma$
and hence we have
\begin{equation*}
\mbox{\rm Ad}_g(0,s)= \left.\frac{d}{dt}\right|_{t=0} g \exp
_{\mathfrak{g}_\Sigma}t (0,s)g ^{-1}= \left.\frac{d}{dt}\right|_{t=0}\exp
_{\mathfrak{g}_\Sigma}t (0,s)=(0,s).
\end{equation*} 
Hence, combining this with~(\ref{almost if f equations}), we get 
$\mbox{\rm Ad} _g(0,s)=(0,s _\Sigma (g) s) =(0,s) $, which implies that
$s _\Sigma (g)=1 $, for any $g  \in G _\Sigma $. This proves point {\bf
(i)}.

\smallskip

\noindent {\bf (ii)} For $(\xi,s) \in  \mathfrak{g}_\Sigma $ we have by
{\bf (i)},
\begin{align*}
\langle  \mbox{\rm Ad}^\ast _{g ^{-1}}(\nu,a),(\xi,s)\rangle&=
\langle  (\nu,a),\mbox{\rm Ad} _{g ^{-1}}(\xi,s)\rangle=
\langle  (\nu,a),\left(\mbox{\rm Ad}_{\pi _G (g)^{-1}}\xi, s+\langle
\mu_\Sigma (g^{-1}) , \xi\rangle\right)\rangle\\
	&=\langle \nu,  \mbox{\rm Ad}_{\pi _G (g)^{-1}}\xi\rangle+ as+a\langle
\mu_\Sigma (g^{-1}) , \xi\rangle=\langle\mbox{\rm Ad}^\ast _{\pi_G(g) ^{-1}}
\nu+a \mu_\Sigma (g  ^{-1}), \xi\rangle+as,
\end{align*} which proves {\bf (ii)}.

\smallskip

\noindent {\bf (iii)} For any $g,h \in G _\Sigma $ and any $(\xi,s)\in
\mathfrak{g}_\Sigma$ we have  by {\bf (i)},
\begin{align*}
\left(\mbox{\rm Ad}_{\pi _G (gh)}\xi, s+\langle
\mu_\Sigma (gh) , \xi\rangle 
\right)&= \mbox{\rm Ad}_{gh}(\xi,s)= \mbox{\rm Ad}_{g}\left(\mbox{\rm
Ad}_{h}(\xi,s)\right)=
 \mbox{\rm Ad}_{g}\left(\mbox{\rm Ad}_{\pi _G (h)}\xi, s+\langle
\mu_\Sigma (h) , \xi\rangle\right)\\
	&= \left( \mbox{\rm Ad}_{\pi_G(g) \pi_G (h)}\xi,s+ \langle  \mu_{\Sigma} (h),
\xi\rangle
 +\langle  \mu_{\Sigma} (g), \mbox{\rm Ad}_{\pi_G(h)}\xi\rangle\right)\\
 	&= \left( \mbox{\rm Ad}_{\pi_G(g) \pi_G (h)}\xi,s+ \langle  \mu_{\Sigma} (h),
\xi\rangle
 +\langle   \mbox{\rm Ad}_{\pi_G(h)}^\ast \mu_{\Sigma} (g),\xi\rangle\right)
\end{align*} which implies {\bf (iii)}.

\smallskip

\noindent {\bf (iv)} By part {\bf (iii)} we can write $\mu_{\Sigma}(e)=
\mu_{\Sigma}(ee)=\mu_{\Sigma}(e)+\mu_{\Sigma}(e) $, and hence $\mu_{\Sigma}(e)=0 $.

\smallskip

\noindent {\bf (v)} By parts {\bf (iii)} and {\bf (iv)} we have  $0=
\mu_{\Sigma}(e)=
\mu_{\Sigma}(g g ^{-1})= \mu_{\Sigma}(g ^{-1})+ \mbox{\rm Ad}^\ast _{\pi_G
(g)^{-1}}\mu_{\Sigma} (g) $.

\smallskip

\noindent {\bf (vi)} For any $(\xi, s), (\eta,t)\in \mathfrak{g}_\Sigma $,
\begin{align*} ([\xi, \eta], -\Sigma(\xi, \eta))&=[(\xi,s),(\eta,t)]= \mbox{\rm
ad}_{(\xi,s)}(\eta,t)=
\left.\frac{d}{d\epsilon}\right|_{\epsilon=0}\mbox{\rm Ad}_{\exp
_{\mathfrak{g}_\Sigma}\epsilon(\xi,s)}(\eta,t)\\
	&= \left.\frac{d}{d\epsilon}\right|_{\epsilon=0} \left(\mbox{\rm Ad}_{\exp
_{\mathfrak{g}}
\epsilon \xi}\eta, t+\langle
\mu_\Sigma (\exp_{\mathfrak{g}_\Sigma }\epsilon(\xi,s)) , \eta\rangle\right) =
\left( [\xi, \eta], \langle T _e \mu_{\Sigma}(\xi,s), \eta
\rangle
\right),
\end{align*} which proves that 
\begin{equation*}
 T _e \mu_{\Sigma}(\xi,s)=- \Sigma(\xi, \cdot ).
\end{equation*} If we now use this equality together with {\bf (iii)} we obtain
\begin{align*} T _g \mu_{\Sigma} \left(T _eL _g(\xi,s)
\right)&=\left.\frac{d}{dt}\right|_{t=0}
\mu_{\Sigma}(g\exp _{\mathfrak{g}_\Sigma}t(\xi,s))=
\left.\frac{d}{dt}\right|_{t=0} \left(
\mu_\Sigma(\exp _{\mathfrak{g}_\Sigma}t(\xi,s))+ \mbox{\rm Ad}^\ast _{\exp
_{\mathfrak{g}}t
\xi}\mu_{\Sigma}(g)
\right)\\
	&=T _e \mu_\Sigma(\xi,s)+ \mbox{\rm ad}^\ast _\xi \mu_{\Sigma}(g)=
 - \Sigma(\xi, \cdot )+ \mbox{\rm ad}^\ast _\xi \mu_{\Sigma} (g). \quad
\blacksquare
\end{align*}

\medskip

\begin{corollary}
\label{extended and projected extended actions} Let $G$ be  a Lie group with  
Lie algebra $\mathfrak{g}$, $\Sigma$  a Lie algebra two-cocycle on
$\mathfrak{g}$, and $\mathfrak{g}_\Sigma$ the one-dimensional central
extension of $\mathfrak{g}$ determined by $\Sigma$. Let $G_\Sigma$ be
the connected and simply connected Lie group whose Lie algebra is
$\mathfrak{g}_\Sigma$ and  $\mu_\Sigma:G _\Sigma \longrightarrow
\mathfrak{g}^\ast$ the extended one-cocycle associated to $\Sigma$. The map
$\overline{\Xi}: G _\Sigma
\times  \mathfrak{g}^\ast \longrightarrow\mathfrak{g}^\ast $ defined by
\[
\overline{\Xi}(g, \mu):=\mbox{\rm Ad}^\ast _{\pi_G(g) ^{-1}}
\mu+\mu_\Sigma (g  ^{-1}), \quad g \in G _\Sigma, \mu \in \mathfrak{g}^\ast,
\] is an action of the Lie group $G _\Sigma$ on $\mathfrak{g}^\ast$. We will refer
to
$\overline{\Xi}$ as the {\bfi  extended affine action} of $G _\Sigma$ on
$\mathfrak{g}^\ast$. This action projects to a $G _\Sigma$-action on the quotient
$\mathfrak{g}^\ast/ \overline{{\cal H}} $ via the map $\Xi : G _\Sigma
\times  \mathfrak{g}^\ast/\overline{{\cal H}} \longrightarrow\mathfrak{g}^\ast/
\overline{{\cal H}} $ defined by
\[
\Xi(g, \mu+ \overline{{\cal H}}):=\mbox{\rm Ad}^\ast _{\pi_G(g) ^{-1}}
\mu+\mu_\Sigma (g  ^{-1})+ \overline{{\cal H}}, \quad g \in G _\Sigma, \mu \in
\mathfrak{g}^\ast,
\] with respect to which $\pi_C $ is $G_{\Sigma} $-equivariant. We will call
$\Xi$  the {\bfi projected extended affine action} of $G _\Sigma$ on
$\mathfrak{g}^\ast/ \overline{{\cal H}}$.
\end{corollary}

\noindent\textbf{Proof.\ \ }The proof of the first statement is a straightforward
verification. The projected extended affine action is well-defined by
Proposition~\ref{the Hamiltonian holonomy is Ad invariant}. \quad $\blacksquare$ 

\begin{remark}
\label{integration and affine actions}
\normalfont For any $\xi \in  \mathfrak{g}_\Sigma $, the associated infinitesimal
generators
$\xi_{\mathfrak{g}^\ast}$ and $\xi_{\mathfrak{g}^\ast/ \overline{{\cal H}}}$ of
the extended and the projected extended affine action of
$G_{\Sigma}
$  on $\mathfrak{g}^\ast $ and
$\mathfrak{g}^\ast/ \overline{{\cal H}} $, respectively, are given by the
expressions,
\begin{align}
\xi_{\mathfrak{g}^\ast} (\mu) &=- \mbox{\rm ad}^\ast
_{\pi_{\mathfrak{g}}(\xi)}\mu+ \Sigma(\pi_{\mathfrak{g}} (\xi), \cdot )=:
\left(\pi_{\mathfrak{g}} (\xi)
\right)_{\mathfrak{g}^\ast}(\mu),\label{infinitesimal generator extended}\\
\xi_{\mathfrak{g}^\ast / \overline{{\cal H}}}(\pi _C(\mu))&= T _\mu \pi _C \left(
-\mbox{\rm ad}^\ast _{\pi _{\mathfrak{g}}(\xi)}\mu+ \Sigma(\pi
_{\mathfrak{g}}(\xi), \cdot )
\right)=T _\mu \pi _C \left(\left(\pi_{\mathfrak{g}} (\xi)
\right)_{\mathfrak{g}^\ast}(\mu)\right),\label{infinitesimal generator extended
projected}
\end{align} for any $\mu \in \mathfrak{g}^\ast$. The second equality
in~(\ref{infinitesimal generator extended}) emphasizes that for any $\zeta \in
\mathfrak{g}$, the expression 
$\zeta_{\mathfrak{g}^\ast}(\mu):= - \mbox{\rm ad}^\ast _{\zeta}\mu+ \Sigma(\zeta,
\cdot )$ defines a Lie algebra action of $\mathfrak{g} $ on
$\mathfrak{g}^\ast$. Unlike the $\mathfrak{g}  _\Sigma$-action defined by the
first equality in~(\ref{infinitesimal generator extended}), the
$\mathfrak{g}$-action cannot, in general, be integrated to a group action. 
\end{remark}

\begin{remark}
\normalfont Suppose that there exists a right group one-cocycle $\sigma:G
\rightarrow  \mathfrak{g}^\ast$ that integrates $\Sigma$, that is, 
\begin{equation}
\label{integration condition sigma} 
T _e \sigma (\xi) = - \Sigma(\xi, \cdot ),
\end{equation} for any $\xi \in \mathfrak{g}$. This happens to be the case when,
for instance, the manifold underlying the group $G_{\Sigma} $ is diffeomorphic to
$G \times  \mathbb{R}  $ (see Chapter 2 of~\cite{roger}).   Then the extended
affine
$G_{\Sigma}$-action on
$\mathfrak{g}^\ast$ drops to the affine $G$-action $\overline{\Theta}:G \times 
\mathfrak{g}^\ast \rightarrow \mathfrak{g}^\ast$ on
$\mathfrak{g}^\ast$ determined by $\sigma$ via the expression
$
\overline{\Theta}(g, \mu)= \mbox{\rm Ad} ^\ast _{g ^{-1}}\mu+ \sigma(g ^{-1})
$, $(g, \mu)\in  G \times  \mathfrak{g}^\ast$. More specifically,
\begin{equation}
\label{relation between various actions}
\overline{\Xi}(\widehat{g}, \mu)= \overline{\Theta}(\pi _G(\widehat{g}), \mu),
\end{equation} for any $\widehat{g} \in G_{\Sigma}$  and  $\mu \in 
\mathfrak{g}^\ast$. This equality implies that
\[
\sigma\circ \pi_G= \mu_{\Sigma}. 
\] Analogously, the affine $G$-action $\Theta:G \times  \mathfrak{g}^\ast/
\overline{{\cal H}}
\rightarrow \mathfrak{g}^\ast/\overline{{\cal H}}$ on  $\mathfrak{g}^\ast/
\overline{{\cal H}} $ induced by $\sigma$   via the expression $\Theta(g, \mu+
\overline{{\cal H}})= \mbox{\rm Ad}^\ast _{g ^{-1}}\mu + \sigma (g ^{-1})+
\overline{{\cal H}} $, $(g, \mu)\in  G \times 
\mathfrak{g}^\ast $, is such that 
\begin{equation*}
\Xi(\widehat{g}, \mu+ \overline{{\cal H}})=\Theta(\pi _G(\widehat{g}), \mu+
\overline{{\cal H}}),
\end{equation*} for any $\widehat{g} \in G_{\Sigma}$  and  $\mu \in 
\mathfrak{g}^\ast$. In order to prove~(\ref{relation between various actions})
notice that~(\ref{infinitesimal generator extended}) and~(\ref{integration
condition sigma}) guarantee that 
\begin{equation}
\label{equality infinitesimal to prove}
\xi_{\mathfrak{g}^\ast} (\mu) =(\pi_{\mathfrak{g}}(\xi))_{\mathfrak{g}^\ast} (\mu),
\end{equation} for any $\xi\in  \mathfrak{g}_\Sigma$, $\mu \in 
\mathfrak{g}^\ast$. The infinitesimal generator on the left hand side of this
expression corresponds to the extended affine $G _\Sigma $-action
$\overline{\Xi} $ on $\mathfrak{g}^\ast$ while the one on the right hand side is
constructed using the affine
$G
$-action $\overline{\Theta} $ induced by $\sigma$. The equality~(\ref{equality
infinitesimal to prove}) implies~(\ref{relation between various actions}). Indeed,
since $G_{\Sigma} $ is connected, any element $\widehat{g} \in G_{\Sigma}$ can be
written as $\widehat{g}=\exp _{\mathfrak{g}_\Sigma} \xi _1 \cdots \exp
_{\mathfrak{g}_\Sigma} \xi _n $, $\xi_1, \ldots, \xi _n \in \mathfrak{g}_\Sigma$.
Hence
$\overline{\Xi}(\widehat{g}, \mu)= \overline{\Xi}(\exp _{\mathfrak{g}_\Sigma} \xi
_1 \cdots \exp _{\mathfrak{g}_\Sigma} \xi _n, \mu)=F _1^{\xi_1} \circ  \cdots
\circ F _1^{\xi_n} (\mu) $, where
$F _t^{\xi _i} $ is the flow of the infinitesimal generator vector field $(\xi
_i)_{\mathfrak{g}^\ast}$. Since by~(\ref{equality infinitesimal to prove})
$(\xi_i)_{\mathfrak{g}^\ast} =(\pi_{\mathfrak{g}}(\xi_i))_{\mathfrak{g}^\ast} $,
the flow $F _t^{\xi _i} $ coincides with the flow $F _t ^{\pi _{\mathfrak{g}}(\xi
_i)} $ of
$(\pi_{\mathfrak{g}}(\xi_i))_{\mathfrak{g}^\ast}$ and hence
\begin{align*}
\overline{\Xi}(\widehat{g}, \mu) &=F _1^{\pi _{\mathfrak{g}}(\xi_1)} \circ  \cdots
\circ F _1^{\pi _{\mathfrak{g}}(\xi_n)} (\mu)= \overline{\Theta}(\exp
_{\mathfrak{g}}(\pi _{\mathfrak{g}}(\xi_1) )\cdots \exp _{\mathfrak{g}}(\pi
_{\mathfrak{g}}(\xi_n)), \mu )\\
	&=\overline{\Theta}(\pi _G(\exp _{\mathfrak{g}_\Sigma}\xi_1) \cdots\pi _G(\exp
_{\mathfrak{g}_\Sigma}\xi_n), \mu )=\overline{\Theta}(\pi _G(\widehat{g}), \mu ),
\end{align*} which proves~(\ref{relation between various actions}). 
\end{remark}

\paragraph{The characteristic distributions of $(\mathfrak{g}^\ast, \{ \cdot ,
\cdot \}_\pm ^\Sigma)$ and 
$\left(\mathfrak{g}^\ast/\overline{{\mathcal H}},
\{ \cdot ,
\cdot
\}_{\mathfrak{g}^\ast/\overline{{\mathcal H}}}^{\pm \Sigma}\right)$.} The
characteristic distribution $\overline{E} $ of $(\mathfrak{g}^\ast, \{ \cdot ,
\cdot \}_\pm ^\Sigma)$ is given by (see for instance Section 4.5.27 in~\cite{hsr})
\begin{equation}
\label{characteristic previous}
\overline{E} (\mu)=\{ \mbox{\rm ad} ^\ast _\xi \mu- \Sigma(\xi, \cdot )\mid \xi\in
\mathfrak{g}\}.
\end{equation} We now show that the characteristic distribution $E$ of the Poisson
manifold
$\left(\mathfrak{g}^\ast/\overline{{\mathcal H}},
\{ \cdot ,
\cdot
\}_{\mathfrak{g}^\ast/\overline{{\mathcal H}}}^{\pm \Sigma}\right)$ is given by 
\begin{equation}
\label{characteristic distribution of projected Poisson Lie} E(\pi_C
(\mu))=\left\{ T _\mu \pi _C \left( \mbox{\rm ad}^\ast _\xi\mu- \Sigma(\xi, \cdot )
\right)\mid \xi\in \left({\rm Lie}\left(\overline{{\cal
H}}\right)\right)^\circ\right\}.
\end{equation} Indeed, since the projection $\pi _C:
\mathfrak{g}^{\ast}\rightarrow \mathfrak{g}^\ast /\overline{{\cal H}} $ is a
Poisson map, for any $f  \in
C^{\infty}(\mathfrak{g}^\ast/\overline{{\cal H}})$ and $\mu \in
\mathfrak{g}^\ast $ we have that $X _f(\pi_C (\mu))=T _\mu
\pi _C \left(X _{f \circ \pi_C} (\mu)\right) $. Since $X _{f \circ
\pi_C} $ is known (see for instance (4.5.20) in~\cite{hsr}) we have 
\begin{equation}
\label{vector field in projected} X _f(\pi_C (\mu))=T _\mu
\pi _C  \left(\mp \mbox{\rm ad}^\ast  _{\frac{\delta(f \circ
\pi_C)}{\delta
\mu}}\mu\pm
\Sigma\left(\frac{\delta(f \circ \pi_C)}{\delta \mu}, \cdot  \right)
\right).
\end{equation}

We now recall a result in~\cite{thesis} (see~\cite[Theorem 2.5.10]{hsr}) that
characterizes the span of the differentials of the invariant functions $f \in
C^\infty(M)^{G} $ with respect to a proper
$G$-action on a manifold $M$ by means of the equality
\begin{equation}
\label{thing in thesis with invariant}
\{ \mathbf{d}f (m)\mid f \in  C^\infty(M)^{G}\}= \left((\mathfrak{g}\cdot
m)^{\circ}\right)^{G _m}.
\end{equation} If we apply this result to the free and proper action of
$\overline{{\cal H}} $ on
$\mathfrak{g}^\ast$ we obtain that
\begin{equation}
\label{thesis for holonomy}
\left\{ \mathbf{d}h (\mu) \mid h \in C^{\infty}(\mathfrak{g}^\ast)^{\overline{{\cal
H}}}\right\}= \left( \left( {\rm Lie}\left(\overline{{\cal H}}\right)\right)
\cdot  \mu\right) ^{\circ}.
\end{equation} However, since for any element $\nu\in {\rm
Lie}\left(\overline{{\cal H}}\right) $ its associated infinitesimal generator
satisfies $ \nu_{ \mathfrak{g}^\ast} (\mu)= \nu $, for any $\mu \in 
\mathfrak{g}^\ast$, we can conclude that ${\rm Lie}\left(\overline{{\cal
H}}\right)\cdot  \mu = {\rm Lie}\left(\overline{{\cal H}}\right)$ and
hence~(\ref{thesis for holonomy}) can be rewritten as 
\begin{equation}
\label{thesis for holonomy short}
\left\{ \mathbf{d}h (\mu) \mid h \in C^{\infty}(\mathfrak{g}^\ast)^{\overline{{\cal
H}}}\right\}= \left({\rm Lie}\left(\overline{{\cal H}}\right)\right)^\circ.
\end{equation}  Now, since for any $\rho \in \mathfrak{g}^\ast$, we have
$\left\langle \rho, \frac{\delta (f \circ  \pi_C)}{\delta
\mu}\right\rangle= \mathbf{d}(f \circ \pi _C)(\mu)(\rho)$ and $f
\circ  \pi_C \in  C^{\infty}(\mathfrak{g}^\ast)^{\overline{{\cal H}}} $,
expressions~(\ref{vector field in projected}) and~(\ref{thesis for
holonomy short}) guarantee the validity of~(\ref{characteristic
distribution of projected Poisson Lie}). 

\paragraph{The symplectic leaves of $(\mathfrak{g}^\ast, \{ \cdot , \cdot \}_\pm
^\Sigma)$ and $\left(\mathfrak{g}^\ast/\overline{{\mathcal H}}, \{ \cdot ,
\cdot
\}_{\mathfrak{g}^\ast/\overline{{\mathcal H}}}^{\pm \Sigma}\right)$.} As we saw in
Lemma~\ref{holonomies are killers} the subspace $\mathfrak{n}:={\rm
Lie}\left(\overline{{\cal H}}\right)^{\circ} $ is an ideal of $\mathfrak{g}$. This implies
that $\mathfrak{n}_\Sigma:=
\mathfrak{n}\oplus \mathbb{R}  $ is an ideal of $\mathfrak{g}_\Sigma $. Let
$N_{\Sigma} $ be the connected and simply connected normal Lie subgroup of
$G_{\Sigma}$ whose Lie algebra is
$\mathfrak{n}_\Sigma$.

\begin{proposition}
\label{symplectic leaves with Hamiltonian holonomy} Let $(\mathfrak{g}^\ast, \{
\cdot , \cdot \}_\pm ^\Sigma)$ be an affine Lie-Poisson structure and
$\left(\mathfrak{g}^\ast/\overline{{\mathcal H}},\{\cdot ,
\cdot\}_{\mathfrak{g}^\ast/\overline{{\mathcal H}}}^{\pm \Sigma}\right)$  the
corresponding
 projected Lie-Poisson structure introduced in Proposition~\ref{projected lie
Poisson structures proposition}. Let $ N_{\Sigma} $ be the normal
subgroup of $G_{\Sigma}$  introduced in the paragraph above. Then for
any   $\mu \in  \mathfrak{g}^\ast $, the symplectic leaves
$\mathcal{L}_{ \mu}$ and $\mathcal{L}_{\pi _C(\mu)}$ of  
$(\mathfrak{g}^\ast, \{ \cdot , \cdot \}_\pm ^\Sigma)$ and
$\left(\mathfrak{g}^\ast/\overline{{\mathcal H}},
\{\cdot ,\cdot\}_{\mathfrak{g}^\ast/\overline{{\mathcal H}}}^{\pm \Sigma}\right)$
that contain
$\mu$ and $\pi_C(\mu) $, respectively, equal
\begin{equation}
\label{leaf of projected expression}
\mathcal{L}_\mu = \overline{\Xi}(G_{\Sigma} ,
\mu)\qquad\text{and}\qquad\mathcal{L}_{\pi _C(\mu)}=
\Xi(N_{\Sigma} , \pi _C(\mu)),
\end{equation} where $\overline{\Xi} $ and $\Xi $ are actions introduced in
Corollary~\ref{extended and projected extended actions}. 
\end{proposition}

\noindent\textbf{Proof.\ \ } If we compare~(\ref{infinitesimal generator
extended}) and~(\ref{infinitesimal generator extended projected})
with~(\ref{characteristic previous}) and~(\ref{characteristic distribution of
projected Poisson Lie}), it is clear that the tangent spaces to the orbits of the
$G_{\Sigma}$ and $N_{\Sigma}$-actions produce  distributions in
$\mathfrak{g}^\ast$ and $\mathfrak{g}^\ast /\overline{{\cal H}} 
$ equal to the characteristic distributions of $(\mathfrak{g}^\ast, \{
\cdot ,\cdot \}_\pm ^\Sigma)$ and
$\left(\mathfrak{g}^\ast/\overline{{\mathcal H}}, \{\cdot ,\cdot
\}_{\mathfrak{g}^\ast/\overline{{\mathcal H}}}^{\pm \Sigma}\right)$,  
respectively, and hence they have the same maximal integral leaves. Since
$G_{\Sigma}$ and  $ N_{\Sigma}  $ are connected these are the
$G_{\Sigma}  $ and  $N_{\Sigma}$-orbits, respectively.
\quad $\blacksquare$

\paragraph{Proof of Theorem~\ref{The Poisson structure of target}.} We will prove
the theorem  by showing that the Poisson structure in~(\ref{bracket on target
space}) is one of the projected Poisson structures presented in
Proposition~\ref{projected lie Poisson structures proposition} when we use a Lie
algebra two-cocycle $\Sigma: \mathfrak{g} \times \mathfrak{g}
\rightarrow  \mathbb{R}  $ that comes naturally with the construction that lead us
to the cylinder valued momentum map. We start by recalling a proposition whose
proof can be found in~\cite{coverings}.

\begin{proposition} Let $(M, \omega)$ be a connected and paracompact symplectic
manifold acted symplectically upon by the Lie group $G$. Let $\mathbf{K}:M
\rightarrow  \mathfrak{g}^\ast/\overline{{\mathcal H}} $ be a cylinder valued
momentum map for this action defined using the holonomy bundle
$\widetilde{M} \subset M \times  \mathfrak{g}^\ast$. Let $p: \widetilde{M}
\rightarrow  M $ be the projection. Then the pair $(\widetilde{M}, \omega
_{\widetilde{M}}:= \widetilde{p} ^\ast 
\omega)$ is a symplectic manifold where the Lie algebra $\mathfrak{g}$ acts
canonically via the map
\[
\xi_{\widetilde{M}}(m, \mu):=(\xi _M (m),-\Psi(m)(\xi, \cdot )),
\] for any $ \xi \in \mathfrak{g} $, $(m, \mu) \in \widetilde{M} $, and where
$\Psi: M
\rightarrow Z ^2 (\mathfrak{g})$ is the Chu map. The projection
$\widetilde{\mathbf{K}} :
\widetilde{M} \subset M \times  \mathfrak{g}^\ast\rightarrow \mathfrak{g}^\ast$ 
is a standard momentum map for this action with infinitesimal non-equivariance
cocycle $\widetilde{\Sigma}:
\mathfrak{g} \times  \mathfrak{g} \rightarrow \mathbb{R} $ given by 
\begin{equation}
\label{cocycle of lifted algebra action}
\widetilde{\Sigma}(\xi, \eta):= \langle \mu, [\xi, \eta]\rangle- \Psi (m) (\xi,
\eta),
\end{equation} for any $(m , \mu)\in  \widetilde{M} $. The value of  $ \widetilde{
\Sigma} $  does not depend on the point $(m, \mu)\in  \widetilde{M } $ used to
define it.
\end{proposition}

The Poisson structure in the statement of Theorem~\ref{The Poisson structure of
target} is one of the projected Poisson structures~(\ref{projected Poisson
structure target}) when we take in~(\ref{affine lie Poisson bracke sigma}) the
two-cocycle  $\widetilde{\Sigma} $ introduced in the statement of the previous
proposition. More specifically, for any
$f,g
\in  C^{\infty}(\mathfrak{g}^\ast)$ we have 
\begin{equation}
\label{characterization with cocycle up}
\{f,g\}_{\mathfrak{g}^\ast/\overline{{\cal H}}}(\pi _C (\mu))=\{f \circ \pi _C, g
\circ
\pi_C\}^{\widetilde{\Sigma}}_+ (\mu).
\end{equation} Indeed, using the independence of
$\widetilde{\Sigma} $ on the point $(m, \mu)$, we can write
\begin{align*}
\{f \circ \pi _C, g \circ\pi_C\}^{\widetilde{\Sigma}}_+ (\mu)&=\left\langle \mu,
\left[\frac{\delta(f \circ  \pi_C)}{\delta \mu},\frac{\delta(g \circ 
\pi_C)}{\delta
\mu}\right]\right\rangle  -\left\langle \mu,
\left[\frac{\delta(f \circ  \pi_C)}{\delta \mu},\frac{\delta(g \circ 
\pi_C)}{\delta
\mu}\right]\right\rangle\\
	&\ \  +
\Psi(m) \left(\frac{\delta(f \circ  \pi_C)}{\delta \mu},\frac{\delta(g \circ 
\pi_C)}{\delta
\mu} \right)= 
\Psi(m) \left(\frac{\delta(f \circ  \pi_C)}{\delta \mu},\frac{\delta(g \circ 
\pi_C)}{\delta
\mu} \right).
\end{align*}

\noindent {\bf Proof of point (i).} Let $f,g \in  C^{\infty}(\mathfrak{g}^\ast/
\overline{{\cal H}})$ be arbitrary functions. We will show that for any $m \in M $
\[
\{f \circ \mathbf{K}, g \circ \mathbf{K}\} (m)=\{ f , g
\}_{\mathfrak{g}^\ast/\overline{{\cal H}}}(\mathbf{K} (m)).
\] Indeed, by definition
\[
\{ f , g \}_{\mathfrak{g}^\ast/\overline{{\cal H}}}(\mathbf{K} (m))= \Psi(m)
\left(\frac{\delta(f \circ  \pi_C)}{\delta \mu},\frac{\delta(g \circ 
\pi_C)}{\delta
\mu} \right)= \omega (m) \left(\left(\frac{\delta(f \circ  \pi_C)}{\delta
\mu}\right)_M (m),\left(\frac{\delta(g
\circ 
\pi_C)}{\delta
\mu}\right)_M (m) \right),
\] 
where $ \mathbf{K} (m) = \pi _C(\mu) $. Since $\{f \circ \mathbf{K}, g
\circ \mathbf{K}\} (m)= \omega (m)(X_{f \circ \mathbf{K}}(m),X_{g \circ
\mathbf{K}}(m))$, it suffices to show that for any $f \in
C^{\infty}(\mathfrak{g}^\ast/
\overline{{\cal H}}) $
\begin{equation}
\label{to be proved vector field} X_{f \circ \mathbf{K}}(m)=\left(\frac{\delta(f
\circ  \pi_C)}{\delta
\mu}\right)_M (m).
\end{equation} This equality holds since by~(\ref{equation that determines nu}) 
\begin{multline*}
\omega (m)(X_{f \circ \mathbf{K}}(m), v _m)
= \mathbf{d}(f \circ \mathbf{K})(m)( v_m)=
\mathbf{d}f  (\mathbf{K} (m))(T _m \mathbf{K}(v _m))\\ 
= \mathbf{d} f (\mathbf{K}(m))(T _\mu \pi _C (\nu))
= \mathbf{d}(f \circ\pi _C)(\mu) (\nu)
=\left\langle\nu, \frac{\delta(f \circ  \pi _C)}{\delta\mu}\right\rangle
= \omega (m)\left(\left( \frac{\delta(f \circ  \pi _C)}{\delta \mu}
\right)_M  (m) , v _m\right)
\end{multline*} 
for any $v _m  \in T _m M $, which proves~(\ref{to be proved
vector field}). 

\medskip

\noindent {\bf Proof of point (ii).} $ \left({\rm Lie}\left(\overline{{\cal
H}}\right)\right)^\circ $ is an ideal in $\mathfrak{g}$ by Lemma~\ref{holonomies
are killers}. In order to prove the second statement note that since $N$  is
connected so are its orbits and hence it suffices to show that the distribution
given by the tangent spaces to the $N$-orbits coincides with
the characteristic distribution of the Poisson bracket~(\ref{bracket on
target space}). By~(\ref{characteristic distribution of projected
Poisson Lie}) and~(\ref{characterization with cocycle up}) the
characteristic distribution is given by 
\begin{equation}
\label{characteristic distribution real one} E(\pi_C (\mu))=\left\{ T _\mu \pi _C
\left(\Psi (m) (\xi, \cdot )
\right)\mid \xi\in \left({\rm Lie}\left(\overline{{\cal
H}}\right)\right)^\circ\right\},
\end{equation} where $\mathbf{K}(m)= \pi _C (\mu) $. Since, by
Proposition~\ref{cocycles for cylinder valued momentum map}, the
infinitesimal generators of the affine $G$-action on
$\mathfrak{g}^\ast/\overline{{\mathcal H}}$ are given by the expression
$\xi_{\mathfrak{g}^\ast/ \overline{{\mathcal H}}} (\pi _C (\mu))=-T _\mu
\pi _C \left(\Psi (m)(\xi, \cdot )\right)$,  for any $\xi \in 
\mathfrak{g}$, $(m, \mu)\in  \widetilde{M} $, the statement follows. 

Point {\bf (iii)} is a straightforward consequence of  {\bf
(ii)}.\quad $\blacksquare$

\section{The reduction theorems}
\label{The reduction theorems section}

The symplectic reduction procedure introduced by Marsden and Weinstein
in~\cite{mwr} consists of two steps. First, one restricts to the level sets of the
momentum  map and second, one projects it to the space of orbits of the group that
leaves that level set invariant.  The elements of the group that leave
a given level set of the momentum map invariant form a closed
subgroup of the original symmetry group. Indeed, if the manifold
is connected, one can always find an action on the dual of the Lie
algebra with respect to which the momentum map is
equivariant~\cite{souriau} and hence this subgroup is the isotropy
subgroup of the momentum value defining the level set.

In the preceding sections we have introduced all the necessary ingredients to
reproduce this construction in the context of the cylinder valued momentum map.
Nevertheless, if we blindly reproduce the Marsden-Weinstein construction it turns
out that we do not obtain symplectic reduced spaces but orbit spaces that are
endowed with a naturally defined Poisson structure that is, in general,
degenerate. As we will see, the reason behind
this surprising phenomenon is the eventual non-closedness of the
Hamiltonian holonomy in the dual of the Lie algebra. Indeed, when this
happens to be the case one may consider two different reduced spaces:
first, one based on the Marsden -Weinstein construction that suggests a
reduction by the $G$-action and that, as we already said, yields a
Poisson manifold; second, the foliation reduction theorem of
Cartan~\cite{cartan} (see Section 6.1.5 in~\cite{hsr} for a
self-contained presentation) imposes the reduction by the group $N$
introduced in Theorem~\ref{The Poisson structure of target} that
integrates $\left({\rm Lie}\left(\overline{{\cal
H}}\right)\right)^\circ$. The resulting reduced space is symplectic. The
two reduced spaces are related via a Poisson reduction procedure that 
will be described in detail and they coincide when the Hamiltonian
holonomy is closed in the dual of the Lie algebra.

Throughout this section all group actions are free and proper. This 
ensures the smoothness of all the orbit spaces that we will  encounter.
The generalization of these results to the context of non-free actions
is the subject of another  paper~\cite{stratified spaces}.

\begin{theorem}[Symplectic reduction]
\label{symplectic reduction theorem cylinder valued momentum map} Let $(M, \omega
)$ be a connected and paracompact symplectic manifold and $G$ a Lie group acting
freely, properly, and symplectically on it. Let
$\mathbf{K}:M\rightarrow  \mathfrak{g}^\ast/\overline{{\cal H}}$ be  a 
cylinder valued   momentum map for this action with associated
non-equivariance one-cocycle $\sigma:M \rightarrow 
\mathfrak{g}^\ast/ \overline{{\cal H}}$. Let $N$ be a normal connected Lie
subgroup of  $G$  that has
$\left({\rm Lie}\left(\overline{{\cal H}}\right)\right)^\circ $  as Lie 
algebra. Then for any $[\mu]\in \mathfrak{g}^\ast/\overline{{\cal H}}$
the orbit spaces $M_{[\mu]}:=\mathbf{K}^{-1}([\mu])/N_{[\mu]}$ are
regular quotient manifolds that are endowed with a natural symplectic
structure $\omega_{[\mu]}$ uniquely determined by the equality
\begin{equation}
\label{symplectic structure reduced cylinder valued momentum map} i ^\ast _{[\mu]}
\omega = \pi_{[\mu]} ^\ast  \omega_{[\mu]},
\end{equation} 
where $N_{[\mu]}$ denotes the isotropy subgroup of $[\mu]$ with respect
to the $N$-action on $\mathfrak{g}^\ast /\overline{{\cal H}}$ obtained by
restriction of the affine $G$-action constructed using the
non-equivariance cocycle $\sigma$ of  ${\bf K}$, $i _{[\mu]}:
\mathbf{K}^{-1}([\mu])\hookrightarrow M$ is the inclusion, and
$\pi_{[\mu]}: \mathbf{K}^{-1}([\mu]) \rightarrow 
\mathbf{K}^{-1}([\mu])/N_{[\mu]} $ is the projection. We will refer to
the spaces $M _{[\mu]}$ as the {\bfi  symplectic reduced spaces}.
\end{theorem}

The proof of this theorem requires  an intermediate result that generalizes in the
context of the cylinder valued momentum map the so called {\bfi   reduction lemma}.

\begin{proposition}[Reduction Lemma]
\label{reduction lemma cylinder valued momentum map} Let $(M, \omega )$ 
be a connected and paracompact symplectic manifold and $G$ a Lie group
acting symplectically on it. Let $\mathbf{K}:M \rightarrow 
\mathfrak{g}^\ast/\overline{{\mathcal H}} $ be  a cylinder valued  
momentum map for this action with associated non-equivariance one-cocycle
$\sigma:M \rightarrow  \mathfrak{g}^\ast/\overline{{\mathcal H}}$ and 
$N$ a normal connected Lie subgroup of  $G$  that has $\mathfrak{n}
:=\left({\rm Lie}\left(\overline{{\cal H}}\right)\right)^\circ $  as Lie
algebra. Then for any $m \in  M $ such that $\mathbf{K}(m)=
\pi_C(\mu)=:[ \mu ] $ we have
\begin{itemize}
\item  [{\rm {\bf (i)}}] $\mathfrak{g}_{[ \mu]} \cdot  m=
\ker T_m\mathbf{K}\cap \mathfrak{g}\cdot m$.
\item [{\rm {\bf (ii)}}] $\mathfrak{n}_{[\mu]}\cdot  m =\ker T _m
\mathbf{K}\cap  \left( \ker T _m
\mathbf{K}\right)^\omega=(\mathfrak{n} \cdot m)^\omega \cap \mathfrak{n}\cdot m$.
\item [{\rm {\bf (iii)}}] If the Hamiltonian holonomy ${\mathcal H} $ is
closed in
$\mathfrak{g}^\ast$ then $\mathfrak{g}_{[ \mu]} \cdot  m=(\mathfrak{g}\cdot m)
^\omega\cap \mathfrak{g}\cdot m $. 
\end{itemize}
\end{proposition}

\noindent\textbf{Proof.\ \ } Let  $\xi\in  \mathfrak{g}$. The equivariance of the
cylinder valued  momentum map ${\bf K}$ with respect to the affine action implies
that
\begin{equation}
\label{chain for reduction lemma} 
T_m \mathbf{K} \left(\xi_M (m)\right) =
\left.\frac{d}{dt}\right|_{t=0} \mathbf{K}(\exp\, t \xi\cdot   m)=
\left.\frac{d}{dt}\right|_{t=0} \Theta_{\exp\, t \xi} \mathbf{K} (m)=
\left.\frac{d}{dt}\right|_{t=0} \Theta_{\exp\, t \xi}[ \mu ]= \xi_{
\mathfrak{g}^\ast/
\overline{{\mathcal H}}}([\mu]).
\end{equation} This chain of equalities  shows that $ \xi_M(m) \in \ker T
_m\mathbf{K}\cap
\mathfrak{g}\cdot m$ if and only if $\xi \in  \mathfrak{g} _{[ \mu]} $ which
proves {\bf (i)}. As to {\bf (ii)}, the second equality is a consequence
of part {\bf (iii)} in Proposition~\ref{properties of the cylinder
valued momentum map in list}. This equality and {\bf (i)} imply the
first one. Part {\bf  (iii)} follows from part {\bf (i)} by noticing
that when ${\cal H}$ is closed in
$\mathfrak{g}^\ast$ then $\ker T _m \mathbf{K}=(\mathfrak{g}\cdot m)^\omega $ by
part {\bf (iii)} in Proposition~\ref{properties of the cylinder valued momentum
map in list}. \quad $\blacksquare$

\begin{remark}
\normalfont Notice that the reduction lemma shows how the reduced space that we
consider in Theorem~\ref{symplectic reduction theorem cylinder valued momentum
map} is the one that is hinted at in the foliation reduction theorem of
Cartan. This theorem studies the leaf space of the characteristic
distribution $\ker T 
\mathbf{K} \cap (\ker T  \mathbf{K} )^\omega$ in the level set
$\mathbf{K}^{-1}([\mu])$ that, by part {\bf  (ii)} of the previous
proposition, coincides with the distribution spanned by the tangent
spaces to the $N_{[\mu]}$-orbits. 
\end{remark}

\noindent {\bf Proof of Theorem~\ref{symplectic reduction theorem cylinder valued
momentum map}.} The freeness of the action implies, via part {\bf (iv)} of 
Proposition~\ref{properties of the cylinder valued momentum map in list}, that
${\bf K}$  is a submersion and hence the fiber
$\mathbf{K} ^{-1}([\mu])$  is a closed embedded submanifold of  $M$. Additionally,
the properness condition guarantees that the orbit space
$M_{[\mu]}:=\mathbf{K}^{-1}([\mu])/N_{[\mu]}
$ is a regular quotient manifold.

We now show that~(\ref{symplectic structure reduced cylinder valued momentum map})
defines a two-form on $M_{[\mu]} $. Let  $m, m' \in  M $, $v _m, w  _m \in T _m
\mathbf{K}^{-1}([\mu]) $,
$v _{m'}, w  _{m'} \in T _{m'} \mathbf{K}^{-1}([\mu]) $ be such that 
\begin{align}
\pi_{[\mu]}(m)&= \pi_{[\mu]}(m'),\label{first thing for well}\\ T _m \pi_{[\mu]}(v
_m) &= T_{m'} \pi_{[\mu]}(v_{m'})\quad\text{and}\quad T _m \pi_{[\mu]}(w _m) =
T_{m'} \pi_{[\mu]}(w_{m'}).\label{second thing for well}
\end{align} The equality~(\ref{first thing for well}) implies that there exists an
element $n \in  N_{[\mu]}$ such that $m'=n \cdot m $ and hence $T _m
\pi_{[\mu]}=\left. T_{m'}\pi_{[\mu]}\circ T _m\Phi\right|_{\ker T _m \mathbf{K}}$,
which, substituted in~(\ref{second thing for well}), yields
$\left(T_{m'}\pi_{[\mu]}\circ T _m\Phi\right)(v _m) =
T_{m'}\pi_{[\mu]}(v_{m'})$ and
$\left(T_{m'}\pi_{[\mu]}\circ T _m\Phi\right)(w _m) =
T_{m'}\pi_{[\mu]}(w_{m'})$. This in turn implies the existence of two
elements
$\xi, \eta \in  \mathfrak{n} _{[\mu]} $ such that $T _m \Phi(v _m)-v_{m'}= \xi_M
(m')$, $T _m
\Phi(w_m)-w_{m'}= \eta_M (m')$. Consequently
\[
\omega(m')(v_{m'},w_{m'})= \omega(\Phi_n (m))(T _m \Phi(v _m)-\xi_M (m'),T _m
\Phi(w_m)-\eta_M (m'))= \omega (m)(v _m, w _m).
\] In the previous chain of equalities we have used three things. First, by the
canonical character of the action, $\omega(\Phi_n (m))(T _m \Phi_n(v_m),
T _m \Phi _n(w _m)) = \omega(m)(v _m , w  _m)$. Second, the canonical
character of the action and the reduction lemma imply that
\begin{align*}
\omega(\Phi_n (m))(T _m \Phi _n(v _m), \eta_M (m'))
&= \omega(\Phi_n (m))(T _m \Phi_n(v _m), T _m \Phi _n\left(\mbox{\rm
Ad}_{n ^{-1}}\eta\right)_M (m))\\ 
&= \omega (m) (v _m, (\mbox{\rm Ad}_{n^{-1}}\eta)_M (m))=0.
\end{align*} 
Finally, for the same reasons, $\omega(\Phi_n (m))(T _m \Phi
_n(w _m), \xi_M (m'))=\omega(\Phi_n (m))(\eta_M (m'), \xi_M (m'))=0 $
which shows that~(\ref{symplectic structure reduced cylinder valued
momentum map}) defines a two-form $\omega_{[\mu]} $ on
$M_{[\mu]}$. Since $\pi_{[\mu]} $ is a surjective submersion and $\omega$ is
closed, so is $\omega_{[\mu]} $. In order to show that $\omega_{[\mu]} $
is non-degenerate let $v _m \in \ker T _m \mathbf{K}  $ be such that 
$\omega_{[\mu]}(\pi_{[\mu]}(m))(T _m \pi _{[\mu]}(v _m),T _m \pi _{[\mu]}(w _m))=0
$  for any $w _m \in  \ker T _m \mathbf{K} $. By~(\ref{symplectic structure
reduced cylinder valued momentum map}),  $\omega (m) (v _m, w _m)=0 $
for any $w _m \in  \ker T _m \mathbf{K} $. Consequently, by the reduction
lemma, $v _m  \in \left(\ker T _m \mathbf{K}\right) ^\omega\cap \ker
T_m\mathbf{K}= \mathfrak{n}_{[\mu]} \cdot m $ and hence $T _m\pi _{[\mu]}
(v_m) = 0$, as required. \quad $\blacksquare$

\begin{theorem}
\label{cylinder valued momentum map symplectic dynamics reduction} Let $(M, \omega
)$ be a connected and paracompact symplectic manifold and $G$ a Lie group acting
freely, properly, and symplectically on it. Let
$\mathbf{K}:M\rightarrow  \mathfrak{g}^\ast/\overline{{\cal H}}$ be  a cylinder
valued   momentum map for this action with associated non-equivariance one-cocycle
$\sigma:M \rightarrow 
\mathfrak{g}^\ast/ \overline{{\cal H}}$. Let $N$ be a normal connected Lie
subgroup of  $G$  that has
$\left({\rm Lie}\left(\overline{{\cal H}}\right)\right)^\circ $  as Lie algebra.
\begin{itemize}
\item [{\rm {\bf (i)}}] Let $F \in  C^\infty(M) ^N  $ be an
$N$-invariant function on $M$ and let $F _t $ be the flow of the
associated Hamiltonian vector field $X _F$. Then 
\[
\mathbf{K} \circ F _t= \mathbf{K}|_{{\rm Dom}(F _t)}.
\]
\item [{\rm {\bf (ii)}}]  Let $F \in  C^\infty(M) ^N  $ and let $M
_{[\mu]} $ be the symplectic reduced space introduced in
Theorem~\ref{symplectic reduction theorem cylinder valued momentum map},
for some $[\mu] \in\mathfrak{g}^\ast/
\overline{{\cal H}} $. Let $f  \in C^{\infty}\left(M _{[\mu]}\right)$
be the function uniquely determined by $f \circ \pi_{[\mu]} = 
F \circ i_{[\mu]}$. Then
\[ 
T \pi_{[\mu]} \circ X _F \circ i _{[\mu]}= X _f \circ \pi _{[\mu]}.
\]
\item  [{\rm {\bf (iii)}}] The bracket $\{\cdot,\cdot\}_{M_{[\mu]}}$
induced by $\omega_{[\mu]} $ on $M _{[\mu]} $ can be expressed as
\[
\{f, h\}_{M _{[\mu]}}(\pi_{[\mu]} (m)) = \{F , H\}(m),
\] 
with $F,H \in  C^\infty(M)^N $ two local $N$-invariant extensions at the
point $m$ of $f\circ\pi_{[\mu]}, h \circ \pi_{[\mu]} \in 
C^{\infty}(\mathbf{K} ^{-1}([\mu]))^{N _{[\mu]}}$ respectively (see the
remark below for an explanation of this terminology).
\end{itemize}
\end{theorem}

\begin{remark}
\label{remark on local extensions}
\normalfont Let $M$ be  a smooth manifold and $S$   an embedded submanifold of
$M$. Let $D\subset TM| _S $  be a subbundle of the tangent bundle of $M$
restricted to $S$ such that 
$D_S:=D\cap TS  $  is a smooth, integrable, regular distribution  on
$S$. Let $\pi _{D_S}:S \rightarrow S/D_S$ be the projection onto the leaf space of
the
$D_S$-distribution in  $S$. It can be proved (see lemma 10.4.14 in~\cite{hsr})
that for any open set
$V
\subset S/ D _S $, any function
$f \in C^{\infty} (V)$, and any $z \in V $  there exist a point $m \in \pi _{D_S}
^{-1} (V) $,  an open neighborhood $U _m$ of $m$ in  $M$, and
$F \in C^{\infty} _M (U _m) $ such that 
\[ 
f \circ  \pi_{D _S}|_{\pi _{D_S} ^{-1} (V)\cap U _m}=F|_{\pi _{D_S} ^{-1}
(V)\cap U _m},
\] and 
\[
\mathbf{d} F (n) |_{D(n)}=0, \quad \mbox{for any}\quad n \in  \pi _{D_S} ^{-1}
(V)\cap U _m.
\]  We say that $F$ is a {\bfi  local $D$-invariant extension} of $f \circ
\pi_{D_S} $ at the point $m \in \pi_{D_S} ^{-1} (V) $. 

The existence of the local $N$-invariant extensions in part {\bf (iii)} of  the
statement of Theorem~\ref{cylinder valued momentum map symplectic dynamics
reduction} is a consequence of the result that we just quoted by taking $S=
\mathbf{K}^{-1}([\mu])$ and $D (m) = \mathfrak{n} \cdot m$, $m \in  
\mathbf{K}^{-1}([\mu]) $. In this setup, the reduction lemma implies
that  $D_S (m)= \mathfrak{n}_{[\mu]} \cdot m$, $m \in
\mathbf{K}^{-1}([\mu])$. 
\end{remark}

\noindent {\bf Proof of Theorem~\ref{cylinder valued momentum map symplectic
dynamics reduction}.} To prove {\bf (i)} it suffices to show that $T
\mathbf{K} (X _F)=0 $. To see this, note that for any $m \in M $  such
that  $\mathbf{K} (m) = \pi  _C(\mu) $ we have that  $T _m\mathbf{K}
(X_F (m))=T _\mu\pi _C (\nu) $ where
\[
\langle \nu , \xi\rangle= \omega (m)\left(\xi  _M (m), X _F (m) \right)=-
\mathbf{d}F (m) \left(\xi _M (m)\right),
\]
for any $\xi \in \mathfrak{g}$. This equality and the $N$-invariance of
$F$ obviously imply that $\nu\in\mathfrak{n}^{\circ}= \left( \left({\rm
Lie}\left(\overline{{\cal H}}\right)\right)^\circ\right) ^{\circ}= {\rm
Lie}\left(\overline{{\cal H}}\right) $. Hence $T _\mu \pi _C (\nu) =0 $, as
required. Part {\bf (ii)} is a straightforward verification. To show
{\bf  (iii)} recall from Remark~\ref{remark on local extensions} that the
existence of the local $N$-invariant extensions is guaranteed by 
Lemma  10.4.14 in~\cite{hsr}. Moreover, by part {\bf (ii)},
\begin{align*}
\{ f  , h \}_{M _{[\mu]}}(\pi_{[\mu]} (m))&= \omega _{[\mu]}(\pi _{[\mu]} (m))(X
_f (\pi _{[\mu]} (m)), X _h(\pi _{[\mu]} (m)))\\
	&=\omega _{[\mu]}(\pi _{[\mu]} (m))(T _m \pi_{[\mu]}(X _F (m)), T _m \pi_{[\mu]}(X
_H (m)))= (\pi _{[\mu]}^\ast  \omega _{[\mu]}) (m)(X _F (m), X _H (m) ) \\
	&= \omega(m) (X _F (m), X _H (m))=\{F , H\}(m). \qquad \blacksquare
\end{align*} 

\paragraph{The optimal momentum map and optimal reduction. A quick overview.} As
stated in the introduction to this section, the analog of the
Marsden-Weinstein reduction  scheme in the cylinder valued momentum map
setup yields a Poisson manifold that is related to the symplectic
reduced space in Theorem~\ref{symplectic reduction theorem cylinder
valued momentum map} via  Poisson reduction. The description of the
symplectic leaves of this Poisson reduced space will require the use of
the   {\bfi  optimal momentum map}~\cite{optimal} and of the reduction
procedure that can be carried out with it~\cite{symplectic reduced} that
we quickly review in the following paragraphs. We refer to these
papers or to~\cite{hsr} for the proofs of the statements quoted below. 

Let $(M, \{ \cdot, \cdot\}) $ be a Poisson manifold and $G$ a Lie
group that acts properly on $M$ by Poisson diffeomorphisms via the left
action $\Phi: G \times M \rightarrow M $. The group of canonical
transformations associated to this action will be denoted by $A_{G}:=\{
\Phi _g: M \rightarrow M \mid g \in G\}$. Let $A_{G}'$ be the 
{\bfi $G$-characteristic\/} or the {\bfi polar distribution} on $M$
associated to $A_{G}$~\cite{dual pairs}  defined for any $m \in M$ by
$A_{G}'(m):=\{ X _f(m)\mid f \in C^\infty(M)^G\}$. The distribution
$A_{G}'$ is a smooth integrable  generalized distribution in the sense
of Stefan and Sussman~\cite{stefan, stefan b, sussman}. The {\bfi
optimal momentum map}
${\cal J}$ is defined as the canonical projection  onto the leaf space of
$A_{G}'$, that is,
$ {\cal J}:M \longrightarrow M / A_{G}'$. By its very definition, the levels sets
of ${\cal J}$ are preserved by the Hamiltonian flows associated to
$G$-invariant Hamiltonian functions and ${\cal J}$ is {\bfi universal}
with respect to this property, that is, any other map whose level sets
are preserved by
$G$-equivariant Hamiltonian dynamics factors necessarily through 
${\cal J}$. By construction, the fibers of
${\cal J}$ are the leaves of an integrable generalized distribution and thereby
{\it initial immersed submanifolds} of
$M$~\cite{dazord 1985}. Recall that $N$ is an {\bfi initial
submanifold\/} of $M$ when  the injection $i:N\rightarrow M$ is a smooth
immersion that satisfies the following property: for any manifold
$Z$, a mapping  $f:Z \rightarrow N $ is smooth if and only if $i\circ f:Z
\rightarrow M $ is smooth.

The leaf space $M / A_{G}'$ is called the {\bfi  momentum space} of ${\cal J}$. We
will consider it as a topological space with the quotient topology. If 
$m\in M$ let $\rho : = {\cal J}(m) \in M/A_{G}'$. Then, for any $g\in
G$, the map $\Psi_g(\rho)={\cal J}(g\cdot m)\in M/A_{G}'$ defines a
continuous $G$-action on $M/A_{G}'$ with respect to which ${\cal J}$ is
$G$-equivariant. Notice that since this action is not smooth and
$M/A_{G}'$ is not Hausdorff in general, there is no guarantee that the isotropy
subgroups $G _\rho$ are closed, and therefore embedded, subgroups of $G$. However,
there is a unique smooth structure on $G _\rho$ for which this subgroup becomes an
initial Lie subgroup of $G$ with Lie algebra $\mathfrak{g}_\rho$ given by 
\[
\mathfrak{g}_\rho=\{ \xi\in\mathfrak{g}\mid  \xi_M(m) \in T _m {\cal J}
^{-1}(\rho),  \mbox{ for all }  m \in {\cal J}^{-1}(\rho)\}.
\] With this smooth structure for $G  _\rho$, the left action
$\Phi^\rho: G _\rho\times {\cal J} ^{-1}(\rho)\rightarrow{\cal J} ^{-1}(\rho)$
defined by
$\Phi^\rho(g, z):= \Phi(g, z)$ is smooth.

\begin{theorem}[\cite{symplectic reduced}]
\label{Symplectic reduction by Poisson actions} Let $(M, \{ \cdot, \cdot\}) $ be a
smooth Poisson manifold and $G$ a Lie group acting canonically and
properly on
$M$. Let ${\cal J}:M
\rightarrow M/ A_{G}'$ be the optimal momentum map associated to this action.
Then, for any $\rho \in M / A_{G}'$ whose isotropy subgroup $G _\rho$ acts
properly on ${\cal J}^{-1}(\rho)$, the orbit space $M _\rho:={\cal J}^{-1}(\rho)/
G _\rho$ is a smooth symplectic regular quotient manifold with symplectic form
$\omega_\rho$ defined by:
\begin{equation}
\label{symplectic 1}
\pi_\rho^\ast\omega_\rho(m)(X _f(m), X _h(m))=\{f, h\} (m), \quad
\mbox{ for any } m \in \mathcal{J}^{-1}(\rho)  \mbox{ and any } f,h
\in C^\infty (M)^G.
\end{equation}  The map $\pi _\rho: \mathcal{J}^{-1}(\rho) \rightarrow 
\mathcal{J}^{-1}(\rho)/G_{\rho} $  is the projection.
\end{theorem}

Suppose now that the $G$-action is free and proper. It is well known that the
orbit space $M/G$ is a Poisson manifold with  the Poisson bracket
$\{\cdot,\cdot\}_{M/G}$, uniquely characterized by the relation
\begin{equation}
\label{characterization reduced free bracket}
\{f,\,g\}_{M/G}(\pi (m))=\{f\circ\pi,\,g\circ\pi\}(m),
\end{equation} for any $m \in M $ and where  $f,\,g:M/G\rightarrow\mathbb{R}$ are
two arbitrary smooth  functions. A fact that we will use in the sequel is that the
symplectic leaves of $(M/G,
\{ \cdot , \cdot \}_{M/G})$ are given by the optimal orbit reduced spaces $ \left(
\mathcal{J}^{-1}(\mathcal{O}_{\rho})/G, \omega _{\mathcal{O}_\rho} \right)$,
$\mathcal{O}_{\rho}:=G \cdot \rho$, $\rho \in M / A_{G}'$,  that are symplectically
diffeomorphic to the optimal point reduced spaces introduced in
Theorem~\ref{Symplectic reduction by Poisson actions} via the map $\pi_\rho (m)
\longmapsto \pi _{\mathcal{O}_\rho} (m)
$, with $m \in   \mathcal{J}^{-1}(\rho) $ and $\pi_{\mathcal{O}_\rho}:
\mathcal{J}^{-1}(\mathcal{O}_{\rho}) \rightarrow 
\mathcal{J}^{-1}(\mathcal{O}_{\rho})/G $ the projection.

\begin{theorem}
\label{all reductions cylinder valued momentum map} Let $(M, \omega )$ be a
connected and paracompact symplectic manifold and $G$ a Lie group acting freely,
properly, and symplectically on it. Let
$\mathbf{K}:M\rightarrow  \mathfrak{g}^\ast/\overline{{\cal H}}$ be  a cylinder
valued   momentum map for this action with associated non-equivariance one-cocycle
$\sigma:M \rightarrow 
\mathfrak{g}^\ast/ \overline{{\cal H}}$. Let $N$ be a normal connected Lie
subgroup of  $G$  that has
$\left({\rm Lie}\left(\overline{{\cal H}}\right)\right)^\circ $  as Lie algebra.
Then for any
$[\mu]\in \mathfrak{g}^\ast/\overline{{\cal H}}$: 
\begin{itemize}
\item [{\rm {\bf (i)}}] The orbit space
$M^{[\mu]}:=\mathbf{K}^{-1}([\mu])/G_{[\mu]}$ is a regular quotient manifold  
endowed with a natural Poisson structure induced by the bracket $\{ \cdot , \cdot
\}_{M ^{[\mu]}}$ determined by the expression
\begin{equation}
\label{Poisson structure reduced cylinder valued momentum map}
\{ f , h \}_{M ^{[\mu]}}(\pi^{[\mu]}(m))=\{F,H\} (m),
\end{equation} where $G_{[\mu]}$ denotes the isotropy subgroup of $[\mu]$ with
respect to the affine $G$-action on $\mathfrak{g}^\ast/ \overline{{\cal H}} $
constructed using the non-equivariance cocycle
$\sigma$ of  ${\bf K}$, $\pi^{[\mu]}:
\mathbf{K}^{-1}([\mu]) \rightarrow  \mathbf{K}^{-1}([\mu])/G_{[\mu]} $ is the
projection, and $F,H \in  C^\infty(M)^{G}  $ are local $G$-invariant
extensions of $f \circ \pi ^{[\mu]} $  and $h
\circ \pi ^{[\mu]} $ around the point $m$, respectively.  We will refer
to the  spaces $M ^{[\mu]}$ as the {\bfi  Poisson
reduced spaces}.
\item [{\rm {\bf (ii)}}] The Lie group $H _{[\mu]}:=G _{[\mu]}/ N
_{[\mu]} $ acts canonically, freely, and properly on  $(M _{[\mu]},
\omega _{[\mu]})$. The reduced Poisson manifold $\left(M _{[\mu]}/ H
_{[\mu]}, \{ \cdot , \cdot \}_{H _{[\mu]}}
\right) $ is Poisson isomorphic to $(M ^{[\mu]}, \{ \cdot  , \cdot  \} _{M
^{[\mu]}})$ via the map 
\begin{equation*}
\begin{array}{cccc}
\Psi:&M _{[\mu]}/ H_{[\mu]}&\longrightarrow &M ^{[\mu]}\\
	&\pi_{H _{[\mu]}}(\pi _{[\mu]}(m))&\longmapsto &\pi ^{[\mu]} (m),
\end{array}
\end{equation*} where $\pi_{H _{[\mu]}}:M _{[\mu]} \rightarrow M _{[\mu]}/ H
_{[\mu]} $  is the projection.
\item [{\rm {\bf (iii)}}] Let $\mathcal{J}_{H _{[\mu]}}: M _{[\mu]}
\rightarrow M _{[\mu]}/ A_{H _{[\mu]}}'$ be the optimal momentum map
associated to the $H _{[\mu]}
$-action on $M _{[\mu]}$. Suppose that for any $\rho \in  M _{[\mu]}/ A_{H
_{[\mu]}}'$, the isotropy subgroup $\left(H _{[\mu]}\right)_\rho $ acts
properly on the  level set $\mathcal{J}_{H _{[\mu]}}^{-1}(\rho)$. Then
the symplectic  leaves of $(M ^{[\mu]}, \{\cdot ,\cdot\}_{M ^{[\mu]}})$ 
are given by the manifolds $\Psi\left(\mathcal{J}_{H_{[\mu]}}^{-1}
(\mathcal{O}_{\rho})/H _{[\mu]}\right)$, for any $\rho
\in  M _{[\mu]}/ A_{H _{[\mu]}}'$.
\item [{\rm {\bf (iv)}}] Let  $\mathcal{J} :M \rightarrow  M/ A_{G}' $
be the optimal momentum map associated to the $G$-action on $M$ and let
$m \in  M $ be such that
$\mathbf{K} (m)= [\mu] $. Then $\mathcal{J}^{-1}(\rho)\subset
\mathbf{K}^{-1}([\mu])$ and $G_{\rho} \subset G _{[\mu]}$. If the 
isotropy subgroup $G_{\rho} $ acts properly on $\mathcal{J}^{-1}(\rho) $
then $(M _\rho:=\mathcal{J}^{-1}(\rho)/ G_{\rho}, \omega_\rho)$ is a
smooth regular quotient symplectic manifold and the map 
\[
\begin{array}{cccc} 
L :&\mathcal{J}^{-1}(\rho)/
G_{\rho}&\longrightarrow&\mathbf{K} ^{-1}([\mu])/ G _{[\mu]}\\
	& \pi_\rho (z)&\longmapsto &\pi ^{[\mu]} (z)
\end{array}
\] is  well defined, smooth, injective, and its image is the
symplectic leaf $\mathcal{L} _{\pi^{[\mu]}(m)}$ of $\left(M ^{[\mu]},\{
\cdot , \cdot \}_{M ^{[\mu]}}\right)$ that contains $\pi^{[\mu]}(m)$. If,
additionally, $\mathcal{J}^{-1}(\rho)$ is a closed subset of
$M$  then $L$ is a Poisson map and the Poisson bracket  $\{ \cdot ,
\cdot\}_{M _\rho}$ induced by the symplectic form $\omega _\rho  $ is determined
by the expression
\begin{equation}
\label{Poisson structure reduced optimal momentum map}
\{ f , h \}_{M_\rho}(\pi_\rho(z))=\{F,H\} (z),
\end{equation} where $F,H
\in  C^\infty(M)^{G}  $ are local $G$-invariant extensions of $f \circ \pi _\rho
$  and $h
\circ \pi _\rho $ around the point $z \in \mathcal{J}^{-1}(\rho)$,
respectively.
\end{itemize}
\end{theorem}

\begin{remark}
\normalfont This theorem links all the three reduction schemes induced
by a free, proper, and canonical  action on a symplectic manifold. On
one hand, we have the two possible reductions using the cylinder valued
momentum map: $M _{[\mu]} $ the symplectic  and $M^{[\mu]}$ the Poisson 
reduced spaces that are related to each other via Poisson reduction by
the quotient group $G _{[\mu]} / N _{[\mu]} $. On the other hand, we can
carry out optimal reduction; point {\bf (iv)} in the previous statement
shows that the optimal reduced spaces are the symplectic leaves of the
Poisson reduced space $M ^{[\mu]}$.

All these three reduced spaces are, in general, different (see the
example below). Nevertheless, we note that  if
${\mathcal H} $ is closed in
$\mathfrak{g}^\ast$ then
$\left({\rm Lie}\left(\overline{{\cal H}}\right)\right)^\circ =\left({\rm
Lie}\left( {\cal H}\right)\right)^\circ = \mathfrak{g}$ and hence $G _{[\mu]}= N
_{[\mu]} $,
$H_{[\mu]}=\{e\} $,   $M _{[\mu]}=M ^{[\mu]} $. Moreover, in this situation, the
closedness hypothesis needed in point {\bf (iv)} of the statement of the theorem
always holds and the optimal reduced spaces are the connected components of $M
_{[\mu]}=M ^{[\mu]}
$. This is so because whenever  ${\mathcal H}
$  is closed in $\mathfrak{g}^\ast$ then, by Proposition~\ref{properties of the
cylinder valued momentum map in list}, $\ker T _z \mathbf{K} =(\mathfrak{g}\cdot 
z)^\omega $, $z \in  M $, and hence $T _z \mathcal{J}^{-1}(\rho) =\ker T
_z
\mathbf{K} $, $ z \in  \mathcal{J}^{-1}(\rho) $, which shows that $
\mathcal{J}^{-1}(\rho) $ is one of the connected components of
$\mathbf{K}^{-1} ([\mu])$ and hence is closed in $M$. 

We emphasize that the closedness of ${\mathcal H} $ in $\mathfrak{g}^\ast$ is a
sufficient, but in general not necessary, condition for the three reduced
spaces to coincide (see Example~\ref{torus example with open holonomy}).
\end{remark}

\begin{example}
\normalfont The following elementary example shows that the three reduced spaces
in the statement of theorems~\ref{symplectic reduction theorem cylinder valued
momentum map} and~\ref{all reductions cylinder valued momentum map}, that is, the
Poisson, the symplectic, and the optimal reduced spaces, are in general distinct.
Let
$M:=\Bbb T^2\times
\Bbb T^2$ be the product of two tori whose elements will be  denoted by the
four-tuples $(e ^{i\theta_1},e ^{i\theta_2},e ^{i \psi_1},
e^{i\psi_2})$. Endow $M$ with the symplectic structure $\omega$ defined
by $\omega:= \mathbf{d} \theta_1\wedge \mathbf{d} \theta_2+ \sqrt{2}\,
\mathbf{d} \psi_1\wedge \mathbf{d} \psi_2. 
$ Consider the canonical circle action given by $e ^{i \phi} \cdot (e ^{i
\theta_1},e ^{i \theta_2},e ^{i \psi_1},e ^{i \psi_2}):=(e ^{i (\theta_1+ \phi)},e
^{i \theta_2},e ^{i (\psi_1+ \phi)},e ^{i
\psi_2})$ and the trivial principal bundle $(M
\times \mathbb{R}) \rightarrow M $ with $(\mathbb{R},+)$ as structure group. It is 
easy to see that the horizontal vectors in
$T(M\times
\mathbb{R})$ with  respect to the connection $\alpha$ defined in~(\ref{definition
of alpha connection}) are of the form $((a _1, a _2,b_1, b _2), -a _2-\sqrt{2}b
_2) $, with
$a_1, a _2, b_1, b _2 
\in
\mathbb{R}$. The surfaces $\widetilde{ M}_\tau \subset M
\times \mathbb{R}$ of the form
$\widetilde{ M}_\tau:=\{((e ^{i
\theta_1},e ^{i \theta_2},e ^{i \psi_1},e ^{i \psi_2}), \tau-
\theta_2-\sqrt{2}\psi_2)\in  M\times \mathbb{R}\mid
\theta_1,
\theta_2, \psi_1, \psi_2 \in \mathbb{R}\}$ integrate the horizontal distribution
spanned by these vectors. 

Take now $\widetilde{M}:= \widetilde{M} _0$ and consider the projection
$\widetilde{p} : \widetilde{M} \rightarrow M $. It is clear that $\widetilde{p}
^{-1}(e ^{i
\theta_1},e ^{i \theta_2},e ^{i \psi_1},e ^{i \psi_2})=\{(e ^{i
\theta_1},e ^{i \theta_2},e ^{i \psi_1},e ^{i \psi_2},- (\theta_2+2 n
\pi)-\sqrt{2}(\psi_2+2 m \pi))\mid m,n \in \Bbb Z\}$. Since the Hamiltonian
holonomy ${\mathcal H} $ coincides with the structure group of the fibration
$\widetilde{p} : \widetilde{M} \rightarrow M $, it follows that 
${\mathcal H}= \Bbb Z+\sqrt{2} \Bbb Z \subset \mathbb{R}$. In this case ${\mathcal
H} $ is indeed not closed; moreover ${\mathcal H}$ is dense in
$\mathbb{R} $, that is, $\overline{{\mathcal H}}= \mathbb{R}$. Therefore, in this
case, the cylinder valued momentum map is a constant map since its range is just a
point and the group
$N=\{e\}$. Hence the symplectic reduced space $M _{[\mu]}$ equals the entire
symplectic manifold $\mathbb{T}^4 $ and the Poisson reduced space $M ^{[\mu]} $
equals the orbit space
$\mathbb{T}^4/S ^1 $.

We now compute the optimal reduced spaces. In this case, 
$C ^{\infty}(M)^{S ^1}$ consists  of all the functions
$f$ of the form $f(e ^{i\theta_1},e ^{i \theta_2},e ^{i \psi_1},e ^{i \psi_2})
\equiv g (e ^{i \theta_2},e ^{i \psi_2}, e ^{i (\theta_1-\psi_1)})$, for some
function $g
\in C^\infty(\mathbb{T}^3)$. An inspection of the Hamiltonian flows associated to
such functions readily shows that the leaves of  $A _{S ^1}' $, that is, the level
sets
${\cal J}^{-1}(\rho)$ of the optimal momentum map ${\cal J}$, are the product of a
two-torus with a leaf of an irrational foliation (Kronecker submanifold)
of another two-torus. The isotropy subgroups $S^1_\rho$ coincide with
the circle $S ^1$, whose compactness guarantees that its action on ${\cal
J}^{-1}(\rho)$ is proper. Theorem~\ref{Symplectic reduction by Poisson
actions} automatically guarantees that the quotients 
\[ 
M _\rho:= {\cal J}^{-1}(\rho)/ S ^1_\rho\simeq \left(\mathbb{T}^2
\times \text{\{Kronecker submanifold of $\Bbb T^2$\}}\right)/S ^1.
\]  
are symplectic and, by Theorem~\ref{all reductions cylinder valued
momentum map}, they are the symplectic  leaves of the quotient Poisson
manifold $\mathbb{T}^4/ S^1$. 
\end{example}

\noindent\textbf{Proof of Theorem~\ref{all reductions cylinder valued momentum
map}.\ \ (i)} The smooth structure of $M ^{[\mu]} $ as a regular quotient 
manifold is a consequence of the freeness and properness of the
$G$-action, using the same arguments as in the proof of
theorem~\ref{symplectic reduction theorem cylinder valued momentum
map}.  Let $S:= \mathbf{K}^{-1} ([\mu]) $   and $D\subset TM| _S$ be
defined by $D(m) :=
\mathfrak{g}\cdot m $, $ m \in  M $. The reduction lemma guarantees that the
distribution $D_S:= D \cap TS $ in $S$ coincides with the tangent spaces to the
$G_{[\mu]}$-orbits and hence lemma 10.4.14 in~\cite{hsr} guarantees that
the local $G$-invariant extensions used in the expression~(\ref{Poisson
structure reduced cylinder valued momentum map}) do exist. Let
$B^\sharp: T ^\ast M \rightarrow  TM$ be vector bundle isomorphism
induced by the symplectic  form $\omega$  on $M$. Since for any $m  \in 
\mathbf{K}^{-1} ([\mu]) $
\[ 
B^\sharp (m)((D(m)) ^{\circ})=( D (m) ) ^\omega= (\mathfrak{g}\cdot 
m)^\omega
\subset (\mathfrak{n} \cdot  m)^\omega= \ker T _m \mathbf{K},
\] 
it is clear that $B ^\sharp (D^{\circ})\subset TS+ D $ and hence the
Marsden-Ratiu Theorem on Poisson reduction~\cite{poissonreduction} guarantees
that~(\ref{Poisson structure reduced cylinder valued momentum map}) is a well
defined bracket.

\smallskip

\noindent {\bf (ii)} The $H _{[\mu]} $-action on  $M _{[\mu]}  $ is given by 
\[ G _{[\mu]} \cdot \pi _{[\mu]} (m):= \pi _{[\mu]}(g \cdot   m), \quad g  \in G
_{[\mu]}.
\] 
This action is obviously canonical, free, and proper, and hence
$M_{[\mu]}/H _{[\mu]}$  is  a smooth Poisson manifold such that  the
projection $\pi _{H_{[\mu]}} $ is  a Poisson surjective submersion. It is
easily verified that $\Psi$ is a well-defined smooth bijective map  with
smooth inverse given by
\[
\pi^{[\mu]} (m) \longmapsto \pi _{H _{[\mu]}}(\pi _{[\mu]} (m)), \quad m  \in 
\mathbf{K} ^{-1}([\mu]).
\] 
Consequently, $\Psi$ is a diffeomorphism. In order to show that $\Psi$ is
Poisson let
$\overline{ \Psi}:= \Psi \circ \pi_{H _{[\mu]}}:M _{[\mu]} \rightarrow  M ^{[\mu]}
$. Notice that $\overline{\Psi} $ makes the diagram
\unitlength=5mm
\begin{center}
\begin{picture}(10,6)
\put(1,5){\makebox(0,0){$M _{[\mu]}$}} 
\put(9.2,5){\makebox(0,0){$\mathbf{K} ^{-1}([\mu])$}}
\put(1,0){\makebox(0,0){$M _{[\mu]}/ H _{[\mu]}$}} 
\put(9.1,0){\makebox(0,0){$M ^{[\mu]}$}} 
\put(1,4){\vector(0,-1){3}}
\put(2,4){\vector(2,-1){6}}
\put(9,4){\vector(0,-1){3}}
\put(7.5,5){\vector(-1,0){5.7}}
\put(2.6,0){\vector(1,0){5.4}}
\put(9.7,2.5){\makebox(0,0){$\pi ^{[\mu]}$}}
\put(0.1,2.5){\makebox(0,0){$\pi_{H _{[\mu]}}$}}
\put(5,5.6){\makebox(0,0){$\pi _{[\mu]}$}}
\put(5,.7){\makebox(0,0){$\Psi$}}
\put(5,3.5){\makebox(0,0){$\overline{\Psi}$}}
\end{picture}
\end{center}
\medskip  commutative. Since $\pi_{[\mu]} $ and $\pi_{H _{[\mu]}}$ are
surjective submersions it follows that $\Psi$  is Poisson if and only if
\[
\{f,h\}_{M ^{[\mu]}}\circ  \Psi \circ \pi_{H _{[\mu]}}\circ  \pi _{[\mu]}=\{ f
\circ \Psi\circ
\pi_{H _{[\mu]}}, h\circ \Psi\circ
\pi_{H _{[\mu]}}\}_{M _{[\mu]}}\circ \pi _{[\mu]},
\] which is equivalent to 
\begin{equation}
\label{thing to prove is equivalent to}
\{f,h\}_{M ^{[\mu]}} (\pi ^{[\mu]}(m) )=\{f \circ \overline{\Psi}, h \circ
\overline{\Psi}\}_{M _{[\mu]} }(\pi _{[\mu]}(m) ),
\end{equation} for any $f,h \in  C^{\infty}(M ^{[\mu]})$ and $ m \in  \mathbf{K}
^{-1} ([\mu])$. By part {\bf (i)}, the left hand side of~(\ref{thing to prove is
equivalent to}) equals $\{F,H\} $, with $F,H
\in  C^\infty(M)^{G}   $ local $G$-invariant extensions at  $m$ of  $f \circ 
\pi^{[\mu]} $ and
$g\circ\pi ^{[\mu]}  $, respectively. Additionally, since the bracket $\{ \cdot  ,
\cdot  \}_{M _{[\mu]}}$ is induced by the symplectic  form $\omega_{[\mu]} $ it
equals, by part {\bf  (iii)} in theorem~\ref{cylinder valued momentum map
symplectic dynamics reduction}, $\{
\overline{F},
\overline{ H}\} (m)$ with  $\overline{F},
\overline{ H} \in C^\infty(M)^{N} $ local $N$-invariant extensions at $m$  of $f
\circ
\overline{\Psi} \circ  \pi _{[\mu]}$ and $h \circ
\overline{\Psi} \circ  \pi _{[\mu]}$, respectively. Since $\overline{\Psi} \circ
\pi _{[\mu]}=
\pi^{[\mu]} $, $\overline{F} $  and $\overline{H} $ can be taken to be $F$ and $H$,
respectively, which proves~(\ref{thing to prove is equivalent to}).

\smallskip

\noindent {\bf (iii)} This follows from part {\bf (iv)} in Theorem 10.1.1
of~\cite{hsr}. 

\smallskip

\noindent {\bf (iv)} By Proposition~\ref{properties of the cylinder valued
momentum map in list} we can think of the connected components of the level sets
of $\mathbf{K} $ as the maximal integral manifolds of the distribution $D$ in $M$
given by $D(z)=(\mathfrak{n}\cdot  z)^\omega$, $z \in  M $. Since
$(\mathfrak{n}\cdot z)^\omega\supset (\mathfrak{g}\cdot z)^\omega= A_{G}'(z)$ and
the level sets of the optimal momentum map $\mathcal{J}$ are the maximal integral
manifolds of $A_{G}' $, we obviously have that $\mathcal{J}^{-1}(\rho)
\subset \mathbf{K}^{-1}([\mu])$. Let now $g \in G_{\rho}$ and $m':=
\mathfrak{g}\cdot m $. Since $\mathcal{J} $ is $G$-equivariant we have
$\mathcal{J}(m')= \mathcal{J} (g \cdot  m)=g\cdot  \rho= \rho$ and
hence $m' \in  \mathbf{K} ^{-1}([\mu])$ which in turn  implies that $
[\mu]= \mathbf{K} (m')= \mathbf{K}(g \cdot   m)= g \cdot 
\mathbf{K} (m)=g \cdot  [\mu] $ and hence guarantees that $g \in  G _{[\mu]} $. 

If  $ G_{\rho} $ acts freely and properly on $ \mathcal{J}^{-1}(\rho) $ then the
hypotheses of the optimal reduction theorem~\ref{Symplectic reduction by Poisson
actions} are satisfied. The map $L $ is just the projection of the $(G_{\rho} , G
_{[\mu]})$-equivariant inclusion
$\mathcal{J}^{-1}(\rho)\hookrightarrow \mathbf{K} ^{-1}([\mu])$ and is hence well
defined and smooth. Injectivity is obvious. We now show that $L (M _\rho)=
\mathcal{L}_{\pi^{[\mu]} (m)}$. Let $\pi _\rho (z) \in M _\rho$. By the definition
of the optimal momentum map there exists a finite composition of Hamiltonian flows
$F_{t _1}, \ldots, F_{t _n }$ corresponding to
$G$-invariant Hamiltonian functions such that  $z=(F^1_{t _1}\circ
\cdots \circ  F^n_{t _n })(m)$. For simplicity in the exposition take
$n=1$ and let $F \in  C^\infty(M)^{G} $ be the function whose Hamiltonian flow is
$F _t $. Let $f \in C^{\infty}\left(M ^{[\mu]}\right)$ be the function
defined by $F \circ  i ^{[\mu]}= f\circ  \pi ^{[\mu]} $. It is easy to
see that the Hamiltonian flow $F _t ^f$ of 
$X _f $ in $M ^{[\mu]} $ is such that  $F _t ^f \circ \pi  ^{[\mu]}= \pi^{[\mu]}
\circ  F _t  \circ  i ^{[\mu]} $. Consequently
\[ L(\pi _\rho(z))=L(\pi _\rho(F _t(m)))= \pi  ^{[\mu]}(F _t  (m))=F _t ^f(\pi
^{[\mu]} (m) )\in 
\mathcal{L}_{\pi^{[\mu]} (m)}.
\] This argument can be reversed by using local $G$-invariant extensions of the
compositions of the functions in $ M ^{[\mu]}  $ with $ \pi ^{[\mu]}$ (that do
exist by part {\bf (i)}) and  hence proving the converse inclusion
$\mathcal{L}_{\pi^{[\mu]} (m)} \subset L(M _\rho)$.

The closedness hypothesis on $\mathcal{J}^{-1}(\rho)  $ implies that this set is
an embedded submanifold of  $M$ (recall that closed integral leaves of constant
rank integrable distributions are always embedded~\cite{camacho neto}).
Additionally, it can be proved (see Proposition 4.5 in~\cite{optimal}) that under
this hypothesis the isotropy subgroup $G_{\rho} $ is closed in $G$ and that
\begin{equation}
\label{closedness and its miracles}
\mathfrak{g}_\rho \cdot  z= T _z \mathcal{J}^{-1}(\rho)\cap  \mathfrak{g}\cdot  z,
\quad z  \in
\mathcal{J}^{-1}(\rho). 
\end{equation} Consequently, the form of the bracket~(\ref{Poisson
structure reduced optimal momentum map})  is obtained by mimicking  the
proof of part {\bf  (iii)} of Theorem~\ref{cylinder valued momentum map
symplectic dynamics reduction}. Notice that the local $G$-invariant
extensions needed in~(\ref{Poisson structure reduced optimal momentum
map}) always exist by Lemma 10.4.14 in~\cite{hsr},  which can be applied
due to the fact that $\mathcal{J}^{-1}(\rho) 
$ is an embedded submanifold of
$M$  and~(\ref{closedness and its miracles}) holds. The Poisson character of $L$ 
is a straightforward consequence of~(\ref{Poisson structure reduced cylinder
valued momentum map}) and~(\ref{Poisson structure reduced optimal momentum map}).
\quad
$\blacksquare$

\section{Example: Magnetic cotangent bundles of Lie groups}
\label{Example: Magnetic cotangent bundles of Lie groups section}

Let $G$  be a finite dimensional Lie group and $T ^\ast  G $  its cotangent bundle 
endowed with the magnetic symplectic structure $\overline{\omega _\Sigma}:= \omega
_{{\rm can}}- \pi ^\ast  B _\Sigma  $, where $\omega _{{\rm can}}$  is the
canonical symplectic form  on $T ^\ast G $,
$\pi :T ^\ast  G \rightarrow  G $  is the projection onto the base, and $B _\Sigma
\in  \Omega ^2(G)^G $ is a left invariant  two-form on $G$  whose value at the
identity is the Lie algebra two-cocycle $\Sigma: \mathfrak{g} \times   \mathfrak{g}
\rightarrow \mathbb{R}$. Since $\Sigma $  is a cocycle it follows that $B _\Sigma 
$ is closed and hence $\overline{\omega _\Sigma}$ is a symplectic  form.

The cotangent lift of the action of $G$ on itself by left translations produces,
due to the invariance of $B _\Sigma $, a canonical
$G$-action on $(T ^\ast  G, \overline{\omega _\Sigma})  $. In the absence of
magnetic terms this action has an associated coadjoint equivariant
momentum map $\mathbf{J}:T ^\ast  G \rightarrow   \mathfrak{g}^\ast$
given by $\langle \mathbf{J}(\alpha_g ),\xi\rangle:= \langle \alpha _g
,T_e R _g (\xi)\rangle $, $ \alpha _g \in T ^\ast  G $, $\xi
\in  \mathfrak{g}$. It is well known that the Marsden-Weinstein reduced spaces
associated to this action are naturally symplectomorphic to the coadjoint orbits
of the $G$-action on $\mathfrak{g}^\ast$, endowed with their canonical
Kostant-Kirillov-Souriau orbit symplectic form. The magnetic term
destroys this picture in most cases. As we will see later on in this
section, $(T^\ast  G,\overline{\omega _\Sigma})$ has, in general, a
non-zero Hamiltonian holonomy ${\mathcal H}$ and hence the lift of left
translation does not admit anymore a standard momentum map. This forces
us, when carrying out reduction, to work in the degree of generality
of the preceding section. We will hence compute in this
setup the cylinder valued and the optimal momentum maps and will
characterize the three reduced spaces introduced in
theorems~\ref{symplectic reduction theorem cylinder valued momentum map}
and~\ref{all reductions cylinder valued momentum map}.  One of the
conclusions of our discussion will be the fact that the resulting reduced
spaces are related not to the coadjoint orbits  but to the orbits of the
extended affine action of $ G_{\Sigma}
$ on $\mathfrak{g}^\ast$, with $ G_{\Sigma} $  the connected and simply
connected  Lie group that integrates $\mathfrak{g}_\Sigma 
$, the one-dimensional central extension of
$\mathfrak{g}$  constructed using the cocycle $\Sigma $. More
specifically, the optimal reduced spaces are naturally symplectomorphic
to these orbits which shows that, in this context, optimal reduction is
the natural generalization of the picture that allows one to see in the
standard setup the Kostant-Kirillov-Souriau coadjoint orbits as
symplectic reduced spaces.

In order to make the problem more tractable we introduce the left trivialization
$\lambda : T ^\ast G \rightarrow  G\times  
\mathfrak{g} ^{\ast}$   of  $T ^\ast   G $ given by  $\lambda(\alpha _g):=(g, T _e
^\ast  L _g(\alpha _g))$, for any $\alpha _g \in T ^\ast  _g G $. Its inverse
$\lambda ^{-1} :G
\times 
\mathfrak{g}^\ast \rightarrow  T ^\ast G $  is  $ \lambda ^{-1}(g, \mu )=T ^\ast 
_g L_{g ^{-1}}(\mu) $, $(g, \mu)\in  G \times \mathfrak{g}^\ast$. The cotangent
lift of left translations on $G$  to  $T ^\ast  G $   is given in this
trivialization  by $h \cdot  (g, \mu )=(h g, \mu )$, $ h,g \in  G $, $\mu \in 
\mathfrak{g}^\ast$, and hence the infinitesimal generators take the form 
\begin{equation}
\label{infinitesimal generators of trivialized}
\xi_{G \times  \mathfrak{g}^\ast}(g, \mu)=(T _e R _g (\xi), 0)=(T _e L
_g(\mbox{\rm Ad}_{g ^{-1}}\xi), 0), \quad \xi \in  \mathfrak{g}.
\end{equation} Additionally, the magnetic symplectic form $\omega_\Sigma:=(\lambda
^{-1})^\ast 
\overline{\omega _\Sigma} $ on $G \times  \mathfrak{g}^\ast$ has the expression
\begin{equation}
\label{trivialized magnetic symplectic form}
\omega _\Sigma (g,  \nu)\left((T _e L _g (\xi), \rho),  (T_e L_g(\eta),
\sigma)\right)
=\langle \sigma, \xi\rangle- \langle  \rho, \eta\rangle+\langle
\nu,[\xi, \eta]\rangle- \Sigma (\xi, \eta).
\end{equation} 
In order to compute the cylinder valued momentum map of the
canonical $G$-action on $(G \times 
\mathfrak{g}^\ast, \omega _\Sigma )$ we note that the horizontal distribution
associated to the connection $\alpha$  on the bundle $(G \times 
\mathfrak{g}^\ast)\times  \mathfrak{g}^\ast \rightarrow  G \times 
\mathfrak{g}^\ast$  defined in~(\ref{definition of alpha connection}) equals,
by~(\ref{infinitesimal generators of trivialized}) and~(\ref{trivialized magnetic
symplectic form}), 
\begin{equation}
\label{horizontal distribution magnetic} H((g, \mu), \nu)=\{(T _e L _g (\xi), \rho
, \mbox{\rm Ad} ^\ast _{g ^{-1}}(\rho-\mbox{\rm ad} ^\ast _\xi \mu+ \Sigma(\xi,
\cdot )) \mid \xi \in 
\mathfrak{g}, \rho \in \mathfrak{g}^\ast
\}.
\end{equation}

\begin{proposition}
\label{intermediated connection magnetic} Let $G \times  \mathfrak{g}^\ast 
\rightarrow  G$ be the trivial principal fiber bundle with structure group
$(\mathfrak{g}^\ast, +) $ and where the action is given by $R _\nu(g, \mu):=(g,
\mu - \nu)$, $g  \in G $, $\mu, \nu \in  \mathfrak{g}^\ast$.
\begin{itemize}
\item [{\rm {\bf (i)}}] The distribution $H _G $ on $G \times 
\mathfrak{g}^\ast$ given by 
\[ H _G(g, \mu)=\{(T _e L _g (\xi), \mbox{\rm Ad}^\ast _{g ^{-1}}(\Sigma(\xi,
\cdot ))) \mid \xi
\in \mathfrak{g}
\}
\] defines a flat connection $\alpha _G  $ on $G \times  \mathfrak{g}^\ast
\rightarrow G $.
\item [{\rm {\bf (ii)}}] The holonomy bundle $\widetilde{(G \times 
\mathfrak{g}^\ast)}_{(g,
\mu, \nu)}$ of the connection $\alpha$   defined in~(\ref{definition of alpha
connection}) that contains the point $(g, \mu, \nu) \in  (G
\times  \mathfrak{g}^\ast)\times  \mathfrak{g}^\ast $ is given by 
\[
\widetilde{(G \times  \mathfrak{g}^\ast)}_{(g, \mu, \nu)}=\{(h, \rho, \mbox{\rm
Ad} ^\ast _{h ^{-1}}\rho + \tau)\mid \rho \in  \mathfrak{g}^\ast, (h , \tau)\in 
\widetilde{G}_{\left(g, \nu-
\mbox{\rm Ad} ^\ast _{g ^{-1}}\mu\right)} 
\},
\] where $\widetilde{G}_{\left(g, \nu-
\mbox{\rm Ad} ^\ast _{g ^{-1}}\mu\right)} $ is the holonomy bundle of the
connection $\alpha _G  $  containing the point $(g, \nu-
\mbox{\rm Ad} ^\ast _{g ^{-1}}\mu)$.
\item [{\rm {\bf (iii)}}] The holonomy group ${\mathcal H}_G $ of 
$\alpha _G  $  and the Hamiltonian  holonomy $ {\mathcal H}  $ of the
$G$-action  on  $(G \times 
\mathfrak{g}^\ast, \omega _\Sigma)$ coincide, that is, ${\mathcal H}= {\mathcal H}
_G$.
\end{itemize}
\end{proposition}

\noindent\textbf{Proof.\ \ (i)} The vertical bundle $V_G$  of  $G \times 
\mathfrak{g}^\ast
\rightarrow  G$ is given by  $V_G(g,  \mu)=\{(0, \rho)\mid \rho \in
\mathfrak{g}^\ast\} $. It easily follows that $T_{(g, \mu)}(G \times 
\mathfrak{g}^\ast)=H _G (g, \mu)\oplus V _G (g, \mu)$ and that $T_{(g, \mu)}R _\nu
(H _G(g, \mu))=H_G(g, \mu -\nu)$, for any $g \in  G $, $\mu,
\nu\in\mathfrak{g}^\ast$, which proves that $H _G $ is the horizontal bundle of a
connection on $G \times  \mathfrak{g}^\ast \rightarrow  G $ whose associated
one-form $\alpha _G \in  \Omega ^1(G \times  \mathfrak{g}^\ast;
\mathfrak{g}^\ast)$ is
\[
\alpha _G(g, \mu)\left(T _e L _g (\xi), \rho \right)= \mbox{\rm Ad}^\ast _{g
^{-1}}\left(
\Sigma(\xi, \cdot )\right)- \rho.
\] In the computation of this expression we used that the horizontal $v_{(g,
\mu)}^H $ and vertical
$v_{(g, \mu)}^V $ components of any vector $v_{(g, \mu)}=(T _e L _g (\xi), \rho )
$ are given by 
$v_{(g, \mu)}^H =(T _e L _g (\xi), \mbox{\rm Ad}^\ast _{g ^{-1}}(\Sigma(\xi, \cdot
)))$ and 
$v_{(g, \mu)}^V= (0, \rho- \mbox{\rm Ad}^\ast _{g ^{-1}}(\Sigma(\xi, \cdot ))) $.
The flatness of $\alpha_G $ will be obtained as a consequence of point {\bf (iii)}
and of the discrete character of ${\mathcal H} $ (see below).

\smallskip

\noindent {\bf (ii)} Let $(h, \tau) \in  \widetilde{G}_{\left(g, \nu-
\mbox{\rm Ad} ^\ast _{g ^{-1}}\mu\right)} $. By definition, there exists
a piecewise smooth horizontal curve $(g (t) , \tau (t) )$ such that  $g
(0) =g $, $g (1)=h $, $\tau(0)= \nu-
\mbox{\rm Ad}^\ast _{g ^{-1}}\mu $, $\tau (1)= \tau$, and $ \tau ' (t)= \mbox{\rm
Ad}^\ast _{g (t)^{-1}}\left(\Sigma(T_{g (t)}L_{g (t)^{-1}}(g' (t)), \cdot )
\right)$. The piecewise smooth curve $\gamma (t):=(g (t) , \mu (t), \mbox{\rm
Ad}^\ast _{g (t)^{-1}}\mu (t)+ \tau (t)) $, with $\mu(t) $ an arbitrary smooth
curve in $\mathfrak{g}^\ast$ such that $\mu (0)= \mu $  and  $\mu (1) = \rho $,
satisfies $\gamma (0)=(g, \mu ,  \nu)$,
$\gamma  (1)=(h, \rho, \mbox{\rm Ad}^\ast _{h ^{-1}}\rho+ \tau)$, and is
horizontal. Indeed,  let
$\xi(t):=T_{g (t)}L_{g (t)^{-1}}(g' (t))\in  \mathfrak{g}$. Then
\begin{align*}
\gamma ' (t)&= \left(T _eL_{g (t)}(\xi (t)), \mu' (t), \mbox{\rm Ad}^\ast _{g (t)
^{-1}}\mu ' (t)- \mbox{\rm Ad}^\ast _{g (t) ^{-1}}\left(\mbox{\rm ad}^\ast _{\xi
(t)}\mu (t)\right)+ \tau' (t)\right)\\
	&=\left(T _eL_{g (t)}(\xi (t)), \mu' (t), \mbox{\rm Ad}^\ast _{g (t)
^{-1}}\left(\mu ' (t)- \mbox{\rm ad}^\ast _{\xi (t)}\mu (t)+ \Sigma(\xi (t), \cdot
)\right)\right),
\end{align*} which belongs to $H(\gamma (t) ) $ by~(\ref{horizontal distribution
magnetic}). This proves that $ (h, \rho, \mbox{\rm Ad}^\ast _{h ^{-1}}\rho+ \tau)
\in  \widetilde{(G \times 
\mathfrak{g}^\ast)}_{(g, \mu, \nu)} $.

Conversely, let $(h, \rho , \sigma) \in  \widetilde{(G \times 
\mathfrak{g}^\ast)}_{(g, \mu, \nu)} $. By definition, there exists a
piecewise smooth horizontal curve $\gamma (t):=(g (t), \mu (t), \nu
(t))$ such that $\gamma (0)=(g, \mu , \nu)$, $\gamma (1)=(h, \rho ,
\sigma)$, and 
\begin{equation}
\label{nu prime equals}
\nu' (t)= 
\mbox{\rm Ad}^\ast _{g (t) ^{-1}}(\mu ' (t)- \mbox{\rm ad}^\ast _{\xi (t)}\mu (t)+
\Sigma(\xi (t), \cdot )),
\end{equation} where
$\xi(t)=T_{g (t)}L_{g (t)^{-1}}(g' (t))$. We now show that $(h, \rho , \sigma)=(h,
\rho, \mbox{\rm Ad}^\ast  _{h ^{-1}}\rho + \tau)$ where $(h, \tau )\in  
\widetilde{G}_{\left(g, \nu-
\mbox{\rm Ad} ^\ast _{g ^{-1}}\mu\right)} $. Let $\tau (t):= \nu (t)- \mbox{\rm
Ad}^\ast _{g (t) ^{-1}}\mu (t) $. Notice that $\tau (0)= \nu - \mbox{\rm Ad} ^\ast
_{g ^{-1}}\mu$, $\tau (1)= \sigma- \mbox{\rm Ad}^\ast _{h ^{-1}} \rho $, and
that, by~(\ref{nu prime equals}), $\tau' (t)= \mbox{\rm Ad}^\ast _{g(t)
^{-1}}\left(
\Sigma(\xi(t), \cdot )\right)$, which shows that $(g (t) , \tau (t))$ is
horizontal. Consequently $(h, \rho , \sigma)=(g (1), \rho, \mbox{\rm Ad} ^\ast _{h
^{-1}}\rho + \tau (1))$  and hence the claim follows.

\smallskip

\noindent {\bf (iii)} Let  $\mu \in  {\mathcal H}  $  be arbitrary. Then there
exists a loop $(g (t) , \rho (t))$ in $G \times  \mathfrak{g}^\ast$ whose
horizontal lift $(g (t), \rho (t), \mu (t))$ satisfies that $\mu= \mu (1)- \mu
(0)$. Consequently, using horizontality
\begin{align}
\label{thing with horizontal integral}
\mu&=\int _0 ^1 \mu ' (t) dt= \int _0^1 \mbox{\rm Ad}^\ast _{g (t) ^{-1}}
\left(\rho ' (t)- \mbox{\rm ad}^\ast _{\xi (t)}\rho (t)+ \Sigma(\xi (t), \cdot
)\right) dt \nonumber\\ 
&= \int_0 ^1 \left(\frac{d}{dt}\left(\mbox{\rm Ad}^\ast _{g  (t)
^{-1}}\rho (t)
\right)+\mbox{\rm Ad}^\ast _{g  (t) ^{-1}} \Sigma(\xi (t), \cdot
)\right) d t=
\int_0 ^1\mbox{\rm Ad}^\ast _{g  (t) ^{-1}} \Sigma(\xi (t), \cdot ) d t,
\end{align} 
where $\xi(t)=T_{g (t)}L_{g (t)^{-1}}(g' (t))$. Therefore the loop
$g (t) $ in $G$ has a horizontal lift $(g (t), \nu (t))$ such that 
\[
\nu (1)- \nu (0)=\int _0 ^1 \nu' (t) dt=\int_0 ^1\mbox{\rm Ad}^\ast _{g  (t) ^{-1}}
\Sigma(\xi (t), \cdot ) d t= \mu,
\] which proves that $\mu \in  {\mathcal H}_G $. The converse inclusion is
obtained by reading backwards the previous argument.

Finally, the equality ${\mathcal H}= {\mathcal H}_G $ implies that the connection
$\alpha _G$ is flat since the tangent spaces to the holonomy bundles equal the
horizontal distribution because ${\rm Lie}\left( {\mathcal H} \right)=\{0\} $.
Indeed, as the distribution associated to the holonomy bundles is integrable by
general theory this shows that the horizontal distribution is integrable and hence
the associated connection is flat. \quad $\blacksquare$

\begin{remark}
\normalfont If the Lie algebra two-cocycle $\Sigma$ can be integrated to a smooth
$\mathfrak{g}^\ast$-valued group one-cocycle $\sigma:G \rightarrow 
\mathfrak{g}^\ast$, that is 
\[
\Sigma(\xi, \cdot )=T _e \sigma (\xi), \qquad \text{for any } \xi \in 
\mathfrak{g},
\] then the holonomy bundles $\widetilde{G}_{(g, \mu)}$ are the graphs of $\sigma$.
More specifically
\[
\widetilde{G}_{(g, \mu)}=\{(h, \sigma(h)+ \mu- \sigma  (g))\mid h \in  G \}.
\] In this particular case ${\mathcal H}_G=\{0\} $  and hence, by 
Proposition~\ref{intermediated connection magnetic},
${\mathcal H}=\{0\}$.
\end{remark} 

We are now going to use the central extensions and their actions introduced in
Section~\ref{Poisson structures on and extensions} to better characterize the
holonomy bundles of $\alpha $ and $\alpha _G $. This will allow us to give an
explicit expression for the cylinder valued momentum map of the $G$-action on the
magnetic cotangent bundle $(G \times  \mathfrak{g}^\ast, \omega _\Sigma)$.

\begin{proposition}
\label{characterization with extension holonomy} Let  $\mathfrak{g}_\Sigma $ be
the one-dimensional central extension of the Lie algebra $\mathfrak{g} $
determined by the cocycle $\Sigma $, $G_{\Sigma} $  the connected and simply
connected Lie group that integrates it, and $ \mu_{\Sigma}: G_{\Sigma}
\rightarrow  \mathfrak{g}^\ast$ the extended $\mathfrak{g}^\ast$-valued one
cocycle introduced in Proposition~\ref{properties of mu and more}. The holonomy
bundle $\widetilde{G}_{(h, \nu)} $ of the connection $\alpha _G  $ that contains
the point $(h, \nu) \in  G \times  \mathfrak{g}^\ast$ equals
\begin{equation}
\label{expression of  bundle with sigma}
\widetilde{G}_{(h, \nu)}=\left\{(\pi_G(\widehat{g}), \mu_{\Sigma}
(\widehat{g}^{-1})-\mu_{\Sigma} (\widehat{h}^{-1})+ \nu)\mid \widehat{g},
\widehat{h} \in  G_{\Sigma}\text{ such that } \pi _G (\widehat{h})=h
\right\},
\end{equation} with $\pi _G : G_{\Sigma} \rightarrow  G $ the projection.
\end{proposition}

\noindent\textbf{Proof.\ \ }Since $\widetilde{G}_{(h, \nu)}=R_{-
\nu}(\widetilde{G}_{(h, 0)})$ it suffices to show that 
\begin{equation}
\label{thing suffices to prove sigma}
\widetilde{G}_{(h, 0)}=\left\{(\pi_G(\widehat{g}), \mu_{\Sigma}
(\widehat{g}^{-1})-\mu_{\Sigma} (\widehat{h}^{-1}))\mid \widehat{g}, \widehat{h}
\in  G_{\Sigma}\text{ such that } \pi _G (\widehat{h})=h
\right\}.
\end{equation} 
We begin with the  inclusion $\supset $. Let $\widehat{g} \in G_{\Sigma}$
and $\widehat{g} (t) $ a piecewise smooth curve in $G_{\Sigma} $ such
that  $\widehat{g}(0)=\widehat{h} $, $\widehat{g} (1)= \widehat{g}$, and
where $\pi _G(\widehat{h})=h  $. We will show that the element $(\pi_G
(\widehat{g}), \mu_{\Sigma}(\widehat{g}^{-1})-
\mu_{\Sigma}(\widehat{h}^{-1}))$ belongs to $\widetilde{G}_{(h,0)}$ by
proving that the curve $\gamma(t):=(\pi_G (\widehat{g}(t)),
\mu_{\Sigma}(\widehat{g}(t)^{-1})- \mu_{\Sigma}(\widehat{h}^{-1}))$ is 
horizontal and connects $(h,0)$ with $(\pi_G (\widehat{g}),
\mu_{\Sigma}(\widehat{g}^{-1})- \mu_{\Sigma}(\widehat{h}^{-1}))$. Indeed,
$\gamma(0)=(h,0)$,  $\gamma (1)=(\pi_G (\widehat{g}),
\mu_{\Sigma}(\widehat{g}^{-1})- \mu_{\Sigma}(\widehat{h}^{-1}))$, and 
$\gamma (t)$ is horizontal because if we write
\[
\frac{d}{dt} \pi _G (\widehat{g} (t))=T _e L_{g (t)}(\xi(t)),
\] 
where $g (t):= \pi _G(\widehat{g} (t))$ and $\xi(t):=T_{g (t)}L_{g
(t)^{-1}}(g' (t))$, then by Proposition~\ref{properties of mu and more} 
\begin{equation}
\label{derivative of mu sigma}
\frac{d}{dt} \mu_{\Sigma}(\widehat{g}(t) ^{-1})=- \frac{d}{dt} \left(\mbox{\rm Ad}
^\ast _{g (t)^{-1}}\mu_{\Sigma}(\widehat{g }(t))\right)=- \mbox{\rm Ad}^\ast _{g
(t)^{-1}}
\left(  - \mbox{\rm ad}^\ast _{\xi (t)}\mu_{\Sigma}
(\widehat{g}(t))+T_{\widehat{g}(t)}\mu_{\Sigma}(\widehat{g}' (t))
\right).
\end{equation} Now let $(\widehat{\xi} (t), s
(t)):=T_{\widehat{g}(t)}L_{\widehat{g}(t)^{-1}}(\widehat{g}' (t))$. It turns out
that $\widehat{\xi} (t)= \xi (t) $ because
\begin{multline*}
\widehat{\xi} (t)= \pi _{\mathfrak{g}} (\widehat{\xi} (t), s (t))=T _e \pi_G
(\widehat{\xi} (t), s (t))= \left.\frac{d}{d\epsilon}\right|_{\epsilon=0}\pi
_G(\widehat{g} (t) ^{-1} \widehat{g}(t + \epsilon ))\\
=\left.\frac{d}{d\epsilon}\right|_{\epsilon=0}\pi _G(\widehat{g} (t) ^{-1})
\pi  _G( \widehat{g}(t + \epsilon ))=\left.\frac{d}{d\epsilon}\right|_{\epsilon=0}
g  (t) ^{-1} g(t+ \epsilon)=T_{g (t)}L_{g (t) ^{-1}}(g' (t))= \xi (t),
\end{multline*} and hence~(\ref{derivative of mu sigma}) equals, by
Proposition~\ref{properties of mu and more}, 
\[
\frac{d}{dt} \mu_{\Sigma}(\widehat{g}(t) ^{-1})=- \mbox{\rm Ad}^\ast _{g (t)^{-1}}
\left(  - \mbox{\rm ad}^\ast _{\xi (t)}\mu_{\Sigma} (\widehat{g}(t))+\mbox{\rm
ad}^\ast _{\xi (t)}\mu_{\Sigma} (\widehat{g}(t))- \Sigma(\xi (t), \cdot )
\right)=\mbox{\rm Ad}^\ast _{g (t)^{-1}}
\left( \Sigma(\xi (t), \cdot )
\right),
\] which shows that $\gamma (t)  $ is horizontal, as required.

We prove next the reverse inclusion 
\[
\widetilde{G}_{(h,0)}\subset \left\{(\pi_G(\widehat{g}), \mu_{\Sigma}
(\widehat{g}^{-1})-\mu_{\Sigma} (\widehat{h}^{-1}))\mid \widehat{g}, 
\widehat{h} \in  G_{\Sigma}\text{ such that } \pi _G (\widehat{h})=h
\right\}
\] 
by showing that any piecewise smooth horizontal curve $\gamma (t) \subset
\widetilde{G}_{(h,0)}$ such that  $\gamma (0)=(h,0)$ satisfies 
$\gamma(t) \in  \left\{(\pi_G(\widehat{g}), \mu_{\Sigma}
(\widehat{g}^{-1})-\mu_{\Sigma} (\widehat{h}^{-1}))\mid \widehat{g},
\widehat{h} \in  G_{\Sigma}\text{ such that } \pi _G (\widehat{h})=h
\right\}$ for all $t \in  I $, where $I$ is the time interval on which
$\gamma$  is defined. Let  $g (t) $ and  $\mu (t) $, $t \in  I $, be two
curves in $G$ and $\mathfrak{g}^\ast$, respectively, such that  $\gamma
(t)=(g (t) , \mu (t))$. The horizontality of $\gamma (t) $ implies that
$\mu ' (t)= \mbox{\rm Ad}^\ast _{g  (t) ^{-1}}\left( \Sigma(\xi (t),
\cdot )\right) $, with $\xi(t)=T_{g (t)}L_{g (t)^{-1}}(g' (t))$. Since
the map $ \pi_G: G_{\Sigma} \rightarrow  G $ is a surjective submersion
it admits local sections. In particular, there exists an open
neighborhood $U  $ of  $h$ in $ G $  and a map $\sigma:U \subset G
\rightarrow  G_{\Sigma}$ such that  $\pi_G \circ \sigma= {\rm id }| _U $  and
$\sigma (h)= \widehat{h} $. The smoothness of the curve  $g (t) $ implies that
there exists $0< t _0 \in  I $ such that $g (t) \in  U $ for any $t \in [0,t _0]$.
Let $\widehat{g} (t):= \sigma (g (t))$ and $\nu (t):=
\mu_{\Sigma}(\widehat{g}(t)^{-1})- \mu_{\Sigma}(\widehat{h}^{-1})$. We will show
that the curve $\gamma _\sigma (t):=(g (t), \nu (t))= (\pi _G(\widehat{g} (t)),
\mu_{\Sigma}(\widehat{g}(t) ^{-1})- \mu_{\Sigma}(\widehat{h}^{-1}))\subset
\left\{(\pi_G(\widehat{g}), \mu_{\Sigma} (\widehat{g}^{-1})-\mu_{\Sigma}
(\widehat{h}^{-1}))\mid \widehat{g}, \widehat{h}
\in  G_{\Sigma}\text{ such that } \pi _G (\widehat{h})=h
\right\}$ is such that
$\gamma _\sigma (t)= \gamma (t) $, for any $t \in [0, t _0]$. Notice that since
$\gamma _\sigma (0)=(g (0), \nu (0))=(h, \mu_{\Sigma}(\sigma(g (0))^{-1})-
\mu_{\Sigma}(\widehat{h}^{-1}))=(h,\mu_{\Sigma}(\widehat{h}^{-1})-\mu_{\Sigma}(\widehat{h}^{-1}))=(h,0)$,
the uniqueness of horizontal lifts guarantees that it suffices to check that
$\gamma_\sigma (t) $ is horizontal. Given that $\pi_G(\widehat{g} (t))=g (t)$ for
any $t \in  [0, t _0]$, an argument similar to the one in the first part of the
proof shows that 
\begin{equation}
\label{argument similar horizontal}
\widehat{g}' (t)=T_{\widehat{g}(t)}L_{\widehat{g}(t) ^{-1}} \left(\xi(t), s (t)
\right),
\end{equation} where $\xi(t)=T_{g (t)}L_{g (t)^{-1}}(g' (t))$ and $s (t) $ is some
piecewise smooth curve in ${\Bbb R}$. The use of~(\ref{argument similar
horizontal}) and of  Proposition~\ref{properties of mu and more} in a
straightforward computation show that $\nu' (t)= \mbox{\rm Ad}^\ast _{g
(t)^{-1}}(\Sigma(\xi(t), \cdot )) $, which proves that $\gamma _\sigma (t) $ is
horizontal and hence $\gamma _\sigma (t)= \gamma (t) \subset
\left\{(\pi_G(\widehat{g}),
\mu_{\Sigma} (\widehat{g}^{-1})-\mu_{\Sigma} (\widehat{h}^{-1}))\mid \widehat{g},
\widehat{h}
\in  G_{\Sigma}, \pi _G (\widehat{h})=h
\right\}$, for any $t \in [0, t _0]$, as required. We can repeat what we just did
by taking a local section $\sigma_1 $ of $\pi_G $ around $g (t _0)$ which will
allow us to construct a curve
$\gamma_{\sigma _1}(t)\subset
\{(\pi_G(\widehat{g}),
\mu_{\Sigma} (\widehat{g}^{-1})-\mu_{\Sigma} (\widehat{h}^{-1}))\mid \widehat{g},
\widehat{h}
\in  G_{\Sigma}, \pi _G (\widehat{h})=h\}$, $t \in [t _0, t _1]$, $t _1> t _0$,
such that  $\gamma _{\sigma _1}(t)= \gamma (t)$, for any $t \in [t _0, t _1]$. The
compactness of $I$ guarantees that by repeating this procedure a finite number of
times we can write $\gamma$ as a broken path made of finite smooth curves included
in $\{(\pi_G(\widehat{g}),
\mu_{\Sigma} (\widehat{g}^{-1})-\mu_{\Sigma} (\widehat{h}^{-1}))\mid \widehat{g},
\widehat{h}
\in  G_{\Sigma}, \pi _G (\widehat{h})=h\}$ which proves that $\gamma$ itself is
included in $\{(\pi_G(\widehat{g}),
\mu_{\Sigma} (\widehat{g}^{-1})-\mu_{\Sigma} (\widehat{h}^{-1}))\mid \widehat{g},
\widehat{h}
\in  G_{\Sigma}, \pi _G (\widehat{h})=h\}$, hence proving the desired inclusion.
\quad $\blacksquare$ 

\begin{corollary}
\label{inclusion with holonomy relevant abelian} In the setup of
Proposition~\ref{characterization with extension holonomy} the following inclusion
holds
\[ {\mathcal H} \subset \mu_{\Sigma}(G_{\Sigma}).
\]
\end{corollary}

\noindent\textbf{Proof.\ \ }Since ${\mathcal H}= {\mathcal H}_G $  by
Proposition~\ref{intermediated connection magnetic} it suffices to show 
that ${\mathcal H}_G \subset \mu_{\Sigma}(G_{\Sigma})$. Let  $\nu \in 
{\mathcal H}_G $. By definition, there exists a loop $g (t)  $
in $G$  such that  $g (0)=g (1)= e
$ with a horizontal lift $ \gamma (t)=(g (t), \mu (t))$ that satisfies 
$\gamma (0)=(e,0)$ and $\gamma (1)=(e, \nu)$. The compactness of the
interval $[0,1]$ implies that we can take local sections
$\sigma_1,\ldots,\sigma _n $ of the projection $\pi_G: G_{\Sigma}\rightarrow G $
and that we can split the interval $[0,1]$ into $n$ intervals of the form 
$[t_{i-1}, t _i] $ with $0= t _0< t _1 <\ldots <t _n=1 $ such that  for any $i \in
\{1, \ldots,n \} $ we can define $\widehat{g}_i(t):= \sigma_i(g (t))$, $t \in 
[t_{i-1}, t _i]$. The sections $\sigma _1, \ldots, \sigma _n $ are chosen in such
a way that $\widehat{g}(0)=e $ and $\widehat{g}_i(t _i)= \widehat{g}_{i+1}(t _i)$.
Moreover, by construction, $g (t)= \pi _G(\widehat{g}_i(t))$ for any $t
\in [t_{i-1}, t _i]$ and hence a strategy similar to the one in the
first part of the proof of Proposition~\ref{characterization with
extension holonomy} shows that
$\gamma_1(t):=(\pi_G(\widehat{g}_i(t)), \mu_{\Sigma}(\widehat{g}_i(t)^{-1}))  $ is
a piecewise smooth horizontal curve for $\alpha_G$ such that $\gamma_1 (0)=(e,0)$
and hence $\gamma _1 (t)= \gamma (t) $, for any $t \in [0,1]$. Consequently, $\nu
= \mu_{\Sigma}(\widehat{g}_n(1))\in  \mu_{\Sigma}(G_{\Sigma})$, as required. \quad
$\blacksquare$

\begin{theorem}
\label{explicit expression of cylinder valued momentum map magnetic} Let $(G \times
\mathfrak{g}^\ast, \omega _\Sigma)$ be a magnetic cotangent bundle of the Lie
group $G$. The map $\mathbf{K}: G \times  \mathfrak{g}^\ast \rightarrow
\mathfrak{g}^\ast/
\overline{{\cal H}} $ given by the expression
\begin{equation}
\label{cylinder valued momentum map for magnetic period}
\mathbf{K}(g, \eta)= \pi_C\left (\mbox{\rm Ad}^\ast _{g ^{-1}}\eta+
\mu_{\Sigma}(\widehat{g} ^{-1})+ \nu _0\right)= \Xi(\widehat{g}, \eta+
\overline{{\cal H}})+ \pi _C(\nu _0)
\end{equation} is a cylinder valued momentum map for the canonical $G$-action  on
$(G \times
\mathfrak{g}^\ast, \omega _\Sigma)$. The element $\nu_0 \in  \mathfrak{g}^\ast$ is
an arbitrary constant, $(g, \eta) \in G \times  \mathfrak{g}^\ast$ is arbitrary, and
$\widehat{g}
\in G_{\Sigma} $ is any element such that  $\pi_G (\widehat{g})=g $. The map
$\mu_{\Sigma}: G_{\Sigma}\rightarrow  \mathfrak{g}^\ast$ is the extended
$\mathfrak{g}^\ast$-valued one-cocycle associated to $\Sigma$ and $\Xi: G_{\Sigma}
\times  \mathfrak{g}^\ast / \overline{{\cal H}} \rightarrow \mathfrak{g}^\ast /
\overline{{\cal H}}$ the associated $G_{\Sigma}$-action on $\mathfrak{g}^\ast /
\overline{{\cal H}}$.

The non-equivariance cocycle $\sigma:G \rightarrow \mathfrak{g}^\ast/
\overline{{\cal H}}$ of ${\bf K}$ is given by 
\begin{equation}
\label{non equivariance cylinder valued momentum map magnetic}
\sigma(g)= \pi _C \left( \mu_{\Sigma}(\widehat{g}^{-1})+ \nu _0- \mbox{\rm
Ad}^\ast _{g ^{-1}}(\nu _0)\right),
\end{equation} with $g \in  G $ and $\widehat{g} \in G_{\Sigma} $ such that 
$\pi_G(\widehat{g})=g
$. Finally
\begin{equation}
\label{inclusion holonomy kernel}
\mu_{\Sigma}(\ker \pi _G)\subset  \overline{{\cal H}}.
\end{equation}
\end{theorem}

\noindent\textbf{Proof.\ \ }If we put together the conclusions of
propositions~\ref{intermediated connection magnetic} and~\ref{characterization with
extension holonomy} we can conclude that for any $(h, \mu , \nu)\in G \times 
\mathfrak{g}^\ast\times  \mathfrak{g}^\ast$
\[
\widetilde{(G \times  \mathfrak{g}^\ast)}_{(h, \mu , \nu)}=\{ (\pi _G
(\widehat{g}), \eta, \mbox{\rm Ad}^\ast _{\pi _G(\widehat{g}^{-1})}\eta+
\mu_{\Sigma}(\widehat{g}^{-1})-\mu_{\Sigma}(\widehat{h}^{-1})+ \nu - \mbox{\rm Ad}
^\ast _{h ^{-1}}\mu)\mid \widehat{g} \in  G_{\Sigma}, \eta \in  \mathfrak{g}^\ast
\}.
\] Hence setting $\nu _0:=\nu - \mbox{\rm Ad} ^\ast _{h ^{-1}}\mu
-\mu_{\Sigma}(\widehat{h}^{-1})$ and  using  the definition of the cylinder valued
momentum map with this holonomy bundle we obtain that
\[
\mathbf{K}(g, \eta)= \pi_C\left (\mbox{\rm Ad}^\ast _{g ^{-1}}\eta+
\mu_{\Sigma}(\widehat{g} ^{-1})+ \nu _0\right)= \Xi(\widehat{g}, \eta+
\overline{{\cal H}})+ \pi _C(\nu _0).
\]   In order to prove~(\ref{non equivariance cylinder valued momentum map
magnetic}) recall that by Proposition~\ref{cocycles for cylinder valued momentum
map}
\[
\sigma(g)= \mathbf{K}(g \cdot  (e,0))- \mathcal{A} d ^\ast _{g ^{-1}} \mathbf{K}
(e,0)= \pi _C(\mu_{\Sigma}(\widehat{g}^{-1})+ \nu _0)- \pi _C(\mbox{\rm Ad}^\ast
_{g ^{-1}}\nu _0)=\pi _C(\mu_{\Sigma}(\widehat{g}^{-1})+ \nu _0-\mbox{\rm Ad}^\ast
_{g ^{-1}}\nu _0).
\]  Finally, let $ \widehat{h} \in  \ker \pi _G $. By~(\ref{non
equivariance cylinder valued momentum map magnetic}), $0= \sigma (e)=
\pi_C(\mu_{\Sigma} (\widehat{h}^{-1}))$. Consequently
$\mu_{\Sigma}(\widehat{h}^{-1})= - \mbox{\rm Ad} ^\ast _{h ^{-1}}\mu_{\Sigma}
(\widehat{h})\in  \overline{{\cal H}} $  and hence, by~(\ref{the
Hamiltonian holonomy is Ad invariant inclusion}), $\mu_{\Sigma}
(\widehat{h})\in 
\overline{{\cal H}} $.

\begin{remark}
\normalfont If there exists a group one-cocycle that integrates $\Sigma$, that is,
the equality~(\ref{integration condition sigma}) holds then we can use the affine
$G$-action $\overline{\Theta}:G \times  \mathfrak{g}^\ast \rightarrow
\mathfrak{g}^\ast$ introduced in Remark~\ref{integration and affine actions} in
order to express the holonomy bundle $\widetilde{(G \times 
\mathfrak{g}^\ast)}_{(h, \mu , \nu)}$ as the graph
\[
\widetilde{(G \times  \mathfrak{g}^\ast)}_{(h, \mu , \nu)}=\{ (g, \eta,
\overline{\Theta}(g, \eta)-\overline{\Theta}(h, \mu)+ \nu)\mid g
\in  G, \eta \in  \mathfrak{g}^\ast
\}.
\] In this case the Hamiltonian holonomy $ {\mathcal H}  $ is obviously trivial
and the cylinder valued momentum map is the standard momentum map given by 
\[
\mathbf{K}(g, \eta)=\overline{\Theta}(g, \eta)+ \nu _0=\Theta(g, \eta)+ \nu _0,
\qquad (g, \eta)\in  G \times  \mathfrak{g}^\ast.
\]
\end{remark}

Once we have computed the cylinder valued momentum map for the symplectic
$G$-action on the magnetic cotangent bundle $(G \times
\mathfrak{g}^\ast, \omega _\Sigma)$ we will carry out reduction in this
context. According to Theorem~\ref{all reductions cylinder valued
momentum map} the optimal reduced spaces provide the symplectic leaves
of the Poisson reduced spaces. In the next two theorems we will describe
these two reduced spaces in the setup of this section. As we will see,
the optimal reduced spaces can be seen as the $G_{\Sigma}$-orbits in
$\mathfrak{g}^\ast$ of the extended affine action, while the Poisson
reduced spaces are, roughly speaking, the $\overline{{\cal
H}}$-saturation of these orbits.

\begin{theorem}[Optimal reduction of magnetic cotangent bundles]
\label{optimal reduction of magnetic cotangent bundles} The optimal reduced spaces
of the canonical $G$-action on the magnetic cotangent bundle $(G \times 
\mathfrak{g}^\ast, \omega _\Sigma )$ are symplectically diffeomorphic to the
orbits corresponding to the extended affine action
$\overline{\Xi}: G_{\Sigma}\times  \mathfrak{g}^\ast \rightarrow \mathfrak{g}^\ast$
of $G_{\Sigma} $ on $\mathfrak{g}^\ast$ endowed with the symplectic structure that
makes them the symplectic leaves of $(\mathfrak{g}^\ast, \{ \cdot , \cdot \}
^\Sigma_-)$ (see Proposition~\ref{symplectic leaves with Hamiltonian holonomy}).
\end{theorem}  

\begin{remark}
\normalfont When the magnetic term is set to zero then $G_{\Sigma}=G $ 
(suppose $G$ is connected) and the extended affine orbits become the
$G$-coadjoint orbits. Consequently, this theorem shows that optimal
reduction, and not the other reduction schemes presented in this paper,
generalizes to the magnetic setup the well known result that says that
the Marsden-Weinstein reduced spaces of the lifted action of a connected
Lie group $G$ on its cotangent bundle (endowed with the canonical
symplectic form) are symplectomorphic to the $G$-coadjoint orbits in
$\mathfrak{g}^\ast$.
\end{remark} 

\noindent\textbf{Proof.\ \ }In order to compute the polar distribution $A_{G}'$ of
the $G$-action on $(G \times  \mathfrak{g}^\ast, \omega _\Sigma )$ notice that for
any $f \in  C^{\infty}(G \times  \mathfrak{g}^\ast)^G\simeq
C^{\infty}(\mathfrak{g}^\ast)$, the corresponding Hamiltonian vector 
field $X _f $ is 
\[ X _f(g, \nu)= \left( T _e L _g \left( \frac{\delta f}{\delta \nu}\right),
\mbox{\rm ad} ^\ast _{\delta f/ \delta \nu }\nu- \Sigma\left(\frac{\delta f}{\delta
\nu}, \cdot \right)\right)
\] 
and hence by~(\ref{infinitesimal generator extended}) we can write
\[ 
A_{G}'(g, \nu)=\{(T _e L _g (\xi), \mbox{\rm ad}^\ast  _\xi \nu-
\Sigma(\xi, \cdot ))\mid \xi \in \mathfrak{g} \} 
=\{(T _e L _g(\pi_{\mathfrak{g}}(\eta)),
-\eta_{\mathfrak{g}^\ast}(\nu))\mid \eta \in  \mathfrak{g}_\Sigma\}.
\] 
Consequently, the leaves of $A_{G}'$ are given by the orbits of the
right $G_{\Sigma}$-action $\Upsilon: G_{\Sigma} \times (G \times 
\mathfrak{g}^\ast)
\rightarrow G \times  \mathfrak{g}^\ast $ defined by $\Upsilon(g,(h, \nu)):=(h
\pi_G (g), \overline{\Xi}(g ^{-1}, \nu))$, $(h, \nu) \in G \times 
\mathfrak{g}^\ast $, $g \in  G_{\Sigma}$. The momentum space $(G \times 
\mathfrak{g}^\ast)/A_{G}' $ can be identified with the orbit space $(G
\times \mathfrak{g}^\ast)/ G_{\Sigma}$ and hence for any $\rho \in  (G
\times  \mathfrak{g}^\ast)/A_{G}'$ there exists an element $\mu\in 
\mathfrak{g}^\ast$ such that  $\mathcal{J}^{-1}(\rho)= G_{\Sigma} \cdot
(e, \mu)$. Moreover
\begin{align*} 
G_{\rho}&=\{g \in  G \mid g \cdot (e, \mu)= \Upsilon(h, (e,
\mu))\text{ for some } h \in G_{\Sigma}\}\\
	&=\{g \in  G\mid (g, \mu)=(\pi _G (h), \overline{\Xi}(h ^{-1},
\mu))\text{ for some } h \in G_{\Sigma}\},
\end{align*} which guarantees that
\begin{equation}
\label{form of isotropy subgroup g rho} G_{\rho}= \pi _G((G_{\Sigma})_\mu).
\end{equation} Consider now the smooth surjective map
\[
\begin{array}{cccc}
\overline{\phi}: & \mathcal{J}^{-1}(\rho)= G_{\Sigma} \cdot (e,
\mu)&\longrightarrow &G_{\Sigma}\cdot  \mu\\
		&(\pi _G(g), \overline{\Xi}(g ^{-1} , \mu))&\longmapsto &\overline{\Xi}(g ^{-1} ,
\mu).
\end{array}
\] The map $\overline{\phi} $ is clearly $G_{\rho} $-invariant and hence it drops
to a smooth surjective map
\[
\phi: \mathcal{J}^{-1}(\rho)/ G_{\rho} \rightarrow G_{\Sigma} \cdot \mu.
\] The map $\phi $ is injective because if $ \phi(\pi _\rho(\pi _G(g),
\overline{\Xi}(g ^{-1} ,
\mu)))= \phi(\pi _\rho(\pi _G(h), \overline{\Xi}(h ^{-1} , \mu)))$ then
$\overline{\Xi}(h g ^{-1}, \mu )= \mu $ and hence $\pi_G (h) \pi_G (g ^{-1})\in 
G_{\rho} $ by~(\ref{form of isotropy subgroup g rho}). Therefore  $\pi_G (h) \pi_G
(g ^{-1}) \cdot (\pi _G(g), \overline{\Xi}(g^{-1} ,\mu))=(\pi _G(h),
\overline{\Xi}(g^{-1} ,\mu))=(\pi _G(h), \overline{\Xi}(h^{-1} ,\mu))$ and hence 
$\pi _\rho(\pi _G(g), \overline{\Xi}(g ^{-1} ,\mu))=\pi _\rho(\pi _G(h),
\overline{\Xi}(h ^{-1} , \mu))$.

The map $\phi $ is a diffeomorphism because for any $\overline{\Xi}(g ^{-1},
\mu)\in  G_{\Sigma} \cdot  \mu $ there is a smooth local section
$\sigma:U_{\overline{\Xi}(g ^{-1},
\mu)}\subset G_{\Sigma}\cdot  \mu \simeq G_{\Sigma}/(G_{\Sigma})_\mu
 \rightarrow G_{\Sigma} $ such that  the map
\[
\begin{array}{cccc}
\phi^{-1}_{U_{\overline{\Xi}(g ^{-1},
\mu)}}: &U_{\overline{\Xi}(g ^{-1},
\mu)}& \longrightarrow &\mathcal{J}^{-1}(\rho)\\
	&\overline{\Xi}(h ^{-1},\mu)&\longmapsto &\left(\pi _G (\sigma(\overline{\Xi}(h
^{-1},\mu))^{-1}) , \overline{\Xi}(h ^{-1},\mu)\right)
\end{array}
\] is a local smooth inverse of $\phi$. Finally, the symplectic character of
$\phi$ is proved via a straightforward diagram chasing exercise. \quad
$\blacksquare$

\medskip

We are now going to provide a model for the Poisson reduced space $M ^{[\mu]}:=
\mathbf{K}^{-1}([\mu])/G _{[\mu]} $ similar to the one provided in the previous
theorem for the optimal reduced space and in which it will be very easy to see how,
as we proved in general in Theorem~\ref{all reductions cylinder valued momentum
map}, the optimal reduced spaces $M _\rho  $ are the symplectic leaves of the
Poisson reduced spaces $M ^{[\mu]} $.

For the sake of simplicity, we will take in our computations the cylinder valued
momentum map in Theorem~\ref{explicit expression of cylinder valued momentum map
magnetic} for which $\nu _0 =0 $, that is,
\begin{equation}
\label{cylinder valued momentum map for computation of reduction}
\mathbf{K}(g, \eta)= \pi_C\left (\mbox{\rm Ad}^\ast _{g ^{-1}}\eta+
\mu_{\Sigma}(\widehat{g} ^{-1})\right)= \Xi(\widehat{g}, \eta+
\overline{{\cal H}}), \qquad (g, \eta) \in  G \times  \mathfrak{g}^\ast,\, \pi_G
(\widehat{g})=g.
\end{equation} In this situation, the non-equivariance cocycle~(\ref{non
equivariance cylinder valued momentum map magnetic}) induces an affine $G$-action
on $\mathfrak{g}^\ast/
\overline{{\cal H}} $ (see Proposition~\ref{cocycles for cylinder valued momentum
map}) given by 
\begin{equation}
\label{affine action magnetic context}
\Theta(g, \pi _C(\mu))= \pi _C(\mbox{\rm Ad}^\ast _{g ^{-1}}\mu+
\mu_{\Sigma}(\widehat{g}^{-1}))= \Xi(\widehat{g}, \pi_C(\mu)),
\end{equation} where $\widehat{g}
\in G_{\Sigma} $ is such that
$\pi_G(\widehat{g})=g$ and hence $\mathbf{K}(g, \eta)=\Theta(g, \pi _C(\eta))$,
$(g,
\eta) \in  G \times  \mathfrak{g}^\ast$. 

It is easy to show that in this specific situation, the Lie algebra two-cocycle
$\widetilde{\Sigma} $ in~(\ref{cocycle of lifted algebra action}) coincides with
$\Sigma$ and hence the Poisson structure on $\mathfrak{g}^\ast/ \overline{{\cal
H}}  $ introduced in Theorem~(\ref{The Poisson structure of target}) with respect
to which $\mathbf{K}:G \times
\mathfrak{g}^\ast
\rightarrow 
\mathfrak{g}^\ast/ \overline{{\cal H}} $ is a Poisson map is
$\left(\mathfrak{g}^\ast/
\overline{{\cal H}}, \{ \cdot , \cdot \}_{\mathfrak{g}^\ast/
\overline{{\cal H}}}^{+ \Sigma}\right) $. 

In the next theorem we will show that the Poisson reduced space $M ^{[\mu]} $ is
Poisson diffeomorphic to $G_{\Sigma}\cdot  \mu+ \overline{{\cal H}}  $ with the
smooth and Poisson structures that we  now discuss.

\paragraph{The smooth structure of $G_{\Sigma}\cdot  \mu+ \overline{{\cal H}}  $.}
Consider the orbit $G_{\Sigma} \cdot  [\mu]:= \Xi(G_{\Sigma}, [\mu])$. By general
theory, $G_{\Sigma} \cdot  [\mu] $ is an initial submanifold of
$\mathfrak{g}^\ast/\overline{{\cal H}}$. Since $\pi_C: \mathfrak{g}^\ast
\rightarrow  \mathfrak{g}^\ast/\overline{{\cal H}} $ is a surjective
submersion, the inverse image $\pi_C ^{-1}(G_{\Sigma} \cdot  [\mu]
)=G_{\Sigma}\cdot  \mu+ \overline{{\cal H}}  $ is an initial submanifold
of  $\mathfrak{g}^\ast $ by the transversality theorem for initial
submanifolds. This is the smooth structure for $G_{\Sigma}\cdot  \mu+
\overline{{\cal H}}  $ that we will consider in what follows.

\paragraph{The Poisson structure of $G_{\Sigma}\cdot  \mu+ \overline{{\cal H}}  $.}
Consider the Poisson manifold $\left(\mathfrak{g}^\ast, \{ \cdot , \cdot
\} ^\Sigma _-\right)$. By Proposition~\ref{symplectic leaves with
Hamiltonian holonomy} the symplectic leaves of this Poisson manifold are
the
$G_{\Sigma} $-orbits of the
$\overline{\Xi} $-action on $\mathfrak{g}^\ast$. Consequently  $G_{\Sigma}\cdot 
\mu+ \overline{{\cal H}}  $ is automatically a quasi-Poisson submanifold of
$(\mathfrak{g}^\ast, \{ \cdot , \cdot \} ^\Sigma _-)$ (see Section 4.1.21
of~\cite{hsr})  and hence a Poisson manifold on its own (see Proposition 4.1.23
of~\cite{hsr}) with the bracket $\{ \cdot , \cdot
\}^{G_{\Sigma}\cdot \mu+ \overline{{\cal H}}}$ given by
\[
\{ f , g
\}^{G_{\Sigma}\cdot \mu+ \overline{{\cal H}}}(\nu)=\{F,G\}_- ^\Sigma (\nu), 
\] where $f,g \in  C^{\infty}(G_{\Sigma}\cdot \mu) $ and $F,G \in
C^{\infty}(\mathfrak{g}^\ast)$  are local extensions of  $f$  and $g$
around $\nu \in G_{\Sigma}\cdot \mu+ \overline{{\cal H}} $,
respectively.

\begin{theorem}[Poisson reduction of magnetic cotangent bundles]
\label{Poisson reduction of magnetic cotangent bundles} Let $(G \times 
\mathfrak{g}^\ast , \omega _\Sigma)$  be the magnetic cotangent bundle of the Lie
group $G$ associated to the Lie algebra two-cocycle $\Sigma$.  Let $\mathbf{K}: G
\times \mathfrak{g}^\ast \rightarrow  \mathfrak{g}^\ast/
\overline{{\cal H}}$  be the cylinder valued momentum map in
\eqref{cylinder valued momentum map for computation of reduction} for
the canonical lifted action of $G$ on  $(G \times  \mathfrak{g}^\ast ,
\omega _\Sigma)$. Then, for any $\mu \in 
\mathfrak{g}^\ast$, the associated Poisson reduced space $M ^{[\mu]}:= \mathbf{K}
^{-1}([\mu])/G _{[\mu]} $ is naturally Poisson diffeomorphic to 
$\left(G_{\Sigma}\cdot \mu+ \overline{{\cal H}},\{ \cdot , \cdot
\}^{G_{\Sigma}\cdot \mu+ \overline{{\cal H}}}\right)$.
\end{theorem}

\begin{remark}
\normalfont In view of Corollary~\ref{inclusion with holonomy relevant abelian}
notice that if
$G$ is Abelian and $\mu_{\Sigma} (G_{\Sigma})$ is closed in $\mathfrak{g}^\ast$
then $\mu_{\Sigma}  $ is a group homomorphism by Proposition~\ref{properties of mu
and more} {\bf (iii)} and hence
\[ G_{\Sigma} \cdot  \mu + \overline{{\cal H}}= \mu + \mu_{\Sigma} (G_{\Sigma})+
\overline{{\cal H}}= \mu + \mu_{\Sigma} (G_{\Sigma})= G_{\Sigma} \cdot  \mu.
\]  
Consequently, by the theorems~\ref{optimal reduction of magnetic
cotangent bundles} and~\ref{Poisson reduction of magnetic cotangent
bundles}, the optimal and Poisson reduced spaces coincide in this case.
\end{remark}

\noindent\textbf{Proof.\ \ }Since $\mathbf{K}(g, \eta)= 
\pi_C\left (\mbox{\rm Ad}^\ast _{g ^{-1}}\eta+
\mu_{\Sigma}(\widehat{g} ^{-1})\right)= \Xi(\widehat{g}, \eta+
\overline{{\cal H}})$, for any $(g, \eta) \in  G \times  \mathfrak{g}^\ast$ and
$\pi_G (\widehat{g})=g$, we have that
\begin{align*}
\mathbf{K}^{-1}([\mu]) &=\{(g, \eta)\in  G \times  \mathfrak{g}^\ast \mid \Xi(
\widehat{g}, \pi _C(\eta))= \pi_C(\mu)\text{ for any }\widehat{g} \in 
G_{\Sigma}\text{ such that }\pi _G(\widehat{g})=g\}\\
	&=\{ (\pi _G(\widehat{g}), \overline{\Xi}(\widehat{g}^{-1}, \mu)+ \nu)\mid
\widehat{g} \in  G_{\Sigma}, \nu \in \overline{{\cal H}}\}.
\end{align*} Consider the smooth surjective map
\[
\begin{array}{cccc}
\overline{\phi} :& \mathbf{K} ^{-1}( [\mu])&\longrightarrow & 
G_{\Sigma} \cdot  \mu + \overline{{\cal H}}\\
&(\pi _G(\widehat{g}), \overline{\Xi}(\widehat{g}^{-1}, \mu)+
\nu)&\longmapsto &
\overline{\Xi}(\widehat{g}^{-1}, \mu)+ \nu.
\end{array}
\] 
The map $\overline{\phi} $ is clearly $G _{[\mu]} $-invariant and hence
it drops to a smooth surjective map
\[
\phi : \mathbf{K} ^{-1}( [\mu])/G _{[\mu]}\longrightarrow G_{\Sigma} 
\cdot  \mu + \overline{{\cal H}}.
\] 
We now show that $\phi  $ is injective. Let  $\pi  ^{[\mu]} \left( \pi
_G(\widehat{g}), \overline{\Xi}(\widehat{g}^{-1}, \mu)+ \nu\right)$ and  
$\pi^{[\mu]} \left( \pi _G(\widehat{h}),
\overline{\Xi}(\widehat{h}^{-1}, \mu)+
\nu' \right) $ be two points in $M ^{[\mu]}$ such that 
\begin{equation}
\label{for injectivity 21}
\overline{\Xi}(\widehat{g}^{-1}, \mu)+
\nu=\overline{\Xi}(\widehat{h}^{-1}, \mu)+ \nu'.
\end{equation} 
Applying $\pi _C $ to both sides of this equality we obtain
$\Xi(\widehat{g}^{-1}, \pi _C(\mu))= \Xi(\widehat{h}^{-1}, \pi _C(\mu))$.
This implies that $\widehat{k}:= \widehat{h} \widehat{g}^{-1} \in 
(G_{\Sigma})_{\pi_C(\mu)}$ and hence by~(\ref{affine action magnetic
context}) 
\[
\Theta(\pi _G (\widehat{k}), \pi _C (\mu))= \Xi(\widehat{k}, \pi _C(\mu))
= \pi_C(\mu)
\] 
which guarantees that  $\pi_G(\widehat{h} \widehat{g}^{-1})\in 
G_{\pi_C(\mu)}$. Therefore, by~(\ref{for injectivity 21}),
$\pi_G(\widehat{h} \widehat{g}^{-1}) \cdot \left( \pi _G(\widehat{g}),
\overline{\Xi}(\widehat{g}^{-1}, \mu)+ \nu\right) 
= \left( \pi_G(\widehat{h}), \overline{\Xi}(\widehat{h}^{-1}, \mu)+
\nu' \right) $ and hence $\pi^{[\mu]} \left( \pi _G(\widehat{g}),
\overline{\Xi}(\widehat{g}^{-1}, \mu)+ \nu\right) =\pi  ^{[\mu]}
\left( \pi _G(\widehat{h}), \overline{\Xi}(\widehat{h}^{-1}, \mu)
+ \nu' \right) $, which shows that $\phi $ is injective. A standard
argument using local sections of $\pi_C $ and of $G_{\Sigma} \rightarrow
G_{\Sigma}/(G_{\Sigma})_{[\mu]} $ (see proof of Theorem~\ref{optimal
reduction of magnetic cotangent bundles}) shows that
$\phi  $ has smooth local inverses and it is hence a diffeomorphism. The Poisson
character of $\phi $ is a straightforward verification. \quad $\blacksquare$

\begin{example}
\label{torus example with open holonomy}
\normalfont
\noindent {\bf An explicit example.} In order to illustrate how to explicitly implement
the constructions introduced in this section we will carry them out for the cotangent
bundle of a four-torus whose canonical symplectic structure has been modified with the
invariant magnetic term induced by the Lie algebra two-cocycle $\Sigma: \mathbb{R}
^4\times  \mathbb{R} ^4\rightarrow \mathbb{R}$ given by the matrix
\[
\Sigma= \left(
\begin{array}{cccc}
0 & 0 &-1 &-\sqrt{2}\\
0 &0 &0 &0\\
1 &0 &0 &0\\
\sqrt{2} &0 &0 &0
\end{array}
\right).
\]
The entries in this matrix have been chosen in such a way that we obtain a Hamiltonian
holonomy group that is not closed in the dual of the Lie algebra.  We will compute a
cylinder valued momentum map for the canonical $\mathbb{T} ^4 $-action on $(T^*
\mathbb{T} ^4, \omega _\Sigma)$ as well as the associated reduced spaces that, as we will
see, are all identical despite the non-closedness of the Hamiltonian holonomy in $
\mathbb{R} ^4 $.

We write $G:= \mathbb{T} ^4 $ and denote by $\mathfrak{g}= \mathbb{R} ^4  $ its Lie
algebra. We start by noting that the one-dimensional central extension
$\mathfrak{g}_\Sigma= \mathbb{R} ^4\oplus \mathbb{R}$ of  $\mathfrak{g}$ is integrated by
the Heisenberg group $G_{\Sigma}= \mathbb{R} ^4\oplus  \mathbb{R} $ with 
multiplication given by 
\[
({\bf u},a)\cdot  ({\bf v},b):=
\left({\bf u}+{\bf v}, a+b- \frac{1}{2} \Sigma ({\bf u},{\bf v})\right)
=\left({\bf u}+{\bf v}, a+b- u _1v _3-\sqrt{2} u _1 v  _4+u _3 v 
_1+\sqrt{2} u _4 v  _2\right),
\]
with $({\bf u},a), ({\bf v},b)\in  G_{\Sigma}  $. An easy calculation shows that for any
$({\bf u},a)\in  G_{\Sigma}  $ and  $(\xi,s )\in  \mathfrak{g}  _\Sigma$,
\[
\mbox{\rm Ad}_{(u,a)}(\xi,s)=( \xi, s - \Sigma(u, \xi)).
\]
Consequently, in view of Proposition~\ref{properties of mu and more}, the extended
$\mathfrak{g}^\ast$-valued one-cocycle $\mu_{\Sigma}: G_{\Sigma} \rightarrow 
\mathfrak{g}^\ast$ associated to $\Sigma$ is given by $ \mu_{\Sigma}(u,a)= \Sigma(\cdot
,u)$, $({\bf u},a)\in  G_{\Sigma}  $. Using Proposition~\ref{intermediated connection
magnetic} {\bf (iii)} and  Proposition~\ref{characterization with extension holonomy}  
it is easy to see that the Hamiltonian holonomy ${\mathcal H} $ of our
setup is given by 
\[
{\mathcal H}= (\Bbb Z +\sqrt{2}\Bbb Z) \times \{0\} \times  \Bbb Z(1, \sqrt{2}),
\]
which is clearly not closed in $\mathbb{R} ^4 $ since $\overline{{\mathcal H}}=
\mathbb{R} \times  \{0\}\times  \Bbb Z(1, \sqrt{2}) $. In order to write down a cylinder
valued momentum map for our example we start by noting that the map 
$\mathbb{R} ^4/ \overline{{\cal H}} \rightarrow \mathbb{R} ^2\times  S ^1$ given by
$[a,b,c,d]\longmapsto (b, \frac{1}{3}(d -\sqrt{2}c),e ^{\frac{2\pi i}{3}(c+\sqrt{2}d )})$
is a group isomorphism and hence we can write $\pi _C(a,b,c,d)=\left(b,
\frac{1}{3}(d -\sqrt{2}c),e ^{\frac{2\pi i}{3}(c+\sqrt{2}d )}\right)$.
Consequently, since the map $\pi _G : G_{\Sigma} \rightarrow  G $ is
given by $({\bf u},a) \mapsto \left(e^{2\pi i u _1},\ldots, e^{2\pi i
u_4}\right)$, we have by Theorem~\ref{explicit expression of cylinder
valued momentum map magnetic} that the map $\mathbf{K}: \mathbb{T}^4
\times 
\mathbb{R} ^4 
\rightarrow  \mathbb{R} ^2\times  S ^1$ defined by
$\mathbf{K}(g,\eta):= \pi _C (\eta + \mu _\Sigma (-u))= \pi _C(\eta - \Sigma(\cdot , u))$
is a cylinder valued momentum map for the $\mathbb{T} ^4 $-action on $(T^*
\mathbb{T} ^4, \omega _\Sigma)$; in this definition, the elements $g \in  \mathbb{T} ^4 $
and $u \in  \mathbb{R}^4 $ are related by the equality $g=\left(e^{2\pi
i u _1},\ldots, e^{2\pi i u _4}\right)$. We can be even more specific by
writing
\[
\mathbf{K}\left(\left(e^{2\pi i u _1},\ldots,e^{2\pi i u _4}\right),
\eta\right)
=  \left(\eta_2, \frac{1}{3} \left(\eta_4-\sqrt{2} u _1-\sqrt{2}
(\eta_3-u _1)\right),e^{\frac{2\pi i}{3}\left(\eta _3-u _1+\sqrt{2}
(\eta_4-\sqrt{2}u _1) \right)}\right).
\]
We now compute the reduced spaces. We start by noticing that in this particular case
the Lie algebra $\mathfrak{n} := \left({\rm Lie}\left(\overline{{\cal
H}}\right)\right)^\circ =\{0\} \times  \mathbb{R} ^3  $. The subgroup
$N=\{e\} \times  \mathbb{T}^3
\subset \mathbb{T} ^4  $ clearly has $\mathfrak{n} $ as Lie algebra. Let  now  $[\mu]\in
\mathbb{R} ^4 / \overline{{\cal H}}\simeq \mathbb{R}^2 \times  S ^1$ be
arbitrary. We have 
\[
N _{ [\mu]}=\{n \in  N\mid \Theta(n, [\mu])= [\mu]
\}=\{n \in  N\mid  \pi _C(\mu_{\Sigma} (\widehat{n}))=0, 
\pi _G(\widehat{n})=n\}=N.
\]
A similar computation shows that $G _{[\mu]}=N= N _{ [\mu]} $. Consequently, the Poisson
$M ^{[\mu]} $   and symplectic $M _{[\mu]} $ reduced spaces coincide. Moreover, by
Theorem~\ref{Poisson reduction of magnetic cotangent bundles} they are naturally Poisson
diffeomorphic to $G_{\Sigma} \cdot  \mu + \overline{{\cal H}} = \mu  + \mu_{\Sigma}
(G_{\Sigma})+ \overline{{\cal H}}$. It can be checked that, in this particular case,
$\overline{{\cal H}}\subset \mu_{\Sigma} (G_{\Sigma})$ and hence, by Theorem~\ref{optimal
reduction of magnetic cotangent bundles}, 
\[
 \mu  + \mu_{\Sigma}(G_{\Sigma})+ \overline{{\cal H}}=\mu  +
\mu_{\Sigma}(G_{\Sigma})= G_{\Sigma} \cdot  \mu \simeq M _\rho. 
\]  
Consequently,
\[
M _{[\mu]}=M ^{[\mu]}\simeq M _\rho \simeq \mu  +
\mu_{\Sigma}(G_{\Sigma})= \mu + \left( \mathbb{R} \times  \{0\} \times  \mathbb{R}
(1,\sqrt{2})\right).
\]
\end{example}

\section{Appendix: The relation between Lie group and cylinder valued momentum
maps}

The cylinder valued momentum maps are closely related to the {\bfi  Lie 
group valued momentum maps}  introduced in~\cite{mcduff 1988,
Ginzburg1992, Huebschmann and Jeffrey 1994, huebschmann 1995 a, lie
group valued maps}. We give the definition of  these objects only for
Abelian symmetry groups because in the non-Abelian case these momentum
maps are defined on spaces that are neither symplectic nor Poisson (they
are referred to as {\bfi  quasi-Hamiltonian spaces}~\cite{lie group
valued maps}).

\begin{definition}
\label{definition of lie group valued momentum map} 
Let $G$ be an Abelian Lie group whose Lie algebra $\mathfrak{g}$ acts
canonically on the symplectic manifold $(M, \omega)$. Let $(\cdot , \cdot
)$ be some bilinear symmetric   nondegenerate form  on the Lie algebra
$\mathfrak{g}$. The map $\mathbf{J}:M \rightarrow G$ is called a
$G$-{\bfi  valued momentum map} for the $\mathfrak{g}$-action on $M$
whenever 
\begin{equation}
\label{defining relation of the Lie group valued momentum map}
\mathbf{i}_{\xi_M} \omega(m) (v_m) 
=\left(T _m (L _{\mathbf{J}(m)^{-1}}\circ \mathbf{J})( v _m),
\xi\right),
\end{equation}  
for any $\xi\in \mathfrak{g}$, $m \in M$, and $v _m \in T _m M$.
\end{definition}

\begin{proposition} Let $G$ be an Abelian Lie group whose Lie algebra
$\mathfrak{g}$ acts canonically on the symplectic manifold $(M, \omega)$. Let 
$\mathbf{J}:M \rightarrow G $ be a $G$-valued momentum map for this
action.
\begin{itemize}
\item [{\rm {\bf (i)}}] The fibers of $\mathbf{J}:M \rightarrow G$
are invariant under the Hamiltonian flows corresponding to
$\mathfrak{g}$-invariant Hamiltonian functions.
\item [{\rm {\bf (ii)}}] $\ker T _m \mathbf{J}= (\mathfrak{g}\cdot
m)^\omega$ for any $m \in M$.
\end{itemize}
\end{proposition}

\noindent\textbf{Proof.\ \ } {\bf (i)} Let $F _t $ be the flow of the Hamiltonian
vector field $X _h$ associated to a
$\mathfrak{g}$-invariant function $h \in C^\infty(M)^{\mathfrak{g}}$. By the
defining relation~(\ref{defining relation of the Lie group valued momentum map})
of the Lie group valued momentum maps we have  for any $m \in M$ and any $\xi\in
\mathfrak{g}$
\begin{align*} ((T_{\mathbf{J}(F _t (m))}L _{\mathbf{J}(F _t (m)) ^{-1}}&\circ T
_{F _t(m)}
\mathbf{J})(X _h(F _t(m))), \xi )=(T_{F _t(m)}(L _{\mathbf{J}(F _t (m)) ^{-1}}\circ
\mathbf{J}))( X _h(F _t(m))), \xi)\\
	&= \mathbf{i}_{\xi_M} \omega(F _t(m)) )( X _h(F _t(m)))=- \mathbf{d}h (F _t(m)))( 
\xi_M(F _t(m)))=0.
\end{align*}  Consequently, 
\[(T_{\mathbf{J}(F _t (m))}L _{\mathbf{J}(F _t (m)) ^{-1}}\circ T _{F _t(m)}
\mathbf{J})(X _h(F _t(m)))=0
\]
and hence $T _{F _t(m)}\mathbf{J} (X _h(F _t(m)))=0  $, which can be rewritten
as 
\[
\frac{d}{dt} (\mathbf{J} \circ F _t) (m)= 0.
\] The arbitrary character of $t$ and  $m$ implies that $\mathbf{J} \circ F _t=
\mathbf{J}|_{{\rm Dom}(F _t)}$, as required. 

\smallskip

\noindent {\bf (ii)} A vector $v _m \in \ker T _m \mathbf{J}$ if and only if  
$T_m\mathbf{J} ( v _m)= 0 $. This identity is equivalent to $((T_{\mathbf{J}(m)}L
_{\mathbf{J}(m) ^{-1}}\circ T _{ m } \mathbf{J})(v _m), \xi) =0$, for any $\xi \in
\mathfrak{g}$ and, by~(\ref{defining relation of the Lie group valued momentum
map}), to
$\mathbf{i}_{\xi_M} \omega( m ) (v _m) =0$, for all $\xi\in \mathfrak{g}$, which in
turn amounts to $v _m \in (\mathfrak{g}\cdot m)^\omega$. \quad $\blacksquare$

\paragraph{Lie group and cylinder valued momentum maps.} We start with a
proposition that states that any cylinder valued momentum map associated to an
Abelian Lie algebra action whose corresponding holonomy group is closed  can be
understood as a Lie group valued momentum map.

\begin{proposition}
\label{first indication of the relation cylinder valued momentum map and lie} Let
$(M, \omega)$ be a connected paracompact symplectic manifold and $\mathfrak{g}$
 an Abelian Lie algebra acting canonically on it. Let ${\cal H}\subset
\mathfrak{g}^\ast$ be the holonomy group associated to the connection
$\alpha$ in~\eqref{definition of alpha connection} and $(\cdot , \cdot ):
\mathfrak{g}\times \mathfrak{g}\rightarrow \mathbb{R}$  some bilinear symmetric 
nondegenerate form  on 
$\mathfrak{g}$. Let $f : \mathfrak{g} \rightarrow \mathfrak{g}^\ast $ be the
isomorphism given by $ \xi\longmapsto (\xi, \cdot )$, $\xi\in \mathfrak{g}$ and
${\cal T}:= f ^{-1}({\mathcal H})$. The map $f$ induces an Abelian group
isomorphism $\bar{f}:
\mathfrak{g}/ {\cal T} \rightarrow \mathfrak{g}^\ast/ {\mathcal H} $ by
$\bar{f}(\xi+ {\cal T}):=(\xi, \cdot )+ {\mathcal H} $. Suppose that ${\mathcal H}
$ is closed in
$\mathfrak{g}^\ast$ and define $\mathbf{J}:=\bar{f}^{-1} \circ \mathbf{K}: M
\rightarrow
\mathfrak{g}/ {\cal T} $, where $\mathbf{K}: M \rightarrow \mathfrak{g}^\ast/
{\mathcal H}$ is a cylinder valued momentum map for the $\mathfrak{g}$-action on
$(M, \omega) $. Then
\begin{equation}
\label{cylinder valued momentum map look like lie}
\mathbf{i}_{\xi_M} \omega(m) (v_m) 
= \left(T_m (L_{\mathbf{J}(m)^{-1}}\circ \mathbf{J})(v_m), \xi\right),
\end{equation} 
for any $\xi\in \mathfrak{g}$ and $v _m \in T _m M$. Consequently,
the map $\mathbf{J}: M \rightarrow \mathfrak{g}/ {\cal T} $ constitutes
a $\mathfrak{g}/ {\cal T} $-valued momentum map for the canonical action
of the Lie algebra $\mathfrak{g}$ of
$(\mathfrak{g}/{\cal T}, +)  $ on $(M, \omega)$. 
\end{proposition}

\noindent\textbf{Proof.\ \ } We start by noticing that the right-hand side
of~(\ref{cylinder valued momentum map look like lie}) makes sense due to the
closedness hypothesis on ${\mathcal H} $. Indeed, this condition and the fact that 
${\mathcal H} $ is discrete due to the flatness of $\alpha$ imply that
$\mathfrak{g} ^\ast / {\mathcal H}
$, and therefore
$\mathfrak{g}/ {\cal T} $, are Abelian Lie groups whose Lie algebras can be
naturally identified with $\mathfrak{g}^\ast  $ and $\mathfrak{g}$, respectively.
This identification is used in~(\ref{cylinder valued momentum map look like lie}),
where we think of   $T _m (L _{ \mathbf{J}(m)^{-1}}\circ 
\mathbf{J})(v _m )\in {\rm Lie} (\mathfrak{g}/ {\cal T})$ as an element of
$\mathfrak{g}$.

In what follows we will use the following notation: if
$\mu\in \mathfrak{g}^\ast$ arbitrary, denote by $\xi_\mu \in
\mathfrak{g} $ the unique element such that $\mu=(\xi_\mu, \cdot )$. 

We now compute $T _m \mathbf{J}(v _m)$. Let
$\mu + {\mathcal H}:=  \mathbf{K}(m)$ and hence $ \mathbf{J}(m)= \xi_\mu+ {\cal
T}$. Now, by Theorem~\ref{properties of the cylinder valued momentum map in list} 
{\bf (ii)} we have 
\begin{equation*} T _m\mathbf{J}( v _m) = T _m(\bar{f}^{-1} \circ \mathbf{K})(v
_m) =T_{\mu+ {\mathcal H}}\bar{f}^{-1} (T _m \mathbf{K} (v _m)) =T_{\mu+ {\mathcal
H}}\bar{f}^{-1} ( T _\mu \pi_C (\nu)),
\end{equation*} where the element $\nu\in \mathfrak{g}^\ast $ is given by 
\begin{equation}
\label{definition of nu for cylinder valued momentum map}
\langle \nu, \eta \rangle= {\bf i} _{\eta _M} \omega (m)(v _m), \quad\text{for
all}\quad\eta \in \mathfrak{g}.
\end{equation} Since  $(\bar{f}^{-1} \circ \pi_C)(\rho)= \xi_\rho+ {\cal T} $ for
any $\rho \in
\mathfrak{g}^\ast $, we can write 
\[ T_{\mu+ \mathcal{H}}\bar{f}^{-1} (T _\mu \pi_C (\nu)) 
=T_{\mu}\left(\bar{f}^{-1}\circ \pi_C\right)(\nu) =
\left.\frac{d}{dt}\right|_{t=0}\left(\bar{f}^{-1} \circ \pi_C
\right)(\mu+ t \nu)= \left.\frac{d}{dt}\right|_{t=0} \left(\xi_\mu+ t \xi_\nu+
{\cal T}
\right).
\] Hence, 
\[ T _m\mathbf{J}( v _m) =
\left.\frac{d}{dt}\right|_{t=0} \left(\xi_\mu+ t \xi_\nu+ {\cal T}
\right) \in T_{\xi_\mu+ {\cal T}}\left(\mathfrak{g}/ {\cal T} \right).
\] Now,
\begin{align*}
\left(T _m (L _{ \mathbf{J}(m)^{-1}}\circ 
\mathbf{J})(v _m), \xi\right)  &=\left(T _{\mathbf{J}(m)} L
_{\mathbf{J}(m)^{-1}}(T _m \mathbf{J} (v _m )),
\xi\right)= \left(\left.\frac{d}{dt}\right|_{t=0}(- \xi_\mu+ {\cal T})+(\xi_\mu+ t
\xi_\nu + {\cal T}), \xi \right)\\
	&=(\xi_\nu, \xi)=\langle \nu, \xi\rangle=  {\bf i} _{\xi_M} \omega(m) (v _m),
\end{align*} where the last equality is a consequence of~(\ref{definition of nu
for cylinder valued momentum map}). \quad $\blacksquare$

\paragraph{Lie group valued momentum maps produce closed Hamiltonian holonomies.} So
far we have investigated how cylinder valued momentum maps can be viewed as 
Lie group valued momentum maps. Now we shall focus on the converse
relation, that is, we shall isolate hypotheses that guarantee that a Lie
group valued momentum map naturally induces a cylinder valued momentum
map.

\begin{theorem} Let $(M, \omega)$ be a connected paracompact symplectic manifold
and $\mathfrak{g}$ an Abelian Lie algebra acting canonically on it. Let ${\cal
H}\subset
\mathfrak{g}^\ast$ be the Hamiltonian holonomy group associated to the connection
$\alpha$ in~\eqref{definition of alpha connection} associated to the
$\mathfrak{g}$-action and let $(\cdot ,\cdot ):
\mathfrak{g}\times \mathfrak{g}\rightarrow \mathbb{R}$ be a bilinear symmetric 
nondegenerate form  on $\mathfrak{g}$. Let $f : \mathfrak{g} \rightarrow
\mathfrak{g}^\ast $, $\bar{f}:
\mathfrak{g}/ {\cal T} \rightarrow \mathfrak{g}^\ast/ {\mathcal H} $, and let
${\cal T}:= f ^{-1}({\mathcal H})$ be as in the statement of Proposition~{\rm
\ref{first indication of the relation cylinder valued momentum map and lie}}. Let
$G $  be a connected Abelian Lie group whose Lie algebra is $\mathfrak{g} $ and
suppose that there exists a
$G$-valued momentum map ${\bf A}:M \rightarrow G $ associated to the
$\mathfrak{g}$-action whose definition uses the form $(\cdot , \cdot )$. 
\begin{itemize}
\item [{\rm {\bf (i)}}] If $\exp : \mathfrak{g}\rightarrow G $ is the
exponential map, then 
\begin{equation}
\label{inclusion of holonomy in kernel} {\mathcal H}
\subset f (\ker \exp ).
\end{equation}
\item [{\rm {\bf (ii)}}] ${\mathcal H} $ is closed in
$\mathfrak{g}^\ast$.
\end{itemize} 
Let  $\mathbf{J}:=\bar{f}^{-1} \circ
\mathbf{K}: M \rightarrow
\mathfrak{g}/ {\cal T} $, where $\mathbf{K}: M \rightarrow \mathfrak{g}^\ast/
{\mathcal H}$ is a cylinder valued momentum map for the $\mathfrak{g}$-action on
$(M, \omega) $. If
$f (\ker \exp )\subset {\mathcal H} $, then $\mathbf{J}:M \rightarrow
\mathfrak{g}/ {\cal T}=\mathfrak{g}/\ker\exp \simeq G $  is a $G$-valued momentum
map that differs from ${\bf A}$ by a constant in $G $.

Conversely, if ${\mathcal H}=f (\ker\exp) $, then $\mathbf{J}:M \rightarrow
\mathfrak{g}/\ker\exp \simeq G$ is a $G$-valued momentum map.
\end{theorem}

\begin{remark}
\normalfont The presence of a Lie group valued momentum map associated to a
canonical Lie algebra action does not imply the reverse inclusion
in~(\ref{inclusion of holonomy in kernel}). A simple example that illustrates this
statement is the canonical action of a two torus
$\mathbb{T} ^2$ on itself  via
\[ (e^{i \phi_1}, e^{i \phi_2}) \cdot (e^{i \theta_1}, e^{i \theta_2}):=(e^{i
(\theta_1+
\phi_1)}, e^{i
\theta_2}),
\] where we consider the torus as a symplectic manifold with the area form. A
straightforward computation shows that the surface
\[
\widetilde{ \mathbb{T}^2}:=\left\{((e^{i \theta_1}, e^{i
\theta_2}),(\theta_2,0))\in
\mathbb{T}^2\times \mathbb{R}^2\mid
\theta_1,
\theta_2 \in \mathbb{R}\right\},
\] is the holonomy bundle containing the point $((e,e),(0,0))\in
\mathbb{T}^2\times \mathbb{R}^2 $ associated to the connection that defines the
corresponding cylinder valued momentum map. This immediately shows that ${\mathcal
H}=
\mathbb Z\times \{0\} $ while $f (\ker \exp)= \mathbb Z \times \mathbb Z $ which
is clearly not contained in ${\mathcal H}$.
\end{remark}

\noindent\textbf{Proof of the theorem.\ \ } We start by assuming that the
$\mathfrak{g}$-action on
$(M,
\omega)$ has an associated $G$-valued momentum map ${\bf A}:M \rightarrow G $ and
we will show that ${\mathcal H} \subset f (\ker \exp)$.

Let $ \mu \in {\mathcal H} $. The definition of the holonomy group ${\mathcal H} $
implies the existence of a piecewise smooth loop  $m:[0,1] \rightarrow M  $ at the
point $m$, that is, $m (0)= m (1)= m \in M $, whose horizontal lift $\widetilde{m}
(t)=(m (t),
\mu(t))$ starting at the point $(m,0)$ satisfies $\mu= \mu(1) $. The horizontality
of
$\widetilde{m} (t) $ implies that 
\begin{equation*}
\langle \dot{\mu}(t), \xi \rangle = \mathbf{i} _{\xi_M} \omega (m(t))(\dot{m}(t))
=\left(T_{m (t)} \left(L _{{\bf A} (m (t))^{-1}}\circ  {\bf A}\right)(\dot{m}(t)),
\xi \right),
\end{equation*} for any $\xi\in \mathfrak{g}$ or, equivalently,
\begin{equation}
\label{horizontal condition for loop}
\dot{\mu}(t)=f \left(\left.\frac{d}{ds}\right|_{s=0} {\bf A} (m (t))^{-1} {\bf
A}(m (s)) \right).
\end{equation} Fix $t _0 \in [0,1]$. Since the exponential map
$
\exp :
\mathfrak{g} \rightarrow G$ is a local diffeomorphism, there exists a smooth curve
$\xi:I _{t _0}:=(t _0- \epsilon, t _0+ \epsilon) \rightarrow \mathfrak{g} $, for
$\epsilon>0 $ sufficiently small, such that for any $s \in (- \epsilon, \epsilon)$
\begin{equation}
\label{definition for function in algebra} {\bf A}(m(t _0+ s ))= \exp \xi(t _0+ s)
{\bf A}(m (t _0)).
\end{equation} We now reformulate locally the expression~(\ref{horizontal
condition for loop}) using the function $\xi:I _{t _0} \rightarrow \mathfrak{g} $.
Let $\tau, \sigma \in (- \epsilon,
\epsilon)$ be such that $t= t _0+ \tau $ and $s= t _0+ \sigma$.
Expression~(\ref{horizontal condition for loop}) can be rewritten in $I _{t _0} $
as
\begin{align*}
\frac{d }{d \tau} \mu (t _0+ &\tau)=f \left(\left. \frac{d }{d
\sigma}\right|_{\sigma= \tau}{\bf A}(m(t _0+ \tau))^{-1}{\bf A}(m(t _0+
\sigma))\right)\\
	&=f \left(\left. \frac{d }{d
\sigma}\right|_{\sigma= \tau}\exp (- \xi(t _0+ \tau)) \exp \xi(t _0+ \sigma){\bf
A}(m(t _0 ))^{-1}{\bf A}(m(t _0 ))\right)\\
	&=f \left(\left. \frac{d }{d
\sigma}\right|_{\sigma= \tau}\exp \left(\xi(t _0+ \sigma)- \xi(t _0+ \tau) \right)
\right)=f \left(T _0 \exp \left(\left. \frac{d }{d
\sigma}\right|_{\sigma= \tau}  \left(\xi(t _0+ \sigma)- \xi(t _0+ \tau) \right)
\right)\right)\\
	&=f \left(\left. \frac{d }{d
\sigma}\right|_{\sigma= \tau}  \left(\xi(t _0+ \sigma)- \xi(t _0+ \tau) \right)
\right)=f \left(\left. \frac{d }{d
\sigma}\right|_{\sigma= \tau}  \xi(t _0+ \sigma)
\right)=f \left(  \frac{d }{d
\tau}   \xi(t _0+ \tau)
\right),
\end{align*} which shows that for any $ t \in I _{t _0} $
\begin{equation}
\label{differential equation for algebra element}
\dot{\mu}(t)=f (\dot{\xi}(t)).
\end{equation} We now cover the interval $[0,1]$ with a finite number of intervals
$I _1, \dots, I _n $ such that in each of them we define a function $\xi_i:I _i
\rightarrow \mathfrak{g}$ that satisfies~(\ref{definition for function in
algebra}) and~(\ref{differential equation for algebra element}). We now write $I
_i=[a _i, a _{i+1}] $, with $i \in \{1,
\dots, n\} $, $a _1 =0 $, and $a _{n+1 }=1 $. Using these intervals, since $\mu(0)
= 0$, we can write
\begin{align}
\mu&= \mu(1)=\int _0^1\dot{\mu}(t) d t=\int _{I _1}\dot{\mu}(t) d t+ \cdots + \int
_{I _n}\dot{\mu}(t) d t=f \left(\int _{I _1}\dot{\xi}_1(t) d t+ \cdots + \int _{I
_n}\dot{\xi}_n(t) d t
\right)\notag\\
	&=f \left( \xi_1(a _2)- \xi_1(a _1)+ \cdots+\xi_n(a _{n+1})-\xi_n(a
_{n})\right).\label{last equality we will see}
\end{align} The construction of the intervals $I _i$, $i \!\in\! \{1, \dots, n\} $
implies that
 ${\bf A}(m(a _i))\!= \!\exp \xi_i(a _i)\times {\bf A}(m(a _i))$. Hence 
\begin{equation}
\label{first guy in the kernel}
\exp \xi_i(a _i)=e
\end{equation}  and hence $\xi_i(a _i) \in \ker \exp $ for all $i \in \{1, \dots,
n\} $. We also have that
\begin{align*} {\bf A}(m (1)) &={\bf A}(m(a _{n+1}))= \exp \xi_n(a _{n+1}){\bf A}(
m( a _n))=\exp
\xi_n(a _{n+1})\exp \xi_{n-1}(a _{n}){\bf A}( m( a _{n-1}))\\
	&=\exp \xi_n(a _{n+1})\exp \xi_{n-1}(a _{n}) \cdots \exp \xi_1(a _2){\bf A}(m (a
_1))
	= \exp(\xi_1(a _2)+ \cdots+ \xi_n(a _{n+1})){\bf A}(m (0)).
\end{align*} Since $m (0)=m (1)=m $ we have ${\bf A}(m (0))={\bf A}(m (1)) $ and
therefore
$\exp(\xi_1(a _2)+ \cdots+ \xi_n(a _{n+1}))= e $, which implies that $\xi_1(a _2)+
\cdots+ \xi_n(a _{n+1}) \in \ker \exp $. If we substitute this relation
and~(\ref{first guy in the kernel})   in~(\ref{last equality we will see}) we
obtain that $\mu \in f (\ker \exp)$, which proves the inclusion ${\mathcal
H}\subset f(\ker \exp)$.

We now show that ${\mathcal H} $ is closed in $\mathfrak{g}^\ast$. The inclusion
${\mathcal H}\subset f(\ker \exp)$, the closedness of $\ker \exp $ in
$\mathfrak{g}$, and the fact that $f$ is an isomorphism of Lie
groups imply that
\[
\overline{{\mathcal H}}\subset \overline{f(\ker \exp)}=f(\ker \exp).
\] Because $G$ is Abelian, $\ker \exp $ is a discrete subgroup of
$(\mathfrak{g},+)$  and hence $\overline{ {\mathcal H}} $ is a discrete subgroup
of $\mathfrak{g}^\ast $. This implies that
$\overline{{\mathcal H}} \subset {\mathcal H}$. Indeed, for any element $\mu\in 
\overline{{\mathcal H}}$ there exists an open neighborhood $U _\mu \subset
\mathfrak{g}^\ast $ such that $U _\mu\cap \overline{{\mathcal H}} =\{ \mu\} $
($\overline{ {\mathcal H}}$ is discrete). As
$\mu \in \overline{ {\mathcal H}}$ we have that $\varnothing \neq U
_\mu\cap {\mathcal H}
\subset U _\mu\cap \overline{{\mathcal H}}=\{ \mu\} $, which implies that $\mu\in
{\mathcal H} $. This shows that $ {\mathcal H}= \overline{{\mathcal H}} $ and
therefore that ${\mathcal H}$ is closed in $\mathfrak{g}^\ast$.

Assume now that $f (\ker \exp )\subset {\mathcal H} $. The hypothesis on the
existence of a Lie group valued momentum map implies, via the
inclusion~(\ref{inclusion of holonomy in kernel}) that we  just proved,
that $f (\ker \exp )= {\mathcal H}$ and that ${\mathcal H} $ is closed in
$\mathfrak{g}^\ast$. Proposition~\ref{first indication of the relation cylinder
valued momentum map and lie} implies that $\mathbf{J}:M
\rightarrow \mathfrak{g}/\ker \exp\simeq G $ is a $G$-valued momentum map for the
$\mathfrak{g}$-action on $(M, \omega)$. We now show that $\mathbf{A}$ and
$\mathbf{J} $ differ by a constant in $G$. The expression~(\ref{defining relation
of the Lie group valued momentum map}) for $\mathbf{A}$ and~(\ref{cylinder valued
momentum map look like lie}) for $\mathbf{J} $ imply that for any $\xi\in
\mathfrak{g} $ and any $v _m \in T _mM
$ we have  
\[
\left(T _m (L _{ \mathbf{A}(m)^{-1}}\circ 
\mathbf{A}) (v _m), \xi\right)=\mathbf{i}_{\xi_M} \omega(m)( v _m)  =\left(T _m (L
_{\mathbf{J}(m)^{-1}}\circ 
\mathbf{J})(v _m), \xi\right),
\] which implies that $T \mathbf{J}=T \mathbf{A}$. Since the manifold $M$ is
connected we have that ${\bf A}$ and ${\bf J}$ coincide up to a constant element
in $\mathbf{G}$.

The last claim in the theorem is a straightforward corollary of 
Proposition~\ref{first indication of the relation cylinder valued momentum map and
lie}.~\ ~$\blacksquare$  

\medskip

\noindent\textbf{Acknowledgments} This research was partially supported by the
European Commission through funding for the Research Training Network
\emph{Mechanics and Symmetry in Europe} (MASIE) as well as the Swiss National
Science Foundation. Part of this work was carried out during the program
``Geometric Mechanics and Its Applications'' at the Bernoulli Center (CIB) of the
\'Ecole Polytechnique F\'ed\'erale de Lausanne that provided an excellent working
environment. We thank Alan Weinstein for encouraging us to look at the ``moment
r\'eduit'' of~\cite{condevaux dazord and molino} and Claude Roger for helpful
discussions on central extensions.

\end{document}